\documentclass[a4paper]{scrartcl}
\usepackage[utf8]{inputenc}
\usepackage[T1]{fontenc}
\usepackage[english]{babel} 
\usepackage{lmodern}
\usepackage[]{algorithm2e}
\usepackage[onehalfspacing]{setspace}
\usepackage{amsmath,amssymb,amsthm,amsfonts,amsbsy,latexsym}
\usepackage{graphicx}
\usepackage{geometry}
\usepackage{natbib}
\usepackage{multirow}
\usepackage{float}
\usepackage{fancyhdr}
\usepackage{makecell}
\usepackage{pgf,tikz}
\usepackage{makeidx}
\usetikzlibrary{shapes,snakes}
\usetikzlibrary{fadings}
\usepackage[pdfpagelabels=true]{hyperref}

\newcommand{\E}{\mbox{I\negthinspace E}}

\newcommand{\R}{\mathbb{R}}

\numberwithin{equation}{section}
\IfFileExists{upquote.sty}{\usepackage{upquote}}{}
\begin{document}

\newtheorem{thm}{Theorem}[section]
\newtheorem{Def}{Definition}[section]
\newtheorem{lem}{Lemma}[section]
\newtheorem{prop}{Proposition}[section]
\newtheorem{rem}{Remark}[section]
\newtheorem{cor}{Corollary}[section]
\newtheorem{ex}{Example}[section]
\newtheorem{ass}{Assumption}[section]
\newtheorem*{bew}{Proof}

\title{Robust Statistics meets elicitability: When fair model validation breaks down}
\author{Tino Werner\footnote{Institute for Mathematics, Carl von Ossietzky Universität Oldenburg, P/O Box 2503, 26111 Oldenburg (Oldb), Germany, \texttt{tino.werner1@uni-oldenburg.de}}}
\maketitle

\abstract{A crucial part of data analysis is the validation of the resulting estimators, in particular, if several competing estimators need to be compared. Whether an estimator can be objectively validated is not a trivial property. If there exists a loss function such that the theoretical risk is minimized by the quantity of interest, this quantity is called elicitable, allowing estimators for this quantity to be objectively validated and compared by evaluating such a loss function. Elicitability requires assumptions on the underlying distributions, often in the form of regularity conditions. Robust Statistics is a discipline that provides estimators in the presence of contaminated data. In this paper, we, introducing the elicitability breakdown point, formally pin down why the problems that contaminated data cause for estimation spill over to validation, letting elicitability fail. Furthermore, as the goal is usually to estimate the quantity of interest w.r.t. the non-contaminated distribution, even modified notions of elicitability may be doomed to fail. The performance of a trimming procedure that filters out instances from non-ideal distributions, which would be theoretically sound, is illustrated in several numerical experiments. Even in simple settings, elicitability however often fails, indicating the necessity to find validation procedures with non-zero elicitability breakdown point.}

\section{Introduction}

\subsection{Comparing competing models}

In a classical statistical decision problem, let $A$ be the decision space and let $\mathcal{A}$ be the corresponding $\sigma$-algebra. The goal is to find the correct decision based on a sample $y=(y_1,...,y_n) \in \mathcal{Y}^n$, for some space $\mathcal{Y}$, generated from some unknown distribution $F_Y \in \mathcal{F}_Y$ for some distribution family $\mathcal{F}_Y$ on $\mathcal{Y}$. Let $S_n: \mathcal{Y}^n \rightarrow A$ be the decision function with which the statistician makes the decision $S_n(y_1,...,y_n)$ if the sample $y$ is available. Assume that the quality of any decision $a \in A$ can be validated by a loss function. In theory, if $F_Y$ was known, one would compute the expected loss, i.e., the risk, of decision $a$, which essentially goes back to \cite{wald}. The action corresponding to the minimum risk is then interpreted as the best action. In fact, we have a three-player game here: The nature that generates the observations, the statistician that decides for some action based on the observations, and the referee that validates the action of the statistician.

Let us illustrate the situation in a simple estimation problem. Consider the situation that one has a sample of observations $y_1,...,y_n \in \mathcal{Y} \subset \R$ which are assumed to be realizations of i.i.d. random variables $Y_1$,...,$Y_n \sim F_Y$ and that the task is to estimate the expected value $\E_{F_Y}[Y_1]$. The action space of the statistician therefore is $A=\R$. The crucial question is whether there exists a method how to assess the quality of the estimation method so that one can decide which estimation method is optimal, provided that one has test data. For the mean estimation, it is well-known that the squared loss $L: \R^2 \rightarrow \R_{\ge 0}$, $L(y,\hat y)=(y-\hat y)^2$ is the correct choice as $\E_{F_Y}[L(Y,\hat y)]$ is minimal if and only if $\hat y=\E_{F_Y}[Y]$, provided that $F_Y \in \mathcal{F}_Y$ for the class $\mathcal{F}_Y$ of distributions on $\mathcal{Y}$ that allow for second moments. 

Denoting the mean functional by $T^{mean}: \mathcal{F}_Y \rightarrow \R$, $T^{mean}(F_Y)=\E_{F_Y}[Y]$, $Y \sim F_Y$, the squared error therefore is strictly consistent for $T^{mean}$ relative to $\mathcal{F}_Y$ since the inequality \begin{equation*} \E_{F_Y}[L(Y,t)]<\E_{F_Y}[L(Y,x)] \end{equation*} for $t=T^{mean}(F_Y)$ and $x \ne T^{mean}(F_Y)$ holds. Then, the mean functional $T^{mean}$ is called elicitable because such a strictly consistent loss function exists. This terminology essentially goes back to \cite{osband85} and has been made accessible for a broader audience by \cite{gnei10}.

\subsection{Example: Stochastic uncertainties from finite samples}

Consider the task of estimating the coefficients of a linear regression model $y_i=x_i\beta+\epsilon_i$ for $x_i \in \mathcal{X} \subset \R^p$, $\beta \in \R^p$, so that $y_i \in \mathcal{Y} \subset \R$ and $\epsilon_i \sim F_{\epsilon}$ i.i.d. for some error distribution $F_{\epsilon}$ with finite second moments, e.g., a Gaussian distribution. Assume that the instances $(x_i,y_i)$ are i.i.d. realizations from some joint distribution $F_{XY}$ on $\mathcal{X} \times \mathcal{Y}$. Then, given some training data, one can estimate the coefficients and evaluate the performance on a test set, quantified in a loss function $L: \mathcal{Y} \times \mathcal{Y} \rightarrow \R_{\ge 0}$. Naturally, it is desirable that the true coefficients correspond to the smallest average loss on the test set. 

The training set and the test set are assumed to be drawn from the same distribution $F_{XY}$, but of course, both sets are finite. The estimated coefficients are usually the minimizers of the average loss, i.e., the empirical risk, on the training set (e.g., \cite{vapnat}). The finiteness of the training set causes the empirical minimizer to differ from the theoretical population risk minimizer, which is well-known and captured by theoretical results on estimation consistency that guarantee, in some cases with a rate, the convergence of the empirical risk minimizer to the population risk minimizer.

As for the validation, the average test loss clearly also is nothing but an empirical version of the population risk corresponding to the selected test loss function. Therefore, on a small test set, it can happen that an optimal model does not appear as optimal, in other words, elicitability can be understood as an asymptotic property where the stochastic fluctuations are averaged out. 

\begin{figure}[H]
\begin{center}
\includegraphics[width=12cm]{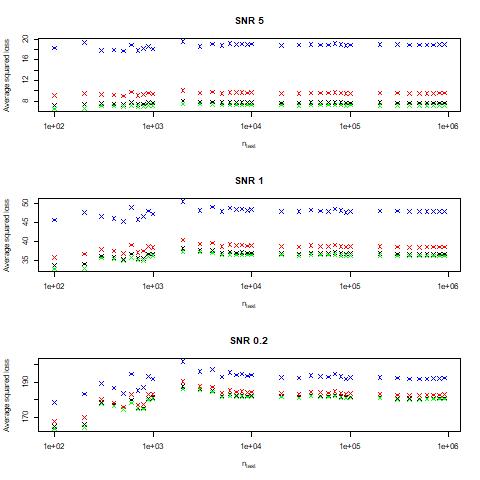}
\end{center}
\caption{Average test losses for an SNR of 5 (top), 1 (middle) and 0.2 (bottom), respectively. The blue, red and black crosses represent the performance of the estimated coefficients on the training data with an SNR of 0.2, 1 resp. 5, the green crosses the performance of the true coefficient vector.} \label{testlossEx0} 
\end{figure}

As an example, we simulate $n_{train}=250$ instances from the linear model $y_i=x_i\beta+\epsilon_i$, for $p=20$, $\beta_j \sim \mathcal{N}(1,1)$ i.i.d. for $j=1,...,20$, $\epsilon_i \sim \mathcal{N}(0,\sigma^2)$ i.i.d. for $i=1,...,n$, and where the $x_i$ are realizations from $X_i \sim \mathcal{N}_{20}((2,...,2),I_{20})$ i.i.d. for the identity matrix $I_{20} \in \R^{20 \times 20}$. The error variance is adapted so that the signal to noise ratio (SNR) is given by 5, 1 and 0.2, respectively. For each case, we compute the least squares coefficient. All three estimated coefficient vectors and the true coefficients are validated by a test set where $n_{val}$ instances are generated as before with the respective SNR. We use \begin{equation*} n_{val} \in \{100,200,...,900,1000,...,9000,10000,...,90000,100000,...,900000\}. \end{equation*} The results are depicted in Fig. \ref{testlossEx0}.

We can observe that, as expected, the average losses fluctuate for small $n_{val}$ but remain nearly constant for large $n_{val}$. Moreover, the true coefficient always corresponds to the best performance, independently from the noise level which solely controls the magnitude of the losses. Additionally, we can observe that the estimated coefficients result in higher average test losses the higher the noise level on the corresponding training data was. Note that elicitability generally does not imply that competing non-optimal models can be ranked in such a way but that this effect comes from the fact that the squared loss is an order-preserving loss function (see, e.g., \cite{stein}), indicating that $\E_{F_Y}[L(Y,t_1)]<\E_{F_Y}[L(Y,t_2)]$ for $t_2<t_1 \le T(F_Y)$ or $t_2>t_1 \ge T(F_Y)$.

\subsection{Example: Issues induced by contamination}

The previous example was based on clean data, in other words, it mainly focused on the two-player game between statistician and referee while not allowing the nature to produce challenging data. Since the validation risk is nothing but an expected value that is empirically approximated by the average test loss, all problems related to empirical means may fail in the case of contaminated test data. Let us illustrate this problem again with a linear regression example. 

\begin{figure}[H]
\begin{center}
\includegraphics[width=10cm]{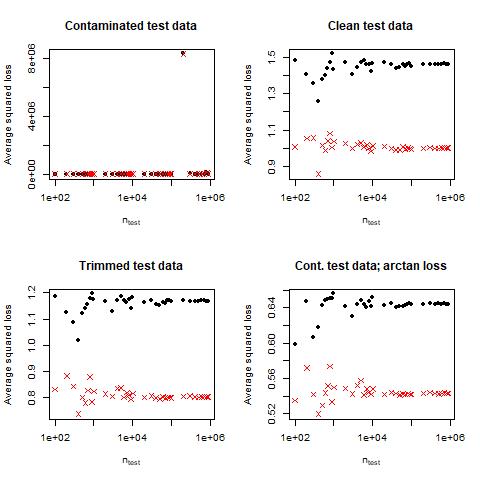}
\end{center}
\caption{Upper left: Mean squared loss on the whole test set; Upper right: An oracle has filtered any contaminated instance so that only non-contaminated instances enter the test loss (the effective test loss size therefore is smaller); Bottom left: Mean $\alpha$-trimmed test losses for $\alpha=0.05$; Bottom right: The losses were transformed via the arctan before being averaged. The losses corresponding to the true coefficient vector $\beta$ are coloured red, the losses corresponding to the least square estimate of $\beta$ in black.} \label{testloss} 
\end{figure}  

Consider a linear regression setting with contaminated responses, i.e., we model the responses by $y_i=x_i\beta+\epsilon_i$ as before but the error distribution $F_{\epsilon}$ is non-standard. Let $F_{\epsilon}=0.95 \mathcal{N}(0,1)+0.05 G$ for a Cauchy distribution $G$. In this example, we just compute the standard least squares estimator and compare its performance on test data with the performance of the oracle $\beta$. The performance of both models is then again evaluated on a test set whose instances are modeled identically as for the training set. The mean squared loss on the test set is depicted in Fig. \ref{testloss} for both models. 

The results are irritating. Seemingly, both models are equally well on the contaminated data (there are some points where the average loss for the true coefficient vector is indeed slightly higher than for the estimated coefficient vector) while when considering the clean data, trimming the results or arctan-transforming the losses, the oracle model is the best model. Disappointingly, all four procedures are either wrong, infeasible, or do not fully address the problem of contamination, letting both models become incomparable for the following reasons.

\begin{figure}[H]
\begin{center}
\includegraphics[width=10cm]{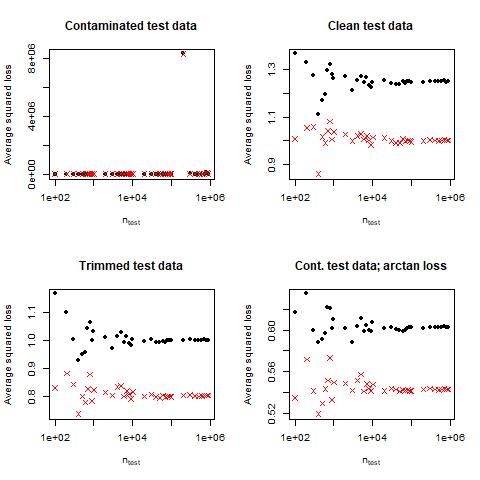}
\end{center}
 \caption{Upper left: Mean squared loss on the whole test set; Upper right: An oracle has filtered any contaminated instance so that only non-contaminated instances enter the test loss (the effective test loss size therefore is smaller but we nevertheless chose $n_{test}$ as $x$-coordinate); Bottom left: Mean $\alpha$-trimmed test losses for $\alpha=0.05$; Bottom right: The losses were transformed via the arctan before being averaged. The losses corresponding to $\beta$ are coloured red, the losses corresponding to the LTS estimate of $\beta$ in black.} \label{testlossslts}
\end{figure}  

The expected test loss which for the contaminated test data does not exist, so the average test loss cannot converge, making both models indeed incomparable through the lens of elicitability. As for the second strategy, the contaminated instances are not identifiable in practice. Although an extreme outlier most probably is most likely a realization from the Cauchy distribution, it is impossible to tell realizations from the normal distribution from ``innocently looking'' realizations from the body of the Cauchy distribution, disallowing to distinguish fully between observations from the normal and the Cauchy distribution, letting a trimming approach only be an approximation of the cleaning procdure. The fourth approach seems to be easily implemented in practice but is generally biased as the oracle model will not turn out as the optimal model which we will explain later. 

We repeat the experiment where we use the least trimmed squares (LTS, \cite{rous85}) estimator with a trimming rate of $\alpha=0.5$ instead of the least squares estimator. The results are depicted in Fig. \ref{testlossslts}. One can observe a similar behaviour here, although the models have been estimated robustly, indicating that contamination in test data is a serious problem against which robust model fitting cannot safeguard.

\subsection{Robust Statistics}

Robust Statistics (see, e.g., \cite{huber}, \cite{hampel}, \cite{maronna}, \cite{rieder}) have provided concepts for estimation under contamination of the underlying data, i.e., if the assumptions on the distribution of the data are not satisfied. The crucial idea is to bound the influence of contaminating data, which is done by outlier detection procedures that identify potential outliers which are trimmed away when computing the estimator, by bounding objective functions that are optimized, or by adaptive weighting strategies that do not trim outlying instances but assign less weight to them.

Robust Statistics provide two main quantifications of robustness of estimation procedures, namely the breakdown point (BDP) and the influence function (IF). The BDP  has been formally introduced \cite[Sec. 6]{hampel71} in a functional version and represents the minimum distance of the ideal and the contaminated distribution so that the estimation becomes unreliable. The finite-sample version, introduced in \cite{huber83}, represents the minimum fraction of outlying instances in a data set in order to arbitrarily perturb the estimator. The IF, introduced in \cite{hampel74}, in contrast is a local robustness measure that quantifies the impact of a single observation onto the estimator, which can be identified with a particular functional derivative (see \cite{averbukh}) of the statistical functional corresponding to the estimation problem. 

So far, to the best of our knowledge, Robust Statistics was mostly applied to estimation procedures or training processes, leading to estimators or models that are usually less perturbed, however, the validation of these estimators did not yet have been fully addressed, at least, the connection to the formal concept of elicitability has not yet been made. In other words, Robust Statistics therefore mainly focused on the two-player game between nature and statistician while paying less attention to the referee that has to validate the statistician's decision while the validation itself is based on challenging data, manifesting a two-player game between nature and referee.

\subsection{Contribution}

The main contribution of this paper is to establish the connection of Robust Statistics, validation and elicitability. 

Although we clearly are not the first ones who consider robust validation methods, we formally relate fair validation in terms of the notion of elicitability to robustness in terms of the notions of contamination and breakdown point.

We discuss different approaches that aim at circumventing the problem, showing why bounding the evaluation loss function or a modified notion of elicitability may fail. We propose a joint trimming approach in order to detect both outlying training and validation batches, which would maintain elicitability but require a perfect detection of contaminating instances. The result of our experiments indicate that fair validation in the presence of contaminated data is a serious problem,  which should be prioritized as future reseach topic.

This paper is organized as follows. In Sec. \ref{sec:prelim}, we compile the necessary background on elicitability and Robust Statistics. Sec. \ref{sec:issues} illustrates the validation problems in the presence of contaminated data and introduces the breakdown point of elicitability. Sec. \ref{sec:trimsec} discusses strategies for trimming test instances or whole test resamples away. These strategies are applied in Sec. \ref{sec:simsec} for different regression and classification scenarios.

\section{Preliminaries} \label{sec:prelim}

\subsection{Elicitability} 

Assume that the goal is to make predictions with a given model. Let $A$ be the action space that contains these predictions. The observations are contained in an observation space $O$ and the question is how to assess the quality of the prediction model. In order to reasonably do that, the following property has to be satisfied (\cite{gnei10}). \

\begin{Def} \label{consloss}  Let $A$ be some action space and let $O$ be some observation space. Let $\mathcal{F}$ be a family of probability distributions and let $L: A \times O \rightarrow \R$ be a loss function. Then, $L$ is called consistent for the statistical functional $T: \mathcal{F} \rightarrow A$ if \begin{equation} \label{cons} \E_F[L(t,Y)] \le \E_F[L(x,Y)]  \end{equation} for all $F \in \mathcal{F}$, $t \in T(F)$, $x \in A$. $L$ is called strictly consistent if equality holds if and only if $x \in T(F)$. \end{Def} \

Exactly as in the work of \cite{wald}, the model that achieves minimum risk, i.e., the model that predicts $x \in T(F)$, is regarded as the best model.   \cite{wald} assumed the existence of a loss function and then interpreted the model corresponding to the minimum expected loss as the best model. However, when designing the loss function, the crucial part is to guarantee that the property that has to be estimated indeed corresponds to the minimum risk and not a different quantity. Although this, at this first glance, looks evident, early results such as \cite[Cor. 2.5.1]{osband85}, where it is proven that the variance is not elicitable, already showed that elicitability is not a trivial property. Let us elucidate the consequence of this result: Although, for a particular data set, one can clearly estimate the variance, there is not chance to objectively compare different variance estimators. The path to a comparison of variance estimators has been paved in \cite{fissler} who have proven that the pair (mean, variance) is jointly elicitable, indicating that there exist loss functions that can objectively assess the quality of a joint mean and variance estimator.  



\subsection{Robust Statistics}

Robust Statistics considers contaminated data. First, we formally introduce the meaning of this notion (cf. \cite[Sec. 4.2]{rieder}). 

\begin{Def} \label{contball} Let $(\Omega, \mathcal{A})$ be a measurable space. Let $\mathcal{P}:=\{P_{\theta} \ | \ \theta \in \Theta\}$ be a parametric family so that $P_{\theta} \in \mathcal{P}_{\theta}$ is a distribution on $(\Omega, \mathcal{A})$. The set $\Theta \subset \R^p$ denotes some parameter space and the ideal distribution is denoted by $P_{\theta_0}$. A \textbf{contamination model} is defined as the family $\mathcal{U}_*(\theta_0):=\{U_*(\theta_0,r) \ |\ r \in [0,\infty[\}$ of \textbf{contamination balls} $U_*(\theta_0,r)=\{Q \in \mathcal{M}_1(\mathcal{A}) \ | \ d_*(P_{\theta_0},Q) \le r \}$ where $\mathcal{M}_1(\mathcal{A})$ denotes the set of probability distributions on $\mathcal{A}$. The radius $r$ the ``contamination radius''. \end{Def}

A standard contamination model is the so-called convex contamination model:

\begin{ex} \label{convcont} The contamination model $\mathcal{U}_c(\theta_0)$ that consists of contamination balls of the form $ U_c(\theta_0,r)=\{(1-r)_+P_{\theta_0}+\min(1,r)Q \ | \ Q \in \mathcal{M}_1(\mathcal{A}) \}$  is called convex contamination model. \end{ex} 

Consider a statistical learning problem where $n$ instances $(x_i,y_i)$ from a contaminated distribution are given. The finite-sample breakdown point from \cite{huber83} is given as the minimum fraction of contaminated instances required in order to achieve a complete distortion of the trained model.

\begin{Def} \label{sec:rob:bdp} Let $Z_n$ be a sample consisting of instances $(x_1,y_1),...,(x_n,y_n)$. The \textbf{finite-sample breakdown point} of the estimator $\hat \theta \in \Theta$ is defined by \begin{equation} \label{fsbdp} \varepsilon^*(\hat \theta,Z_n)=\min\left\{\frac{m}{n} \ \bigg\vert \ \sup_{Z_n^m}(||\hat \theta(Z_n^m)||)=\infty \right\} \end{equation} for the set $Z_n^m$ that has exactly $(n-m)$ instances in common with the original sample $Z_n$ and where $\hat \theta(Z_n^m)$ denotes the estimated coefficient on $Z_n^m$. \end{Def} 



\section{Robustness and elicitability}  \label{sec:issues}

\subsection{Issues arising from contamination}

Consider a standard normal distribution $F$ and a Cauchy distribution $G$ so that the convex contamination model considers the distributions $F_t=(1-t)F+t G$. Then, as moments of $F_t$ do not exist for any $t>0$, properties such as the mean are not elicitable relative to the class of such distributions. This may sound awkward as the mean of $F_t$ does not exist in the first place, making a discussion on the validation of a mean estimator seemingly obsolete. This is wrong, as the goal is to estimate the mean of the ideal distribution $F$, which exists by assumption. 

The consequence is that although one can robustly estimate the mean, elicitability is not valid in this setting, as the test data are usually also assumed to be realizations from $F_t$, or, at least, also from some contaminated distribution. Hence, competing (robust and non-robust) estimators for the mean are incomparable if $G$ is not sufficiently regular. Therefore, computing the test losses as in the upper left parts of Fig. \ref{testloss} and Fig. \ref{testlossslts} does not allow for comparison of competing estimators.

A standard procedure to assess whether a robust estimation is necessary is to compute a robust and non-robust estimator and to compare both. This can be done either by the estimates themselves so that one would decide for the non-robust one if they differ only marginally, or by computing in-sample or out-of-sample losses. Our argumentation elucidates that the latter strategy (including the evaluation by an in-sample loss) is in general not meaningful as, by the lack of elicitability, the decision which of the estimators is the best is not reliable.

\subsection{Trimming and regression}

Trimming is a standard technique in Robust Statistics which already has entered regression, e.g., in the least trimmed squares (LTS, see \cite{rous06}) or the sparse LTS (SLTS, see \cite{alfons14}). We distinguish between two fundamentally different strategies. 

Consider a linear regression setting with the model $y_i=x_i\beta+\epsilon_i$, for realizations $y_i$ of $Y_i \sim F_{Y|X=x_i}$, realizations $x_i$ of $X_i \sim F_X$, $\beta \in \R^p$, and an error term $\epsilon_i \sim F_{\epsilon}$ for some error distribution $F_{\epsilon}$.

\textbf{1.)} If the goal is to model a trimmed expectation of $Y_i$ given $x_i$, the instance-specific trimmed means $\E[Y_i|X=x_i]$ are modeled. The trimmed mean is covered by the more general range value at risk (RVaR) in \cite{fissler21b}, and shown to be not elicitable. In this trimmed regression setting, one would therefore have to jointly estimate the trimmed mean and two quantiles, as proven in \cite{fissler21b}.

\textbf{2.)} Aiming at validating an estimator for $\E[Y|X]$, one would consider the mean squared error $\E_X[\E_{Y|X}[(T^{mean}(F_{Y|X})-Y)^2]]$ (cf. \cite{brehmer}, \cite{fissler21b}), so one would empirically approximate this quantity by computing \begin{center} $ \displaystyle \frac{1}{n_{test}}\sum_{i=1}^{n_{test}} (y_i-x_i\hat \beta)^2 $ \end{center} on a test set with $n_{test}$ instances $(x_i,y_i)$, and where $\hat \beta$ is an estimate for $\beta$. If the goal is now to trim suspicious test instances away, for example, by computing \begin{equation} \label{trimsq} \frac{1}{\lceil (1-\alpha)n_{test} \rceil}\sum_{i \in \{1,...,n_{test}\} \setminus I^{trim}(\alpha)} r_i^2, \ \ \ r_i:=y_i-x_i\hat \beta \end{equation} for trimming rate $\alpha$ and the set \begin{equation*} I^{trim}(\alpha):=\left\{i \in \{1,...,n_{test}\} \ \bigg| \ \sum_{k=1}^{n_{test}} I(r_k^2 \ge r_i^2) \le \lfloor n_{test}\alpha \rfloor \right\}, \end{equation*} one would still validate the mean (with a trimmed test set), and not a trimmed mean as in variant 1.). Therefore, the issue that a trimmed mean is not elicitable alone does not apply here, so elicitability of the mean by evaluating the strict consistent squared loss is still valid. Here, the worst test instances (w.r.t. the estimated coefficient vector $\hat \beta$) are trimmed away in the two-player game between nature and referee, similarly as the worst training instances are trimmed away in the two-player game between nature and statistician. 


\subsection{Possible ``robust'' evaluation methods and problems} \label{sec:methods}

Now, we come back to the different evaluation methods that we used in the experiment corresponding to Fig. \ref{testloss}, namely: i) removing contaminated instances; ii) trimming; iii) robustifying, i.e., bounding, the test loss function. 

As for i), for the moment, assume that an oracle that correctly identifies instances from $G$ is indeed given. Then, the following simple example reflects the situation in the upper right part of Fig. \ref{testloss}.

\begin{ex} Let $Y_i \sim F_t=(1-t)F+tG$ for some distribution $G$ that does not allow for second moments and a distribution $F$ that allows for second moments, and let $t \in ]0,1[$. If an oracle procedure is able to set weights $w_i=I(Y_i \ \text{stems from}\ F)$, the average loss \begin{equation*} \frac{1}{n}\sum_{i=1}^n w_i (Y_i-\mu)^2 \end{equation*} converges indeed to $\E_F[(Y_1-\mu)^2]$, so the desired type of elicitability, i.e., eliciting the mean of the non-contaminated distribution $F$, is valid. In general, any property that is elicitable w.r.t. $F$ remains elicitable after this cleaning step.\end{ex}

As such an oracle is not given in practice, removing all instances from $G$ is of course not feasible on real data.  Existing outlier detection approaches for trimming do not have the property of reliably flagging out instances that stem from $G$, in fact, they are more likely to also flag instances as outliers that stem from the tails of $F$ while realizations from the the body of $G$ may not be detected. The former property would not be problematic as one could think of weights $w_i$ in the example above that are zero for all instances from $G$ but also for some (not all) instances from $F$, which would only reduce the convergence speed but not still allow for strict consistency.



As for ii) and iii), first consider the following example.

\begin{ex} In a rather artificial example, consider a linear regression model $y_i=\beta_0+x_i\beta+\epsilon_i$ with $\epsilon_i \sim 0.5U([-20,-10])+0.5U([10,20])$. Then, applying unsuitable robust losses of the form $L(y_i,\hat \beta_0+x_i\hat \beta)=L(|y_i-\hat \beta_0-x_i\hat \beta|)$ that remain constant if the input exceeds the value 5 would make all regression lines which are affine shifts of the ideal regression line by a shift of at most 5 incomparable as the losses on all instances would be trimmed.  \end{ex}

The example shows that a clipped loss function that is constant outside some finite interval could be misleading. Alternatively, one can strictly isotonically transform the losses using some transformation $h$ with a bounded image space. This has been done in the bottom right part of Fig. \ref{testloss} for $h=\tan^{-1}$. More formally, one computes $ \frac{1}{n}\sum_i h(L(y_i,\hat y_i))$ so that, for a fixed $h$, strict consistency of $h \circ L$ is guaranteed if \begin{equation} \label{elictrafo} \E_{F_Y}[h(L(Y,t))]<\E_{F_Y}[h(L(Y,x))] \end{equation} for $t=T(F_Y)$ and $x \ne T(F_Y)$.

However, due to $\E_{F_Y}[h(L(Y,t))] \ne h(\E_{F_Y}[L(Y,t)])$, having a strictly consistent loss function $L$, one cannot conclude that inequality \ref{elictrafo} holds. Using an affine linear function $h$ would not solve the problem that regularity conditions may not be satisfied by the contaminated distribution.


In the remainder of this paper, we formalize a breakdown point for elicitability and empirically assess the attained BDP for the trimming strategy on simulated data. As for trimming, one could, from the perspective of model validation, interpret a referee that uses a trimmed empirical mean for validation as a supporter of a statistician in the presence of a challenging nature in the sense of a three-player game where the pair (statistician, referee) plays against the nature.

\subsection{Breakdown of elicitability}


\begin{Def} Let $\mathcal{F}$ be a class of distributions on an observations space $O$ and let $T: \mathcal{F} \rightarrow A$ be an elicitable statistical functional for an action space $A$. Let $\mathfrak{L}(T,\mathcal{F})$ be the set of all loss functions $L: A \times O \rightarrow \R$ that are strictly consistent for $T$ w.r.t. $\mathcal{F}$. \begin{itemize}

\item[a)] Strict consistency of $L \in \mathfrak{L}(T,\mathcal{F})$ for $T$ w.r.t. $\mathcal{F}$ breaks down for contamination radius $r \in [0,\infty[$ if $\exists F \in U_*(\mathcal{F},r): \E_F[L(t,Y)] \ge \E_F[L(x,Y)]$ for $t \in T(F)$, $x \notin T(F)$, the set $U_*(\mathcal{F},r)=\bigcup_{F \in \mathcal{F}} U_*(F,r)$  of contamination balls $ U_*(F,r)=\{Q \in \mathcal{M}_1(\mathcal{A}) \ \vert \ d_*(F,Q) \le r \}$ for $F \in \mathcal{F}$.

\item[b)] Elicitability of $T$ w.r.t. $\mathcal{F}$ breaks down for contamination radius $r$ if strict consistency breaks down for all $L \in \mathfrak{L}(T,\mathcal{F})$. \end{itemize} \end{Def}

\begin{Def}\label{bdpcons} Let $k$ random variables $Y_1,...,Y_k$ be given so that $0 \le m \le k$ of them are distributed according to some $F \in U_*(\mathcal{F},r)$ and $(k-m)$ according to some $F_0 \in \mathcal{F}$. The \textbf{breakdown point for strict consistency for a loss function} $L \in \mathfrak{L}(T,\mathcal{F})$, $L: A \times O \rightarrow \R$ for action domain $A$ and observation domain $O \subset \mathcal{Y}^k$, is given by \begin{equation} \label{fsbdpL}\begin{split} \varepsilon^*(L,k,\mathcal{F},r)=\min\{\frac{m}{k} \ \bigg\vert \ \exists Y_{i_s} \sim F_r^{(s)} \in U_*(\mathcal{F},r), s=1,...,m, \\ 1\le i_1<i_2<...<i_m \le k, Y_{j_{s'}} \sim F_0^{(s')}, F_0^{(s')} \in \mathcal{F}, s'=1,..., k-m, \\ 1\le j_1<j_2<...<j_{k-m} \le k, i_s \ne j_{s'}: \\ \E[L(t,(Y_1,...,Y_k)] \ge \E[L(x,Y_1,...,Y_k)] \ \exists x \notin T(F), t \in T(F) \}. \end{split} \end{equation} \end{Def}

In other words, the BDP for strict consistency for $L$ is the minimum relative number of input objects for $L$ from the observation domain that are realizations from a contaminated distribution that can destroy the strict consistency of $L$. One statement trivially follows for unbounded loss functions:

\begin{prop} In the setting of Def. \ref{bdpcons}, an upper bound for the BDP for strict consistency of an unbounded loss function $L$ is given by \begin{equation} \label{fsbdpLE} \begin{split} \varepsilon^*(L,k,\mathcal{F},r)=\min\{\frac{m}{k} \ \bigg\vert \ \exists Y_{i_s} \sim F_r^{(s)} \in U_*(\mathcal{F},r), s=1,...,m, 1\le i_1<i_2<...<i_m \le k,\\ Y_{j_{s'}}  \sim F_0^{(s')}, F_0^{(s')} \in \mathcal{F}, s'=1,...,k-m,1\le j_1<j_2<...<j_{k-m} \le k, i_s \ne j_{s'}: \\ \sup_{F_r^{(1)},...,F_r^{(s)},F_0^{(1)},...,F_0^{(s')}}(\E[L(t,(Y_1,...,Y_k)])=\infty \}. \end{split} \end{equation} \end{prop}

\begin{ex}\label{bdpex}  Loss functions of the type considered in Def. \ref{bdpcons} where more than two input arguments enter arise in the context of ranking problems (\cite{TW19b}). In a pair-wise comparison, one essentially has four input arguments, namely the predictions $\hat y_i$ and $\hat y_j$ for the true responses $y_i$ and $y_j$, respectively, and the true responses, so that the loss function is given by \begin{equation}  \label{rankloss} L(y_i,y_j,\hat y_i,\hat y_j)=I((y_i-y_j)(\hat y_i-\hat y_j)<0),  \end{equation} i.e., a misranking occurs if the ordering of $\hat y_i$ and $\hat y_j$ does not match the ordering of $y_i$ and $y_j$. 

The existence of the expectation of the individual random variables $Y_i$ is not necessary for strict consistency since the loss function itself may be bounded, which is true for the loss function from Eq. \ref{rankloss}, which has been shown to be strictly consistent for ranking in \cite{TW21}. Here, the expectations of $Y_i$ or $Y_j$ need not exist. However, the BDP for strict consistency is not $1$ but $0.5$, which follows from \cite{TW21}, as strict consistency is not given if $Y_i$ and $Y_j$ are identically distributed, i.e., contaminating one of both random variables suffices in order to let strict consistency break down.  \end{ex}

An interesting question is whether it is possible to achieve a breakdown of elicitability by other means than by contaminating the distributions in such a way that the necessary expectations of the loss function no longer exist. A characterization of contamination balls that let elicitability break down for each functional that already has been proven to be elicitable is out of scope for this work, so we leave this topic open for future research.

\subsection{Empirical breakdown point of elicitability}

Ex. \ref{bdpex} motivated that the BDP for strict consistency of a loss function can be larger than 0. Therefore, we assume in general that a loss function has a BDP for strict consistency of $c \in ]0,1]$. However, this BDP is only one component of a true BDP for elicitability, which has to consider the aggregation over the test data.

Assume that, on a test data set with $n_{test}$ instances, one evaluates the loss function $c(n_{test},k)$ times. For example, a ranking loss function with $k$ input pairs of the form $(y_i,\hat y_i)$, a natural strategy is to approximate the test risk via a U-statistic (e.g., \cite{TW19b}), i.e., \begin{equation} \label{ustat} \frac{1}{\binom{n_{test}}{k}}\sum_{i_1<i_2<...<i_k \in \{1,...,n_{test}\}} L(y_{i_1},...,y_{i_k},\hat y_{i_1},...,\hat y_{i_k}), \end{equation} corresponding to $c(n_{test},k)=\binom{n_{test}}{k}$. The standard validation procedure where only one pair $(y_i,\hat y_i)$ enters the loss function corresponds to $c=1$, $k=1$, and $c(n_{test},k)=n_{test}$. Due to the arithmetic mean aggregation in Eq. \ref{ustat}, the empirical BDP of strict consistency on the test data using loss function $L$ is given by $\frac{c}{c(n_{test},k)}$, which, for the natural property that $c(n_{test},k) \rightarrow \infty$ for $n_{test} \rightarrow \infty$, converges to zero. Here, one does not even have to distinguish between a breakdown of strict consistency or of elicitability, as the empirical BDP of elicitability on the test data w.r.t. a  family $\mathcal{L}$ of loss functions would be given by $\sup_{L \in \mathcal{L}}\left(\frac{c(L)}{c(n_{test},k)}\right)$, where $c(L)$ denotes the BDP of loss function $L$ for explicitness, which still converges to zero for growing $n_{test}$.

For a loss function with $k>1$, in the context of contaminated data, there may be $k$-tuples with $j=0,1,...,k$ contaminated instances. As strict consistency would only break down for $j \ge \lceil ck\rceil$, one could ask whether the problematic $k$-tuples can be ignored, i.e., trimmed. Before we go into details, we first consider resampling. Due to resampling in cross-validation procedures, contaminated instances may appear in multiple training batches and in multiple validation batches. Here for $B$ resamples with $n_{val}$ instances, as the cross-validated error is the mean of all mean test losses on the resamples, one has $Bc(n_{val},k)$ evaluations of the loss function in total, which still corresponds to an asymptotic BDP of zero. Similarly as for the $k$-tuples, some resamples may contain a larger fraction of contaminated instances than other. This idea has been used in \cite{TW21c} for trimming training resamples away in a robust Stability Selection. In the context of validation, one could ask whether whole test resamples can be trimmed away, which can be interpreted as an outer trimming, in contrast to inner trimming when ignoring single $k$-tuples.

\section{Illustration of a trimming procedure} \label{sec:trimsec}

As cleaning the test data preserves elicitability, as discussed in Sec. \ref{sec:methods}, in contrast to a transformation of the loss function or a trimmed loss function, we want to illustrate how a trimming procedure that should approximate the cleaning procedure can be done in practice, where both the training and the test set can be assumed to be contaminated.

\subsection{Impact of contamination in training and test data}

We make one standard assumption: The contaminated instances follow an unknown distribution that differs from the ideal distribution. Most importantly, outlying instances generally do not have any structure in the sense that the response variable and the predictor variable are related by some model. 

Then, we can distinguish the following four cases: \begin{enumerate}

\item \textbf{Clean training and test data:} A correctly specified model should correspond to low training and test losses. As usual, overfitted models may be more likely to correspond to very low training but larger test losses. Robust models are expected to show the same behaviour due to the absence of contamination. 

\item \textbf{Contaminated training data, clean test data:} Here, due to the contaminated instances in the training data that do not follow the model that relates the predictor and response variables on the clean subset of the data, both classical and robust models should correspond to large in-sample losses here. In contrast, the test loss should be large for classical models due to having been distorted by the contaminated instances, while the robust model should correspond to a low test loss. 

\item \textbf{Clean training data, contaminated test data:} As the training data are clean, robust and classical models should not differ much and therefore correspond to similar (low) training losses, but to large test losses. 

\item \textbf{Contaminated training and test data:} For both models, both the training and the test losses are expected to be high. Robust models should approximate the map between the predictor and response variables well for the instances from the ideal distribution. Due to outlying instances, the average test losses should also be large for the robust model, but it should perform rather well on the clean test instances.  

\end{enumerate}

\subsection{Trimming strategy}

Outlier detection procedures aim at identifying instances in a data set that differ from the majority of the data by some metric. In particular, approaches such as LTS and SLTS use the in-sample losses for trimming, i.e., one iteratively searches for the $(1-\alpha)$-fraction of instances for which the current in-sample loss is minimal, i.e, which are fitted best by the current model. 

This strategy can be used for trimming test data by computing a robust model on the training data and by evaluating a loss function on each test instance. Then, the test instances with the largest losses are trimmed away and do not enter the aggregated test loss.

A similar strategy can be used in order to clean the training data. We propose a model-based leave-one-out (LOO) approach. The idea is to train a model on all training instances except one, and compute the loss for the remaining one. Iterating this procedure $n$ times where the $i$-th instance, $i=1,...,n$, is left out, respectively, the training instances with the highest losses are trimmed. 

In cross-validation, where one does not consider one test set but multiple test sets, \cite{khan10} proposed that for each of the test folds, a trimmed loss is computed. Finally one averages these trimmed losses over all test folds. Effectively, in each of the folds, a $\alpha$-fraction of the instances is ignored when computing the loss. The problem, through the lens of elicitability, is that once at least one fold has a higher contamination rate than $\alpha$, the whole procedure breaks down, as the batch-wise losses are averaged by the arithmetic mean, which has a BDP of zero. In the context of resampling, even setting $\alpha=0.5$ would not remedy this problem since the contamination rate on a resample can exceed $0.5$, even if the contamination rate on the original data set is at most $0.5$. 

Note that the LOO approach is not affected by this issue and allows to keep the assumption that the contamination rate on the complete data set is smaller than 0.5, safeguarding robust algorithms with a BDP of 0.5 from being distorted.

We aim at comparing several heuristics, which, due to neglecting whole resamples, could potentially even work in situations where the procedure from \cite{khan10} would fail. 


One strategy for trimming whole resamples away is to compute in-sample or out-of-sample losses after having estimated a model. Therefore, trimming away the $\alpha$-fraction of models with the highest in-sample or out-of-sample loss can be regarded as a natural way to extend the application of trimming from instances to whole data resamples, similarly as in \cite{TW21c}. This strategy is not conflicting with the lack of elicitability as the models are not compared, but implicitly the data resamples (i.e., one would never compare the performance of two or more models but the performance of one model on different data resamples).

Another idea arises from the standard diagnostic concept of Robust Statistics when facing the question whether to use the robust model or the classical model on a given data set after having trained both. Usually, one compares both models and if the deviation between the models is small, one assumes that the training data is (sufficiently) clean, allowing for using the classical model. A crucial question is how to quantify these deviations. A natural way would be the direct comparison of model parameters, i.e., the deviation $||\hat \beta-\tilde \beta||$ in linear regression model, where $\hat \beta$ and $\tilde \beta$ are the classical and robust regression coefficient estimate, respectively. Such deviations have, for example, even been quantified for neural networks with the same architecture, see \cite{goldberger}. However, in particular for non-linear models, slight changes in the model parameters can have a significant impact on the predictions, making it difficult to order the model deviations when identifying the most deviating $\alpha$-fraction. 

Summarizing, we propose the following approach: \ \\

\textbf{1a) Trimming of training instances:} Trim away training instances according to the LOO-strategy. \\
\textbf{1b) Trimming of training batches (for cross-validation):} Compute the classical and robust models on each training resample and trim away training resamples. \\
\textbf{2a) Trimming of test instances (for a static test set):} Use the remaining training samples, train a robust model, compute the out-of-sample losses on the test data, and trim away test instances. \\
\textbf{2b) Trimming of test batches (for cross-validation):} Compute the average test loss on each individual test set and robustly aggregate them by trimming away test sets.

\section{Simulation study} \label{sec:simsec}

\subsection{Data} 

For both linear regression and classification, the regressors $X_i \in \R^p$ are always generated according to a $\mathcal{N}_p(\mu_p,I_p)$-distribution where $\mu_p=(\mu,...,\mu) \in \R^p$ and where $I_p$ is the identity matrix of dimension $p$. The true non-zero coefficients $\beta_j$ are $\mathcal{N}(1,1)$-distributed. The number $s^0$ of non-zero coefficients is set in advance, and the non-zero components of the coefficient vector $\beta$ are randomly selected.  As for the noise variables $\epsilon_i$, $i=1,...,n$, for linear regression, they are normally distributed with mean $0$ and a variance such that a given signal-to-noise ratio (SNR) is attained. For binary classification, we first compute $\check Y_i=X_i\beta$ and $\bar Y=\frac{1}{n}\sum_i \check Y_i$ and define $\tilde Y_i=\check Y_i-\bar Y_i$. We then compute $p_i=\exp(\tilde Y_i)/(1+\exp(\tilde Y_i))$ so that the binary response $Y_i$ is generated according to a $B(1,p_i)$-distribution, $i=1,...,n$. The centering step, i.e., the computation of the $\tilde Y_i$, is done in order to encourage a balanced data set, avoiding that either class may be underrepresented. 

There are two types of contamination in our experiments: Case-wise $Y$-contamination and cell-wise $X$-contamination. As for case-wise $Y$-contamination, we generate a realization from a $B(n,r)$-distribution where $r$ is the contamination radius so that all responses corresponding to a component 1 are replaced by a gross outlier, which is the value 50 in all cases. As for cell-wise $X$-contamination, we first randomly select $\lfloor rn \rfloor$ instances and for each of the selected instances, $\lfloor 0.1p\rfloor$ of the components are randomly selected for replacement with a gross outlier, which again is the value 50.

For each scenario, we consider different data configurations that are specified in Tab. \ref{conftablereg} and Tab. \ref{conftableclass}. Here, $n_{test}$ denotes the number of test instances in randomized cross-validation settings, $n_{sub}$ denotes the size of the batches in randomized cross-validation. For $K$-fold cross-validation, we always use $K=10$ and $K=5$. Each scenario is repeated $V=1000$ times for the data with $p=20$ and $p=250$, and $V=500$ times for $p=500$. When trimming away instances, we always use the trimming rate $\alpha=r$.

\begin{table}
\begin{center}
\begin{tabular}{|p{0.4cm}|p{0.4cm}|p{0.5cm}|p{0.5cm}|p{0.4cm}|p{0.4cm}|p{2cm}|p{4cm}|} \hline 
 $p$ & $n$ & $n_{test}$ & $n_{\text{sub}}$ & $s^0$ & $\mu$ & SNR & $r$   \\ \hline
 20 & 250 & 100 & 125 & 20 & 2 & $\{0.5,2,5\} $ & $\{0.05,0.15,0.25,0.5\}$   \\ \hline
250 & 100 & 50 & 50 & 15 & 2 & $\{0.5,2,5\} $ & $\{0.05,0.15,0.25,0.5\}$  \\ \hline
 500 & 100 & 50 & 50 & 15 & 2 & $\{0.5,2,5\} $ & $\{0.05,0.15,0.25,0.5\}$  \\ \hline
\end{tabular}
\end{center}
\caption[Scenario specification]{Scenario specification for regression data} \label{conftablereg}
\end{table} 

\begin{table}
\begin{center}
\begin{tabular}{|p{0.4cm}|p{0.4cm}|p{0.5cm}|p{0.5cm}|p{0.4cm}|p{2cm}|p{4cm}|} \hline 
  $p$ & $n$ & $n_{test}$ & $n_{\text{sub}}$ & $s^0$ & $\mu$ & $r$  \\ \hline
 20 & 250 & 100 & 125 & 20 & \{0.5,3,8\} & $\{0.05,0.15,0.25,0.5\}$  \\ \hline
250 & 100 & 50 & 50 & 15 & \{0.5,3,8\} & $\{0.05,0.15,0.25,0.5\}$  \\ \hline
500 & 100 & 50 & 50 & 15 & \{0.5,3,8\} & $\{0.05,0.15,0.25,0.5\}$  \\ \hline
\end{tabular}
\end{center}
\caption[Scenario specification]{Scenario specification for classification data} \label{conftableclass}
\end{table}

\subsection{Algorithms} 

For data with $p=20$, we apply standard linear regression using the \texttt{lm}-function in $\mathsf{R}$ and LTS, provided by the \texttt{ltsReg} function from the package \texttt{robustbase} (\cite{robustbase}, \cite{robustbase2}). We set $\alpha=0.5$ for LTS. For classification, we apply the functions \texttt{glm} and the robust counterpart \texttt{glmrob} from the package \texttt{robustbase}, using the robust fitting method \texttt{Mqle}. 

For data with $p=250$ and $p=500$, we apply the Lasso, using the \texttt{glmnet}-function from the \texttt{glmnet}-package (\cite{friedman10}), and the robust counterpart SLTS, provided by the function \texttt{sparseLTS} from the package \texttt{robustHD} (\cite{alfons16}). As for classification, we apply LogitBoost, provided by the function \texttt{glmboost} from the package \texttt{mboost (\cite{mboost})}. As there does not yet seem to be an implemented robust counterpart for high-dimensional data, we restrict ourselves to the AUC-Boosting, setting the \texttt{family}-argument in \texttt{glmboost} to \texttt{AUC()}, being aware that this algorithm is generally not robust either.

\subsection{Trimming resamples} \label{sec:trimresample}

\subsubsection{Contaminated training data, clean test data}  \label{sec:valclean}

In this preliminary setting, we have clean test data but contaminated training data. This test set is also used for $K$-fold CV, as the usual approach where the complement of the training batch is used for validation would not allow for clean test data here. The goal is to trim away CV-batches. Ideally, from an algorithmic perspective the CV-batches with the largest relative fraction of outlying instances would be trimmed away. In other words, we in fact have a ranking problem here, w.r.t. the outlyingness of the individual batches. 

For ranking problems, \cite{TW21b} showed that the hard ranking loss function (cf. Ex. \ref{bdpex}, now in a more general notation) \begin{center} $ \displaystyle L(\pi,\hat \pi)=\frac{1}{K(K-1)}\sum \sum_{k \ne l} I((\pi_k<\pi_l)(\hat \pi_k<\hat \pi_l)) $ \end{center} for a vector $\pi$ of decreasing ranks (i.e., $i$ with $\pi_i=1$ is the index of the most outlying batch) and a predicted counterpart $\hat \pi$ is strictly consistent. Therefore, given $K$ CV-resamples, we first order the individual resamples according to their outlyingness in terms of the relative fraction of contaminated instances, resulting in the true vector $\pi$ of ranks. Then, we order the batches according to the average in-sample loss, the average test loss, or the norm differences of the coefficient vectors of a robust and non-robust model, respectively, resulting in a predicted vector $\hat \pi$ of ranks for each outlyingness measure. We finally evaluate the hard ranking loss function w.r.t. $\pi$ and each of the predicted counterparts, and plot the average hard ranking losses over all $V$ repetitions individually for each configuration of contamination radius, contamination scheme, SNR, CV-scheme, and regression/classification algorithm. 

As for the training and test losses, we use both the usual squared loss for regression and binomial likelihood loss for classification, respectively, but also a trimmed version where the largest $0.5$-fraction of losses is ignored on each CV-resample.

The figures are provided in the Appendix in Sec. \ref{app:valclean}.

For regression, we can observe that the difference between the results for $5$-fold and $10$-fold CV as well as for randomized CV with $10$ and $100$ batches is negligible, respectively. 

For randomized CV, when using the training loss in the context of $Y$-contamination, the best results are achieved when using a non-trimmed loss and a robust regression model, the worst results when combining both a robust regression model and a trimmed loss. When using non-trimmed losses, the hard ranking errors monotonically increase as $r$ increases, while slightly decreasing for $r=0.15$ and $r=0.25$ and increasing again for higher $r$. When using the test losses, the performance unexpectedly increases for all strategies mostly for growing $r$, while using a non-robust regression algorithm provides the best results.

 In the presence of $X$-contamination, the combination of a trimmed loss and a robust regression model provides the most stable results for $p=20$, being much worse than those for the other loss and regression combinations for low $r$, and better for large $r$. For larger $p$, robust regression with non-trimmed losses works best. When using test losses, using non-robust regression and a non-trimmed loss leads to the best results.
 
 For $K$-fold CV, the hard ranking errors generally decrease with for large $r$ when using in-sample losses, and increase when using test losses. As for test losses, using non-trimmed losses leads to the best performance in the presence of $Y$-contamination, while for $X$-contamination, all four strategies lead to comparable results.
 
 The large difference in the results for $K$-fold and randomized CV may be explained by the size of the batches, which is $n/2$ for the latter strategy, but only $n/K$ for $K$-fold CV.

For classification, the results when using trimmed and non-trimmed training or test losses are comparable, often identical, therefore, the figures often contain only two graphs. 

For randomized CV, there is neither a considerable difference between the case of $10$ and $100$ batches, nor between LogitBoost and AUC-Boosting. At least, LogitBoost works slightly better when using training losses than AUC-Boost. This also holds for $K$-fold CV.

Note that these preliminary experiments have yet nothing to do with elicitability, but should serve as a baseline if one would like to trim CV-procedures by trimming CV-batches away. Given that the test set is clean, the results are indeed disappointing.

\subsubsection{Contaminated training data, post-trimming}  \label{sec:traintrim}

In this subsection, we consider only the $k$-fold CV in the context of contaminated training data. In contrast to the experiments before in Sec. \ref{sec:valclean}, we first trim the data according to the model-based leave-one-out (LOO) approach. We repeat the analysis of the ranking errors of the batches based on the trimmed training data. 

The goal is to investigate whether the identification of outlying batches is more successful if outlying instances are removed beforehand.

The figures are provided in Sec. \ref{app:traintrim}.

The results differ from those in Sec. \ref{sec:valclean} as the shapes of the curves do not match. It seems that the trimming of the training data according to a robust regression model leads to a more coherent training set in the sense that the identification of outlying batches is facilitated, or that even batches are not contaminated due to having filtered out many outlying instances in advance, which would explain why the ranking error is very low for small contamination radii. In all cases, the four curves are much closer to each other than in Sec. \ref{sec:valclean}. For $r=0.5$, the ranking errors are again around 0.5 for all methods.

\subsubsection{Contaminated test data, post-trimming} \label{sec:valtrim}

Now, we only consider the randomized CV in the context of both contaminated training and test data. As for the data configurations, we must therefore include a contamination radius $r_{val}$ for the test data. We consider the 9 configurations $(r,r_{val}) \in \{0.1,0.25,0.5\}^2$. The other data configuration components remain as specified in Tab. \ref{conftablereg} and Tab. \ref{conftableclass}.

The idea is to first train a model on the training data and to compute the model-based outlyingness for each single test instance, and to trim the most outlying test instances. Afterwards, randomized CV is performed by sampling batches from the training data and training a model on each batch. Then, each batch is validated according to the test loss on the cleaned test data. 

The goal of this set of experiments is to study whether the performance is worse when the identification of the most outlying batches is done w.r.t. the cleaned data than with truely clean data as in the first set of experiments in Sec. \ref{sec:valclean}. Since $K$-fold CV normally does not include a separate test set, we restrict the experiments here on randomized CV. The evaluation methods are the same as for the first set of experiments.

The figures are provided in the Appendix section \ref{app:valtrim}.

As for the test losses, one can again observe well that non-robust regression leads to the best results in the context of $Y$-contamination, while the ordering varies for $X$-contamination. For $p=20$ robust methods often lead to slightly better results than non-robust methods, in contrast to the experiments in Sec. \ref{sec:valclean}.

For classification, there are only negligible differences between LogitBoost and AUC-Boosting and Logit and GLM-rob, respectively. The differences to the experiments in Sec. \ref{sec:valclean} are also negligible.

In summary, the results are comparable with those on clean test data, however, one should more frankly say ``comparably bad''.





\subsection{Trimming single instances} \label{sec:triminstance}

\subsubsection{Contaminated training and test data, pre-trimming} \label{sec:valcont}

Again, both the training and the test data are contaminated. The goal is two-fold: \textbf{i)} Identify the outlying test instances based on a model trained on the contaminated training data; \textbf{ii)} Identify the contaminated training instances by the LOO approach. Both identifications are done separately here, i.e., the trimming of the test data is based on the full training data.
 
As for the performance metric, we only need an identification of the outlying training or test instances. Therefore, we have a so-called weak ranking problem (see \cite{TW19b} for the terminology). In the previous section, all batches were potentially contaminated, moreover, there were ties, therefore, we considered finding the ranking of the batches as a hard ranking problem. Here, we use the weak ranking loss function \begin{center} $ \displaystyle L_K^{weak;M}(\pi,\hat \pi)=\frac{2}{K}\sum_{k \in Best_M} I(\hat \pi_k>M), $ \end{center} where $M$ is the number of instances that is flagged as outlying and $Best_M$ is the index set corresponding to the true $M$ outlying instances. In other words, one just counts the relative fraction of true outlying instances that are not flagged as such. This loss function has been shown to be strictly consistent for weak ranking in \cite{TW21b}.  

The figures are provided in the appendix section \ref{app:valcont}.

As for the identification of the contaminated training instances and $p=20$, as expected, the ranking errors increase with increasing contamination radius $r$. In addition, robust regression leads generally to smaller ranking errors, except for the case of $Y$-contamination and $r=0.5$, where they are comparable with those from non-robust regression. $X$-contamination leads to worse results than $Y$-contamination, except for robust regression and $r=0.5$. A lower SNR generally leads to worse results. These observations are not surprising as one can expect that an outlying instance can be better identified with a good model than with a model that is distorted itself, and low contamination radii, robust algorithms and a high SNR make it more likely to get a good model. For $p=250$ and $p=500$, the differences between the performance of robust and non-robust models for identifying contaminated training instances are negligible. 

The results for the identification of the contaminated test instances are similar. As we have two contamination radii here, we can observe increasing ranking errors wit increasing $r$, but not necessarily for increasing $r_{val}$, on the contrary, the ranking errors often even decrease with increasing $r_{val}$, both for robust and non-robust regression. The empirical BDP however always increases with increasing $r$ and $r_{val}$, and is often significantly higher for non-robust regression than for robust regression. Note that for $X$-contamination and large $p$, the empirical BDP is often close to 1, indicating that in nearly all repetitions, contaminated test instances survived.

In classification, we can essentially observe no difference between the performance of LogitBoost and AUC-Boosting, both for the identification of contaminated training and test instances, which is not surprising, but the difference between the performance of Logit regression and \texttt{glmrob} is also not considerable. The behaviour that the ranking error decreases both for the identification of the contaminated training and test instances decreases with increasing $r_{val}$ also occurs in this setting. The empirical BDP is always very close to 1.

\subsubsection{Contaminated training and test data, post-trimming} \label{sec:valconttraintrim}

Finally, we repeat the steps as in Sec. \ref{sec:valcont}, but we first trim the training data according to the LOO strategy. Then, the trimming of the test instances is based on a model trained on the trimmed training data.

The figures can be found in Sec. \ref{app:valconttraintrim}. We do not repeat the figures for the ranking errors on the training data, as there is no difference to the technique in Sec. \ref{sec:valcont}. However, as the training data are trimmed both according to a non-robust and to a robust model, we depict the results for each of both settings.

The results are similar as in Sec. \ref{sec:valcont}, and there is no clear tendency whether trimming the training data first improves the ranking performance. One can observe a slightly better performance when using a non-robust model for trimming the training data than a robust model. The difference can be seen most prominently in Fig. \ref{valconttraintrimp20bdp}, where the mean BDP of elicitability is nearly zero in the context of $Y$-contamination, $SNR=0.5$ and $r=0.25$, when using a non-robust model both for trimming the training data and for the identification of the most outlying test data, while the corresponding curve is significantly higher when using a robust model for trimming the training data. Another interesting example is the situation for $X$-contamination in Fig. \ref{valconttraintrimglmbdp} and Fig. \ref{valcontglmbdp}, where the robust model at least was able to identify all outlying test instances in around 40\% of the repetitions, while in Fig. \ref{valconttraintrimglmbdp}, trimming the training data in advance resulted in a success of this form in nearly no repetition.

\section{Discussion and conclusion}

We assumed for simplicity that the data with which we train and validate a statistical estimator are generated from a convex contamination model $F_t=(1-t)F+tG$, for the ideal distribution $F$. Other types of contamination can be found in \cite{rieder}. From the perspective of the ideal distribution, the goal is to infer a statistical property w.r.t. $F$, e.g., the mean. 

\textit{What can we learn from the simulation results?} We restricted our simulations to rather simple regression and classification settings, where the contamination just appears as an additive shift. Of course, we do not claim any form of optimality regarding our trimming approach, and the performance could be improved by using outlier detection procedures first on the data and by iterating multiple times between training and test data in the sense that one iteratively removes a certain fraction of instances until an $\alpha$-fraction is removed, instead of directly removing an $\alpha$-fraction in one iteration. Nevertheless, the results are discouraging. In a cross-validation setting, one can easily generate situations where a fold contains a large fraction of outlying instances. Such a fold is problematic for validation, even if one would use a robust aggregation procedure for the losses, as suggested in \cite{khan10}. The identification of such batches was often not better than random guessing, in particular, for high contamination radii. As for static test sets, the identification of the outlying test instances has a high success rate for some configurations, in other, it is nearly zero. Here, one has to take into account the simplicity of the contamination, which makes an identification rather easy, and the fact that we used an oracle trimming rate, which equals the true contamination radius, which is not known in practice. One could argue that using a large trimming rate, maybe even more than 0.5, could work in practice, assuming that the ideal and contaminating distribution have a sufficiently low overlap, however, one ultimately sacrifices a lot of (expensive) test data, without any guarantee that only instances from the ideal distribution have survived.

\textit{Which issues are implied to validation by contaminated distributions?} We have shown that even if the property is elicitable w.r.t. the family of distributions that includes $F$ but not $G$, using a loss function from the corresponding family of strictly consistent loss functions fails in eliciting the property w.r.t. $F$. As trimming is part of some successful robust estimation procedures, we aimed at removing outlying instances from the test data, ideally only those from $G$, since a perfect identification and trimming of those instances would be theoretically sound in the sense that elicitability w.r.t. $F$ is maintained. The main difference between trimming in training and testing is that during training, the model is not necessarily distorted if contaminating instances survive, which does not carry over to validation because the presence of realization from $G$ in the test data does not allow for eliciting the property of interest w.r.t. $F$, although the test loss is not necessarily distorted. If the test data are realized from $F_t$, then a strictly consistent loss function for the property of interest, e.g., the mean, such as the squared loss function, elicits $\E[Y_t]$ for $Y_t \sim F_t$, but not $\E[Y]$ for $Y \sim F$, i.e., $\E[Y]$ would appear to be a sub-optimal estimate.

\textit{Why would that be problematic?} In the simplest case where $G$ is contained in the same family of distributions as $F$, for which elicitability holds, one could take the perspective of robust estimation and argue that the population minimizer of a robust loss function such as a redescender, or the trimmed population minimizer of a non-robust loss function, is essentially a trimmed mean (w.r.t. $F_t$). Would it not be natural to include $G$ into validation? The answer should be ``no'', because it would even put robust estimation procedures in question. If some kind of trimmed mean would be used for validation, it would in fact just be a construction in order to let the robust procedure, that searches for the minimizer of this criterion, appear to be optimal. If, in the simple case where $G$ is sufficiently regular to allow for elicitability, one would just use a standard strictly consistent loss function and indeed approximate $\E_{F_t}[L(Y_t,x)]$ for $Y_t \sim F_t$ and some prediction $x$, why would one have applied a robust estimation procedure in the first place?

\textit{Could one just change the aggregation procedure of the test losses?} In robust estimation, robust aggregation procedures such as medians or general trimmed means are popular. We already pointed out that that the required moment inequality $\E_F[L(Y,t)]<\E_F[L(Y,x)]$, for the quantity of interest, $t \in T(F)$, and any prediction $x \notin T(F)$, would not be transferable to ``robust'' transformations of the loss function the form $h(L)$, which restrict the image space in order to circumvent the issue that the contaminated distribution $F_t$ may not allow for moments. Similarly, robustly aggregating the test losses, e.g., by a median as suggested in \cite{khan10}, would be in conflict with the established notion of elicitability. In fact, this idea would require a modified notion of elicitability, namely the property that \begin{equation} \label{medelic} med_F(L(Y,t))<med_F(L(Y,x)).  \end{equation} Unfortunately, inequalities in terms of expectation values generally do not carry over to median inequalities, necessitating to find new families of strictly consistent loss functions that satisfy an inequality of the form \ref{medelic}.  

\textit{How helpful is a modified notion of elicitability for robust aggregation?} Should one be able to find loss functions that satisfy inequalities such as \ref{medelic}, one would clearly remedy the problem that expectations may not exist if the contaminated distribution is not sufficiently regular. However, the problem that, on contaminated data, the corresponding inequality would consider medians or trimmed means of the contaminated distribution $F_t$ instead of the ideal distribution $F$ would still remain, so one can suspect that even for strictly consistent loss functions of that kind, the true quantity of interest w.r.t. $F$ would appear to be sub-optimal.

Therefore, the question how to properly evaluate and compare models in the presence of contaminated test data should be prioritized for future research.

\section*{Acknowledgements}

The author would like to thank Prof. P. Ruckdeschel for helpful discussions. Of course, the author is solely responsible for any error.

\renewcommand\refname{References}
\bibliography{Biblio}
\bibliographystyle{abbrvnat}
\setcitestyle{authoryear,open={((},close={))}}

\appendix

\section{Trimming resamples: Figures} \label{app:trimresample}

\subsection{Contaminated training data, clean test data: Figures}  \label{app:valclean}

In the following figures, the left part always corresponds to an SNR of 5, the middle part of an SNR of 2, and the right part of an SNR of 0.5, respectively, for regression, and to $\mu=8$, $\mu=3$, and $\mu=0.5$, respectively, for classification.

\subsubsection{$p=20$, regression, loss-based}

\begin{figure}[H]
\begin{center}
\includegraphics[width=5cm,height=4cm]{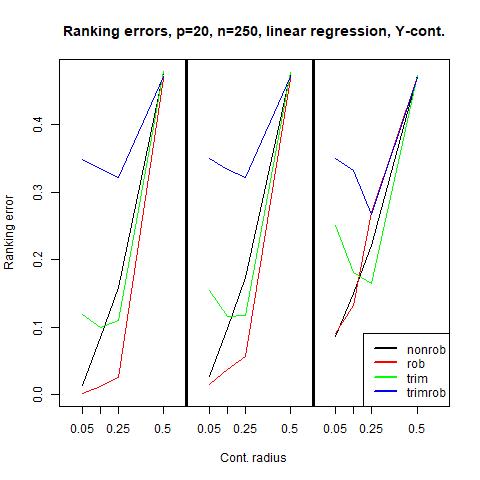} 
\includegraphics[width=5cm,height=4cm]{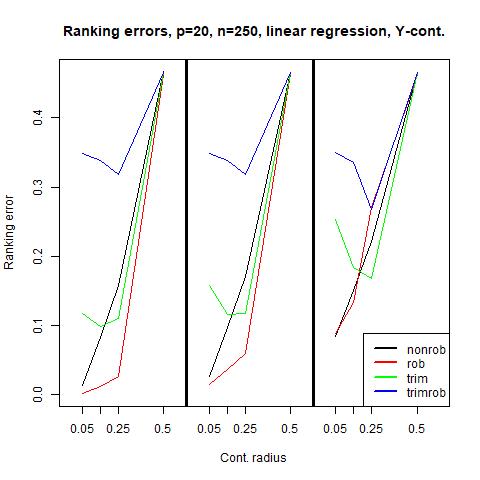} \\
\includegraphics[width=5cm,height=4cm]{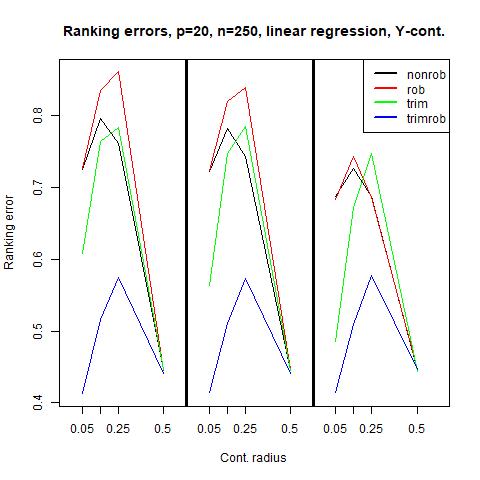} 
\includegraphics[width=5cm,height=4cm]{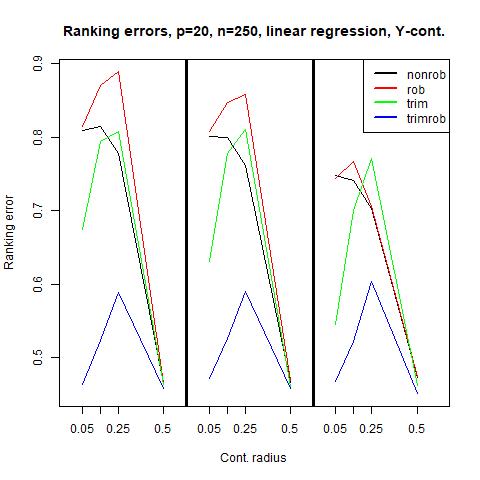} 
\caption{\tiny{Hard ranking errors in the context of $Y$-contamination for randomized cross validation with 10 (upper left) and 100 batches (upper right) and for $k$-fold cross validation with $K=10$ (bottom left) and $K=5$ (bottom right). The black lines correspond to the ordering of the instances according to non-trimmed losses of non-robust regression, the green lines to trimmed losses of non-robust regression, the red lines to non-trimmed losses of robust regression, and the blue lines to trimmed losses of robust regression. All losses are training losses.}} \label{valcleanp20Ycont}
\end{center}
\end{figure}

\begin{figure}[H]
\begin{center}
\includegraphics[width=5cm,height=4cm]{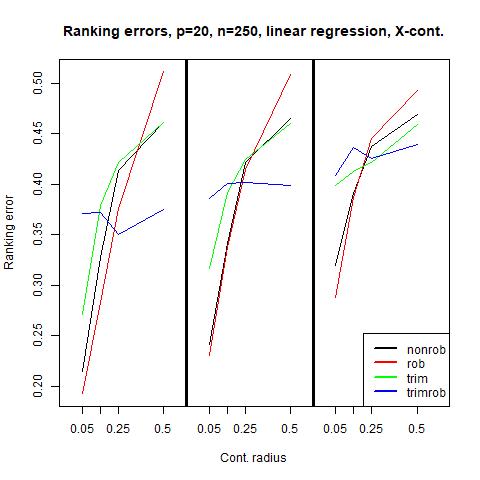} 
\includegraphics[width=5cm,height=4cm]{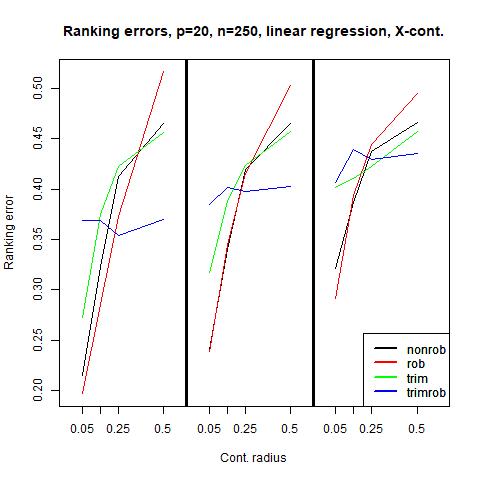} \\
\includegraphics[width=5cm,height=4cm]{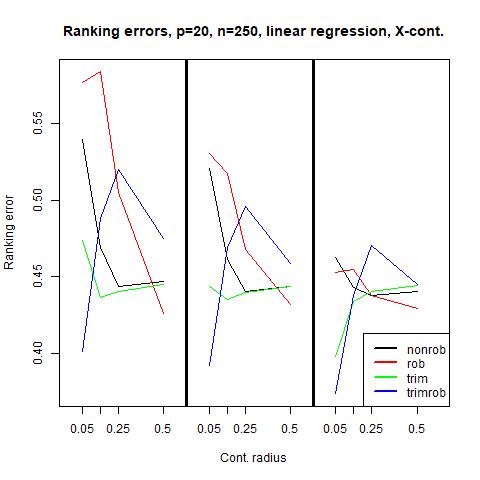} 
\includegraphics[width=5cm,height=4cm]{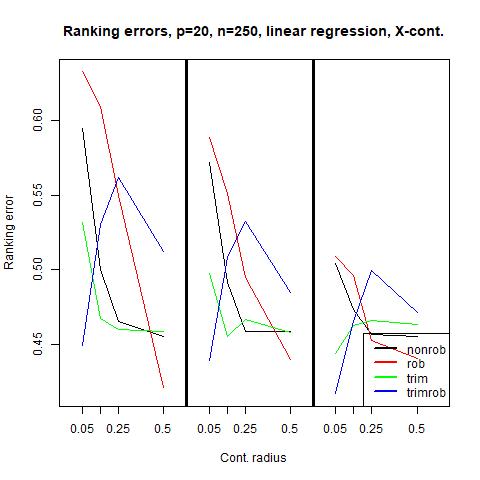} 
\caption{\tiny{Hard ranking errors in the context of $X$-contamination for randomized cross validation with 10 (upper left) and 100 batches (upper right) and for $k$-fold cross validation with $K=10$ (bottom left) and $K=5$ (bottom right). The black lines correspond to the ordering of the instances according to non-trimmed losses of non-robust regression, the green lines to trimmed losses of non-robust regression, the red lines to non-trimmed losses of robust regression, and the blue lines to trimmed losses of robust regression. All losses are training losses.}} \label{valcleanp20Xcont}
\end{center}
\end{figure}

\begin{figure}[H]
\begin{center}
\includegraphics[width=5cm,height=4cm]{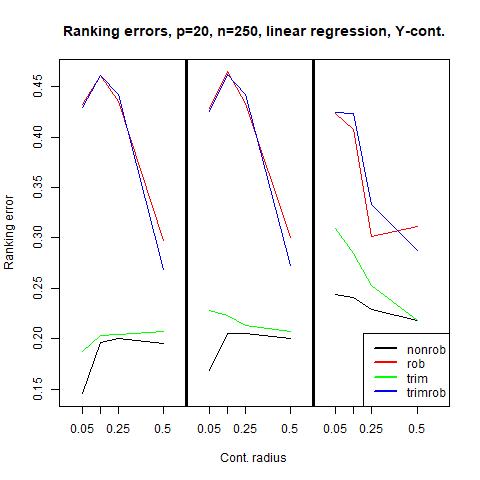} 
\includegraphics[width=5cm,height=4cm]{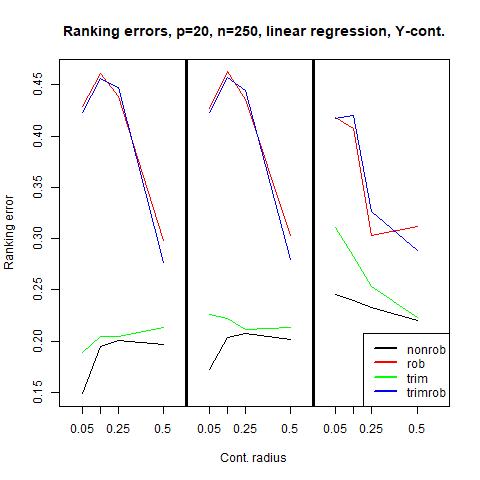} \\
\includegraphics[width=5cm,height=4cm]{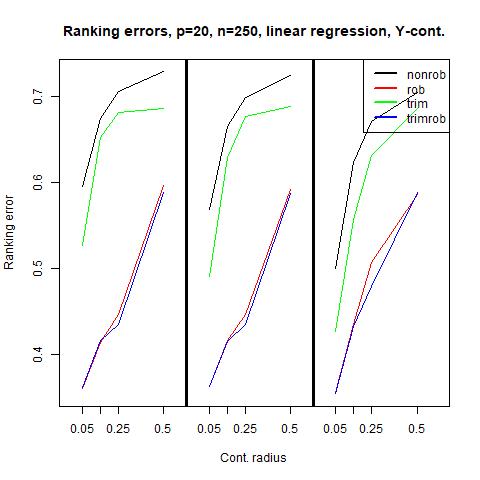} 
\includegraphics[width=5cm,height=4cm]{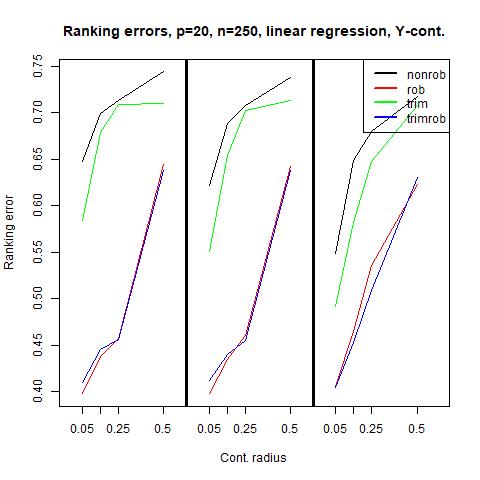} 
\caption{\tiny{Hard ranking errors in the context of $Y$-contamination for randomized cross validation with 10 (upper left) and 100 batches (upper right) and for $k$-fold cross validation with $K=10$ (bottom left) and $K=5$ (bottom right). The black lines correspond to the ordering of the instances according to non-trimmed losses of non-robust regression, the green lines to trimmed losses of non-robust regression, the red lines to non-trimmed losses of robust regression, and the blue lines to trimmed losses of robust regression. All losses are test losses.}} \label{valcleanp20Ycontvalloss}
\end{center}
\end{figure}

\begin{figure}[H]
\begin{center}
\includegraphics[width=5cm,height=4cm]{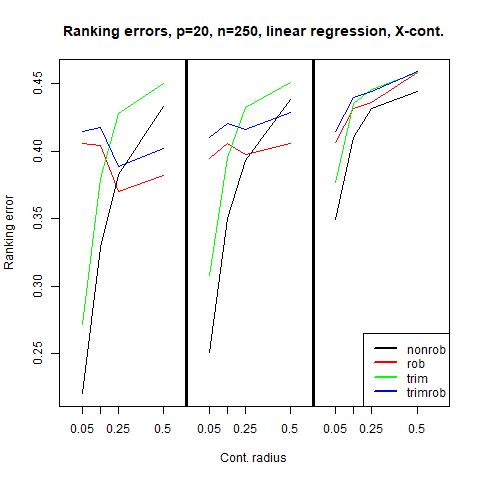} 
\includegraphics[width=5cm,height=4cm]{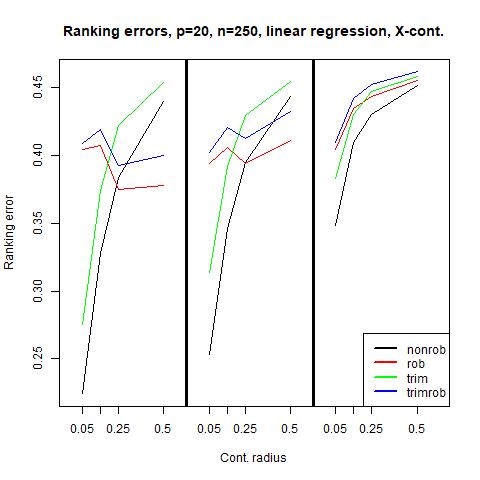} \\
\includegraphics[width=5cm,height=4cm]{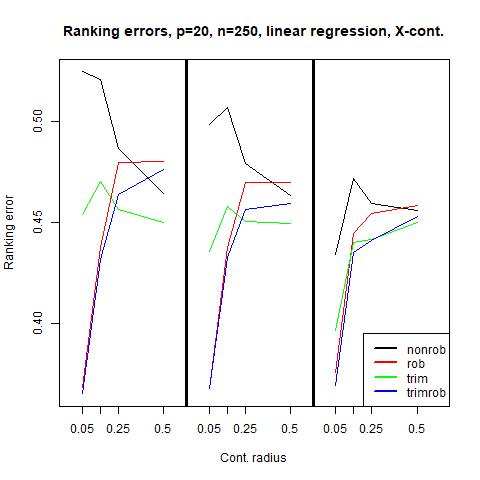} 
\includegraphics[width=5cm,height=4cm]{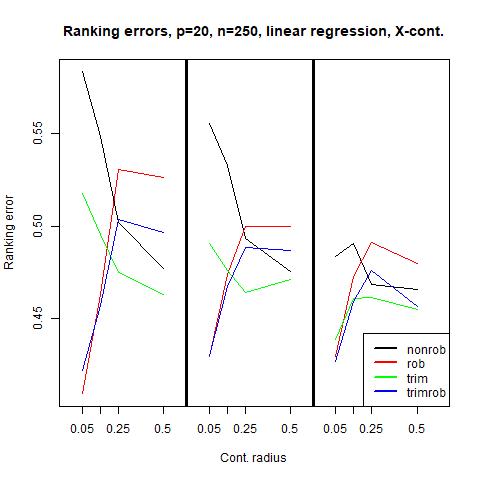} 
\caption{\tiny{Hard ranking errors in the context of $X$-contamination for randomized cross validation with 10 (upper left) and 100 batches (upper right) and for $k$-fold cross validation with $K=10$ (bottom left) and $K=5$ (bottom right). The black lines correspond to the ordering of the instances according to non-trimmed losses of non-robust regression, the green lines to trimmed losses of non-robust regression, the red lines to non-trimmed losses of robust regression, and the blue lines to trimmed losses of robust regression. All losses are test losses.}} \label{valcleanp20Xcontvalloss}
\end{center}
\end{figure}

\subsubsection{$p=250$, regression, loss-based}

\begin{figure}[H]
\begin{center}
\includegraphics[width=5cm,height=4cm]{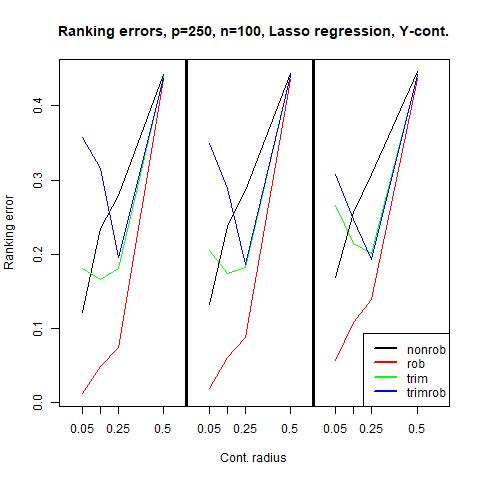} 
\includegraphics[width=5cm,height=4cm]{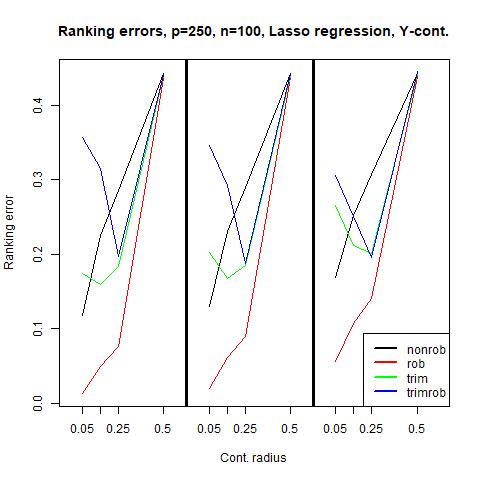} \\
\includegraphics[width=5cm,height=4cm]{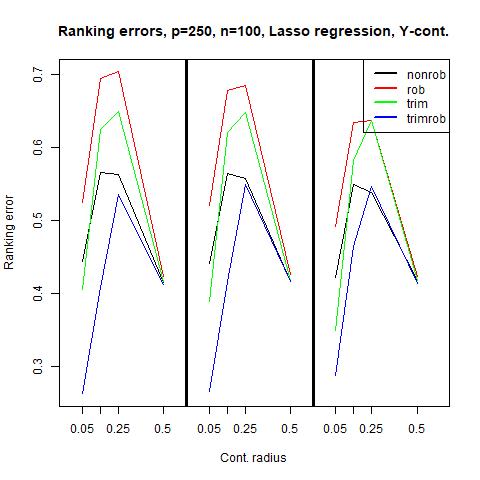} 
\includegraphics[width=5cm,height=4cm]{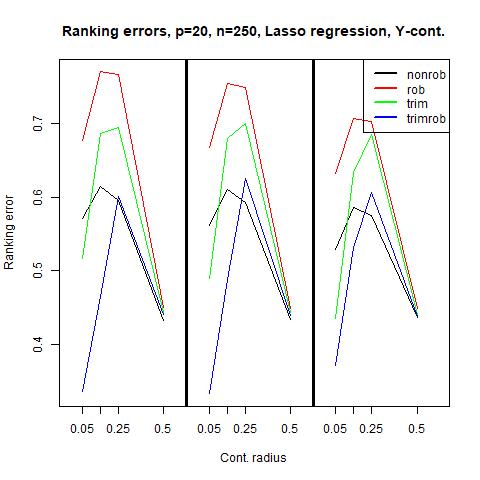} 
\caption{\tiny{Hard ranking errors in the context of $Y$-contamination for randomized cross validation with 10 (upper left) and 100 batches (upper right) and for $k$-fold cross validation with $K=10$ (bottom left) and $K=5$ (bottom right). The black lines correspond to the ordering of the instances according to non-trimmed losses of non-robust regression, the green lines to trimmed losses of non-robust regression, the red lines to non-trimmed losses of robust regression, and the blue lines to trimmed losses of robust regression. All losses are training losses.}} \label{valcleanp250Ycont}
\end{center}
\end{figure}

\begin{figure}[H]
\begin{center}
\includegraphics[width=5cm,height=4cm]{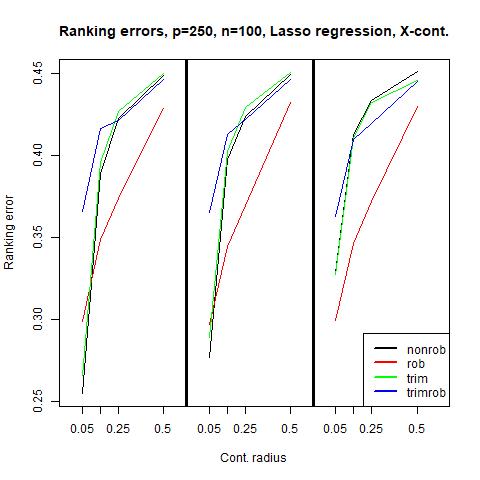} 
\includegraphics[width=5cm,height=4cm]{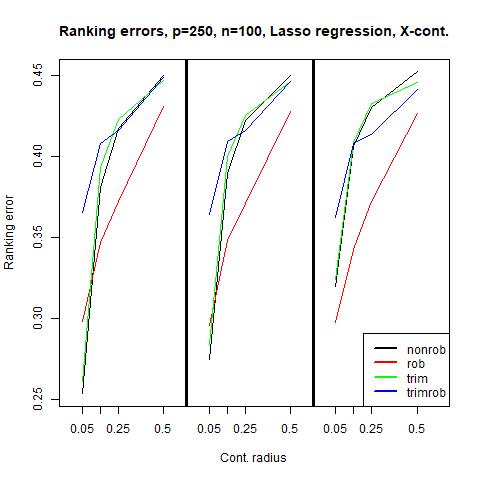} \\
\includegraphics[width=5cm,height=4cm]{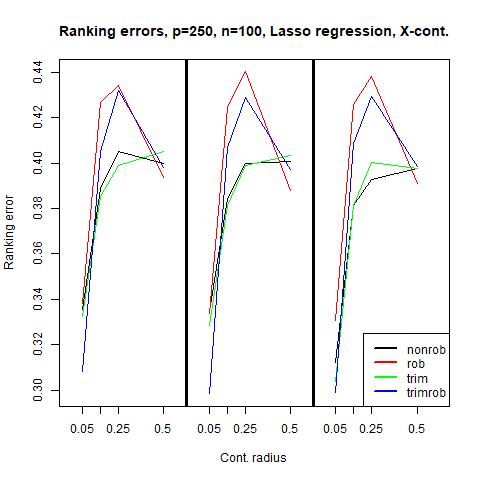} 
\includegraphics[width=5cm,height=4cm]{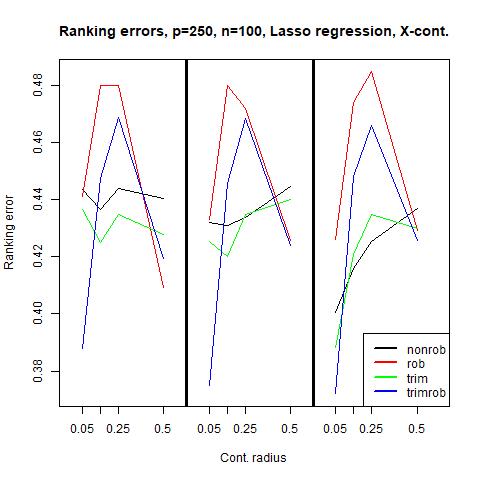} 
\caption{\tiny{Hard ranking errors in the context of $X$-contamination for randomized cross validation with 10 (upper left) and 100 batches (upper right) and for $k$-fold cross validation with $K=10$ (bottom left) and $K=5$ (bottom right). The black lines correspond to the ordering of the instances according to non-trimmed losses of non-robust regression, the green lines to trimmed losses of non-robust regression, the red lines to non-trimmed losses of robust regression, and the blue lines to trimmed losses of robust regression. All losses are training losses.}} \label{valcleanp250Xcont}
\end{center}
\end{figure}

\begin{figure}[H]
\begin{center}
\includegraphics[width=5cm,height=4cm]{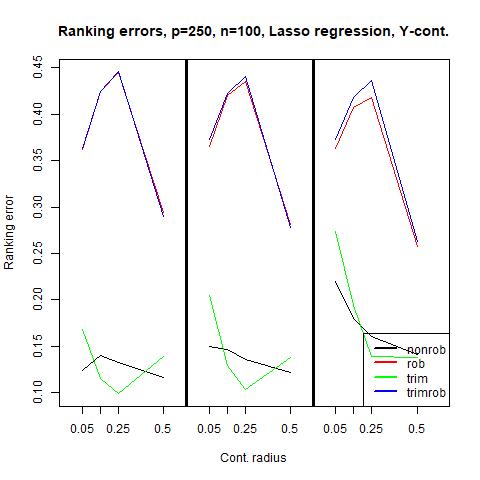} 
\includegraphics[width=5cm,height=4cm]{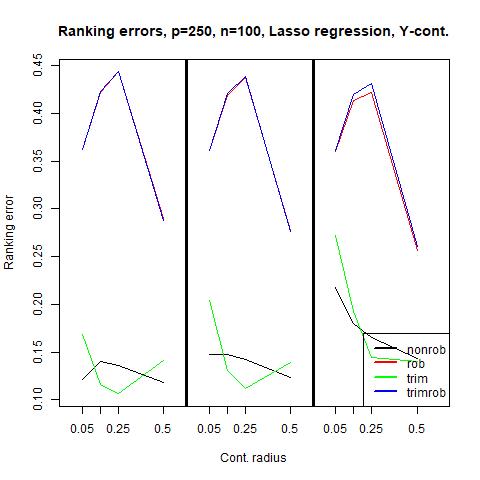} \\
\includegraphics[width=5cm,height=4cm]{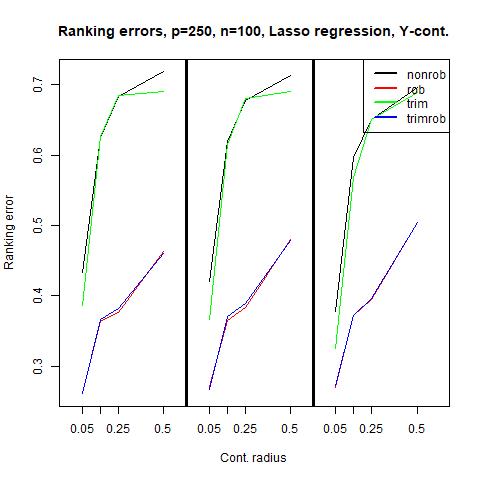} 
\includegraphics[width=5cm,height=4cm]{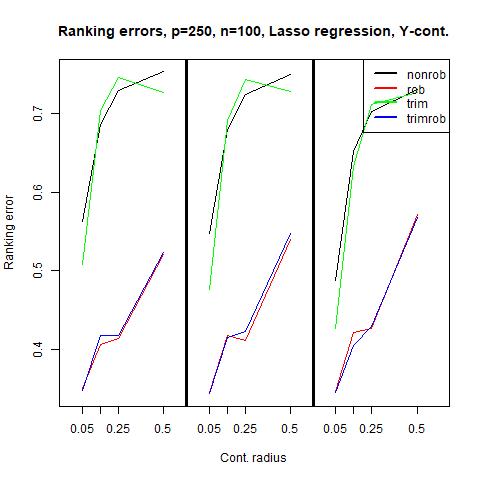} 
\caption{\tiny{Hard ranking errors in the context of $Y$-contamination for randomized cross validation with 10 (upper left) and 100 batches (upper right) and for $k$-fold cross validation with $K=10$ (bottom left) and $K=5$ (bottom right). The black lines correspond to the ordering of the instances according to non-trimmed losses of non-robust regression, the green lines to trimmed losses of non-robust regression, the red lines to non-trimmed losses of robust regression, and the blue lines to trimmed losses of robust regression. All losses are test losses.}} \label{valcleanp250Ycontvalloss}
\end{center}
\end{figure}

\begin{figure}[H]
\begin{center}
\includegraphics[width=5cm,height=4cm]{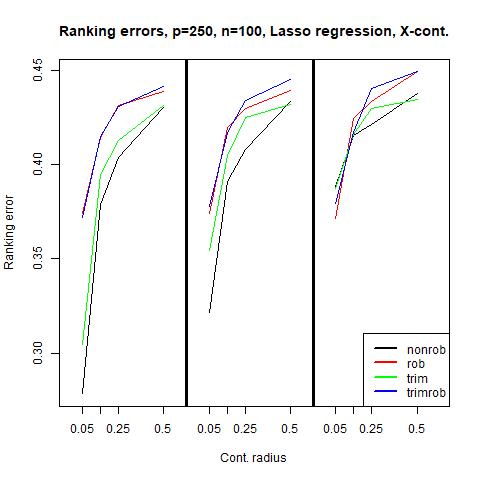} 
\includegraphics[width=5cm,height=4cm]{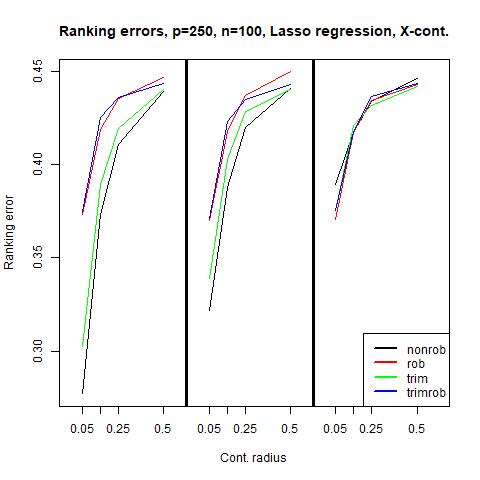} \\
\includegraphics[width=5cm,height=4cm]{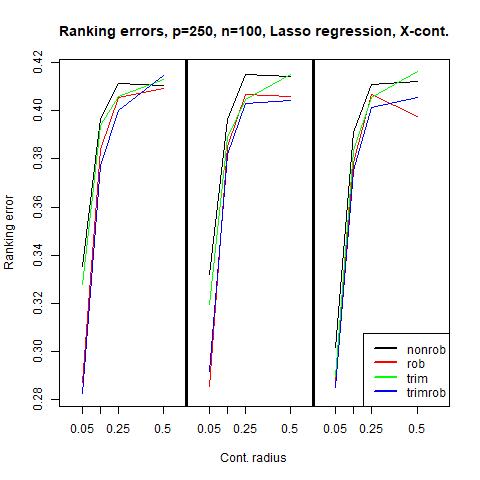} 
\includegraphics[width=5cm,height=4cm]{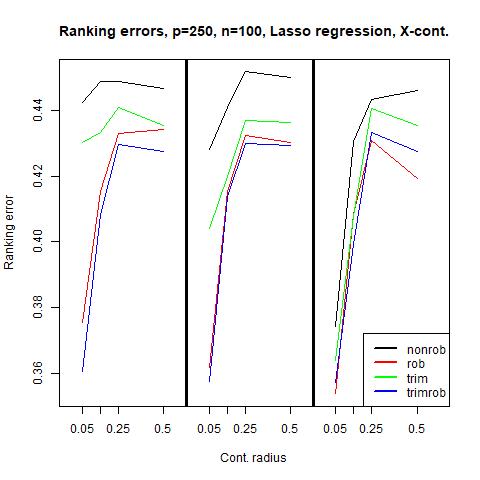} 
\caption{\tiny{Hard ranking errors in the context of $X$-contamination for randomized cross validation with 10 (upper left) and 100 batches (upper right) and for $k$-fold cross validation with $K=10$ (bottom left) and $K=5$ (bottom right). The black lines correspond to the ordering of the instances according to non-trimmed losses of non-robust regression, the green lines to trimmed losses of non-robust regression, the red lines to non-trimmed losses of robust regression, and the blue lines to trimmed losses of robust regression. All losses are test losses.}} \label{valcleanp250Xcontvalloss}
\end{center}
\end{figure}

\subsubsection{$p=500$, regression, loss-based}

\begin{figure}[H]
\begin{center}
\includegraphics[width=5cm,height=4cm]{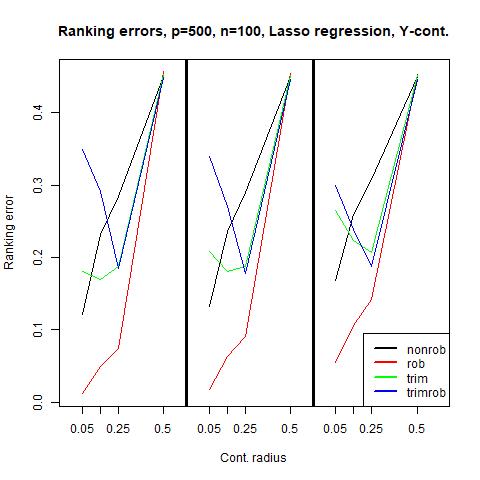} 
\includegraphics[width=5cm,height=4cm]{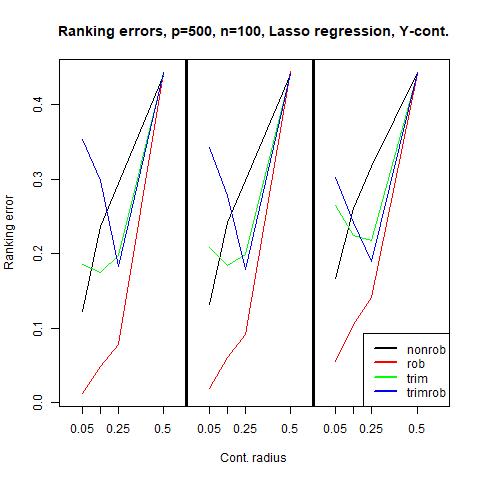} \\
\includegraphics[width=5cm,height=4cm]{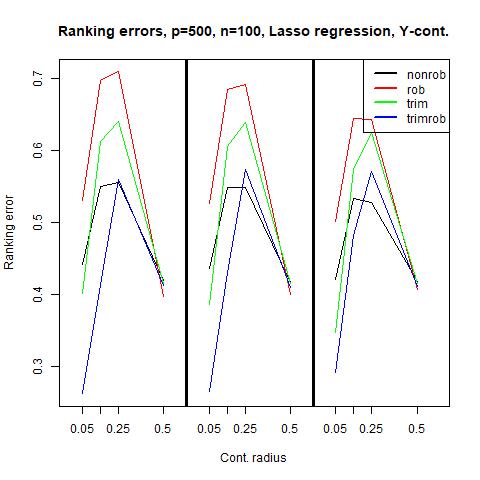} 
\includegraphics[width=5cm,height=4cm]{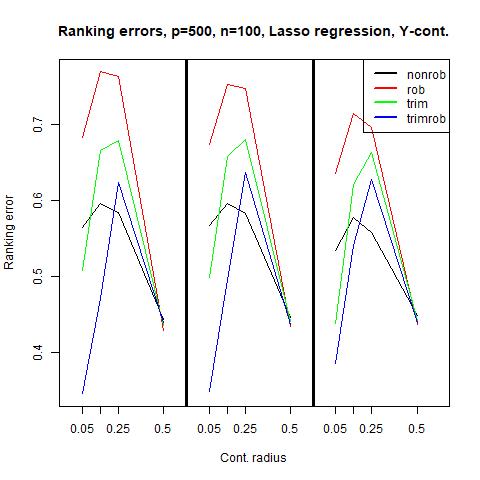} 
\caption{\tiny{Hard ranking errors in the context of $Y$-contamination for randomized cross validation with 10 (upper left) and 100 batches (upper right) and for $k$-fold cross validation with $K=10$ (bottom left) and $K=5$ (bottom right). The black lines correspond to the ordering of the instances according to non-trimmed losses of non-robust regression, the green lines to trimmed losses of non-robust regression, the red lines to non-trimmed losses of robust regression, and the blue lines to trimmed losses of robust regression. All losses are training losses.}} \label{valcleanp500Ycont}
\end{center}
\end{figure}

\begin{figure}[H]
\begin{center}
\includegraphics[width=5cm,height=4cm]{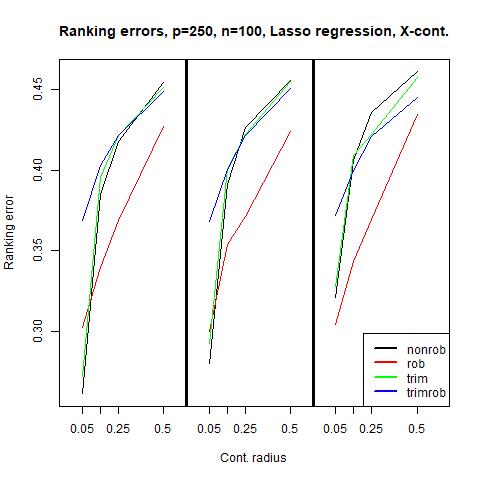} 
\includegraphics[width=5cm,height=4cm]{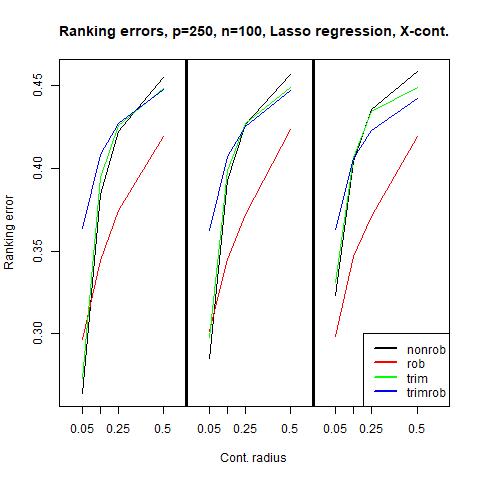} \\
\includegraphics[width=5cm,height=4cm]{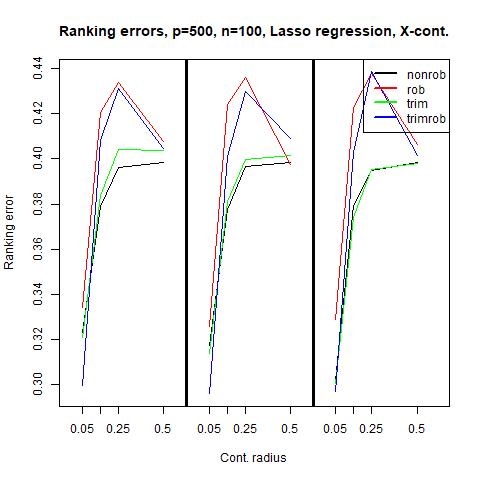} 
\includegraphics[width=5cm,height=4cm]{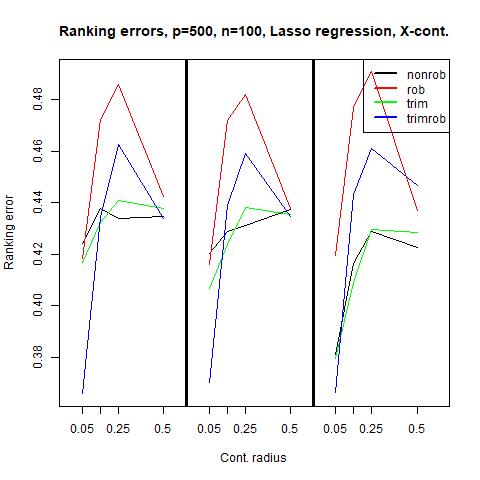} 
\caption{\tiny{Hard ranking errors in the context of $X$-contamination for randomized cross validation with 10 (upper left) and 100 batches (upper right) and for $k$-fold cross validation with $K=10$ (bottom left) and $K=5$ (bottom right). The black lines correspond to the ordering of the instances according to non-trimmed losses of non-robust regression, the green lines to trimmed losses of non-robust regression, the red lines to non-trimmed losses of robust regression, and the blue lines to trimmed losses of robust regression. All losses are training losses.}} \label{valcleanp500Xcont}
\end{center}
\end{figure}

\begin{figure}[H]
\begin{center}
\includegraphics[width=5cm,height=4cm]{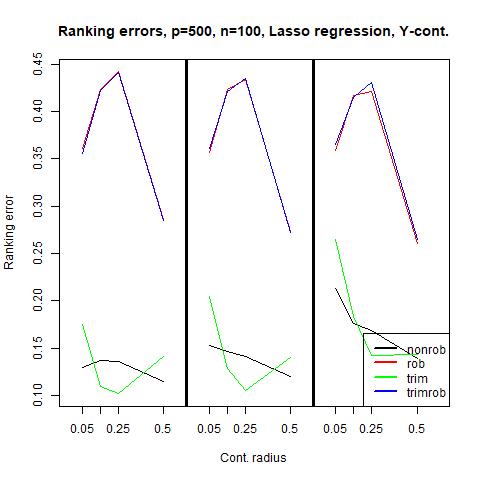} 
\includegraphics[width=5cm,height=4cm]{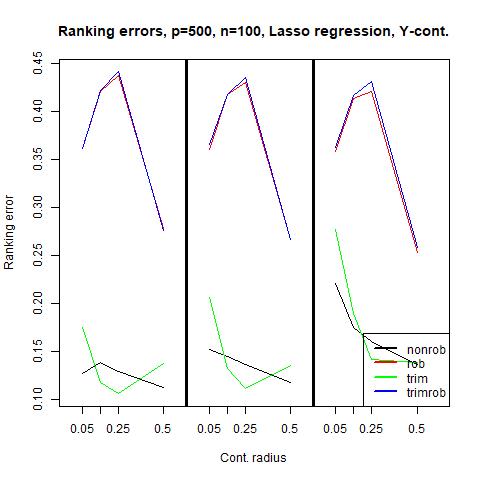} \\
\includegraphics[width=5cm,height=4cm]{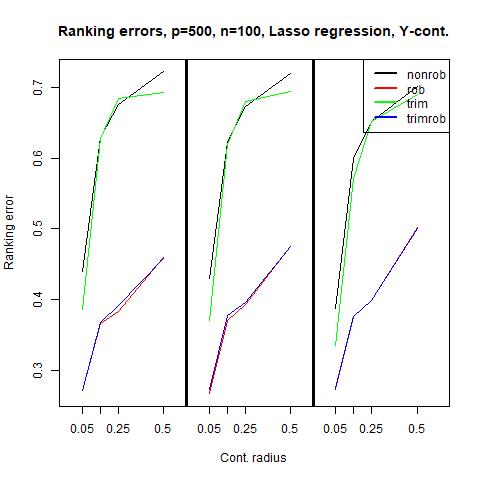} 
\includegraphics[width=5cm,height=4cm]{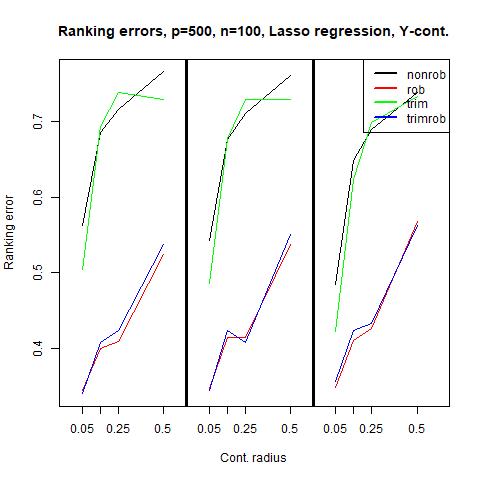} 
\caption{\tiny{Hard ranking errors in the context of $Y$-contamination for randomized cross validation with 10 (upper left) and 100 batches (upper right) and for $k$-fold cross validation with $K=10$ (bottom left) and $K=5$ (bottom right). The black lines correspond to the ordering of the instances according to non-trimmed losses of non-robust regression, the green lines to trimmed losses of non-robust regression, the red lines to non-trimmed losses of robust regression, and the blue lines to trimmed losses of robust regression. All losses are test losses.}} \label{valcleanp500Ycontvalloss}
\end{center}
\end{figure}

\begin{figure}[H]
\begin{center}
\includegraphics[width=5cm,height=4cm]{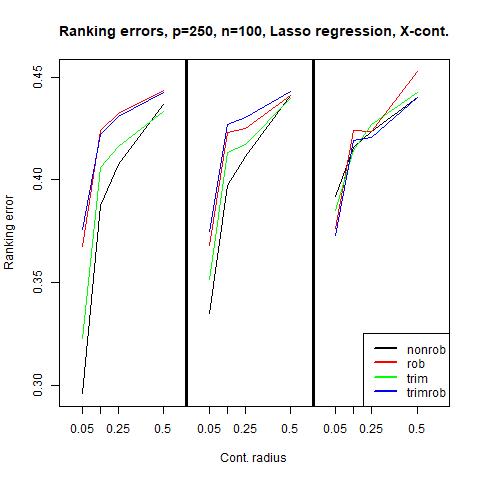} 
\includegraphics[width=5cm,height=4cm]{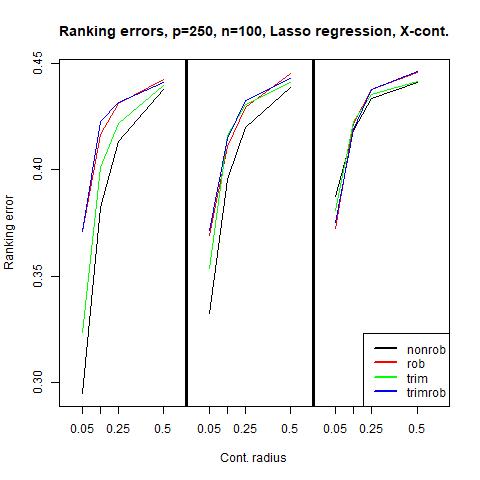} \\
\includegraphics[width=5cm,height=4cm]{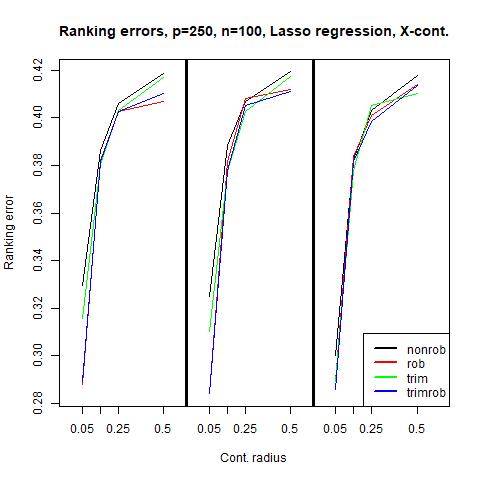} 
\includegraphics[width=5cm,height=4cm]{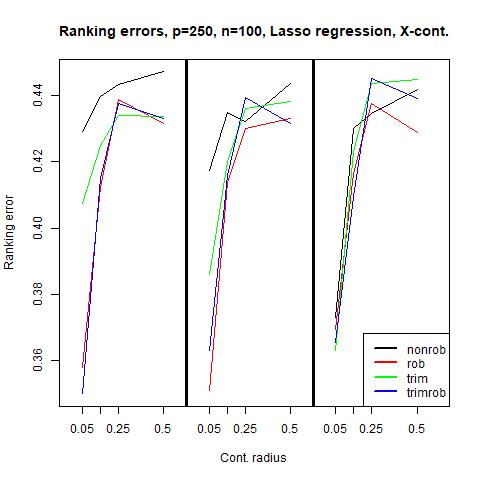} 
\caption{\tiny{Hard ranking errors in the context of $X$-contamination for randomized cross validation with 10 (upper left) and 100 batches (upper right) and for $k$-fold cross validation with $K=10$ (bottom left) and $K=5$ (bottom right). The black lines correspond to the ordering of the instances according to non-trimmed losses of non-robust regression, the green lines to trimmed losses of non-robust regression, the red lines to non-trimmed losses of robust regression, and the blue lines to trimmed losses of robust regression. All losses are test losses.}} \label{valcleanp500Xcontvalloss}
\end{center}
\end{figure}

\subsubsection{Regression, coefficient-based}

\begin{figure}[H]
\begin{center}
\includegraphics[width=5cm,height=4cm]{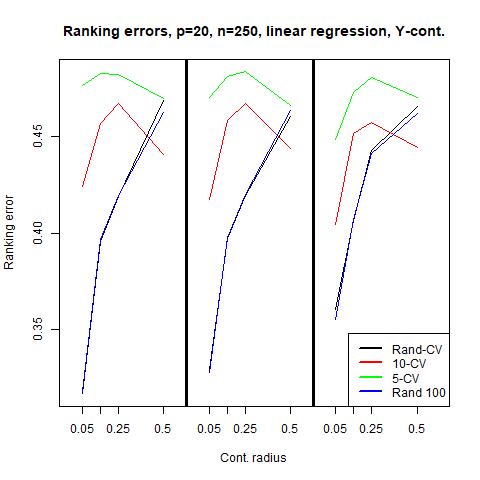} 
\includegraphics[width=5cm,height=4cm]{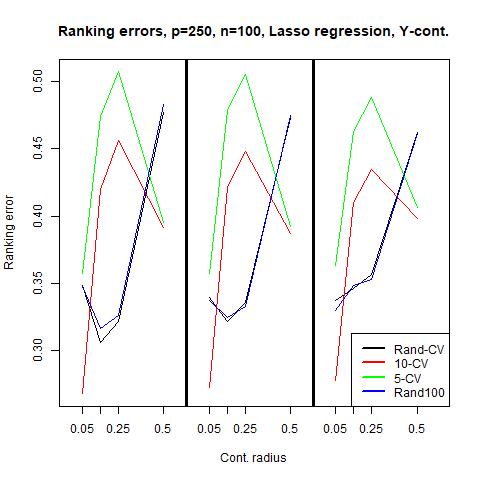} 
\includegraphics[width=5cm,height=4cm]{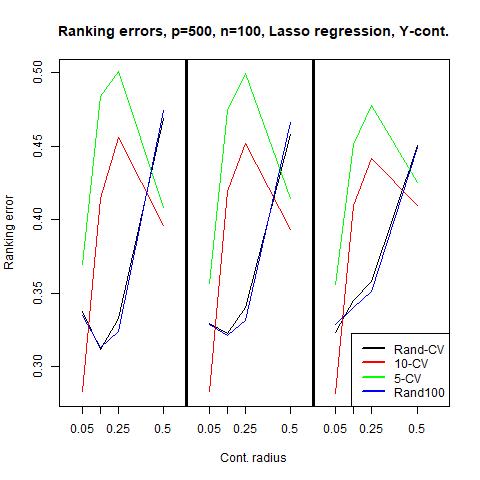} \\
\includegraphics[width=5cm,height=4cm]{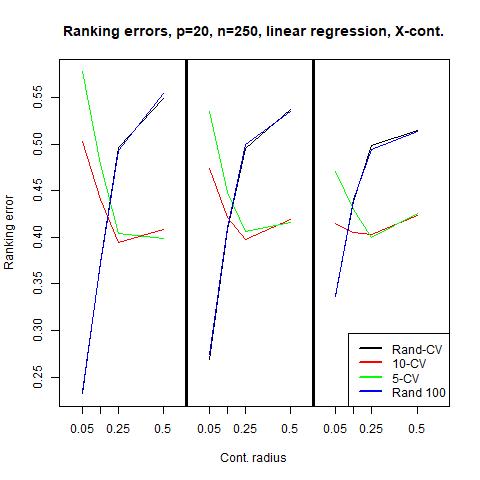} 
\includegraphics[width=5cm,height=4cm]{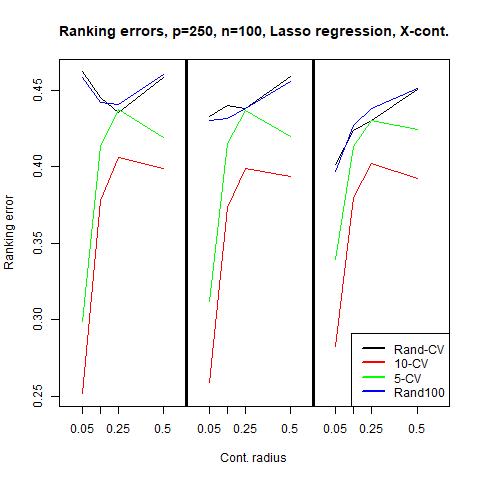} 
\includegraphics[width=5cm,height=4cm]{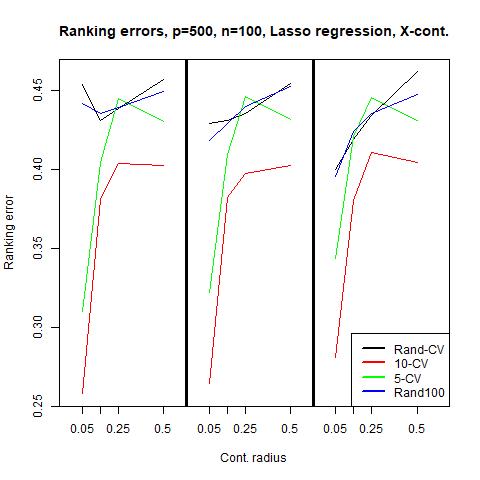} 
\caption{\tiny{Hard ranking errors in the context of $Y$-contamination (top row) and $X$-contamination (bottom row), respectively, for randomized cross validation for $p=20$ (left), $p=250$ (middle), and $p=500$ (right) according to the differences in the coefficient vectors of non-robust and robust regression. The black lines correspond to the ordering of the instances according to the Euclidean distances of the coefficient vectors for randomized cross-validation with 10 batches, the blue lines for randomized cross validation with 100 batches, and the red and green lines for $k$-fold cross validation with $K=10$ and $K=5$, respectively.}} \label{valcleancoefdiff}
\end{center}
\end{figure}

\subsubsection{$p=20$, classification, loss-based}

\begin{figure}[H]
\begin{center}
\includegraphics[width=5cm,height=4cm]{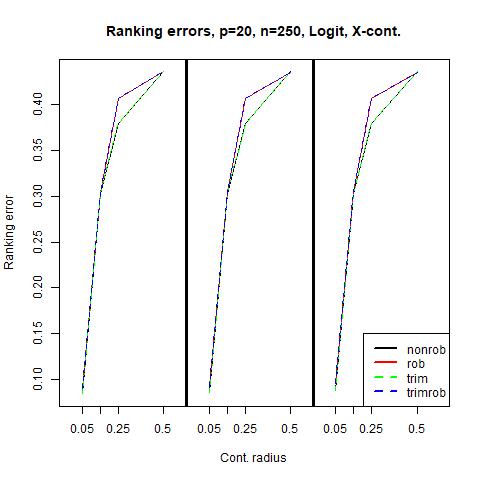} 
\includegraphics[width=5cm,height=4cm]{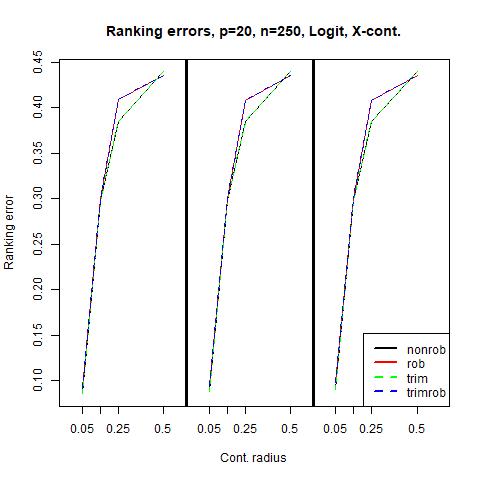} \\
\includegraphics[width=5cm,height=4cm]{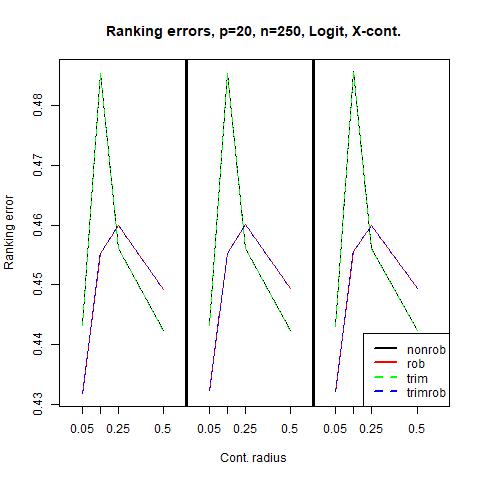} 
\includegraphics[width=5cm,height=4cm]{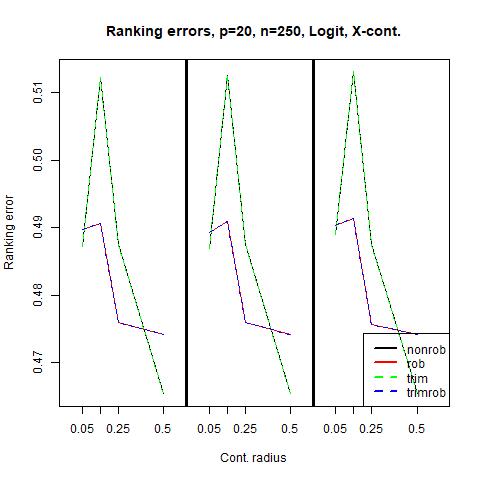} 
\caption{\tiny{Hard ranking errors in the context of $X$-contamination for randomized cross validation with 10 (upper left) and 100 batches (upper right) and for $k$-fold cross validation with $K=10$ (bottom left) and $K=5$ (bottom right). The black lines correspond to the ordering of the instances according to non-trimmed losses of Logit regression, the green lines to trimmed losses of Logit regression, the red lines to non-trimmed losses of robust binary classification, and the blue lines to trimmed losses of robust binary classification. All losses are training losses.}} \label{valcleanp20Xcontglm}
\end{center}
\end{figure}

\begin{figure}[H]
\begin{center}
\includegraphics[width=5cm,height=4cm]{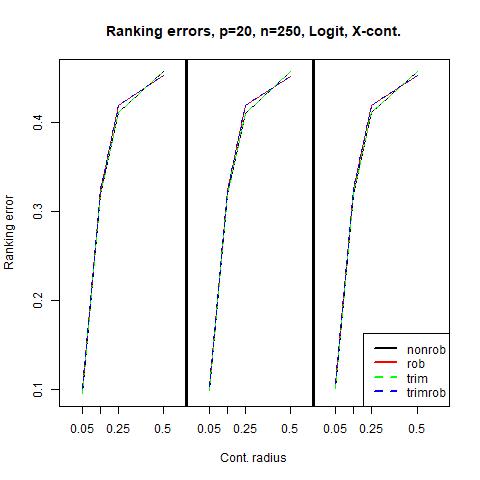} 
\includegraphics[width=5cm,height=4cm]{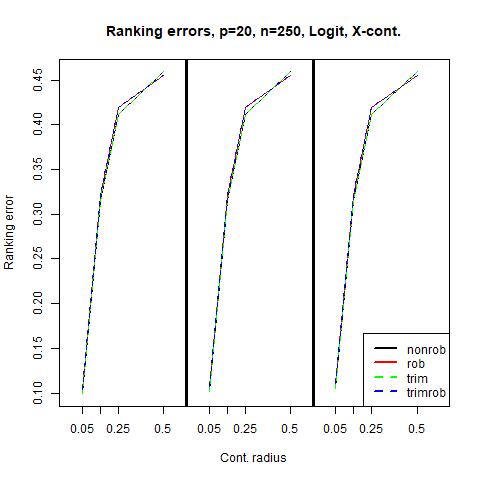} \\
\includegraphics[width=5cm,height=4cm]{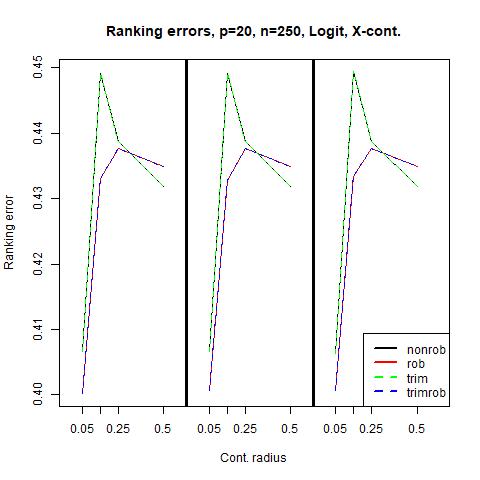} 
\includegraphics[width=5cm,height=4cm]{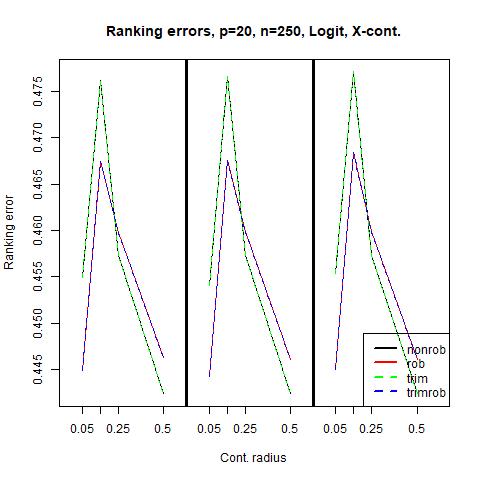} 
\caption{\tiny{Hard ranking errors in the context of $X$-contamination for randomized cross validation with 10 (upper left) and 100 batches (upper right) and for $k$-fold cross validation with $K=10$ (bottom left) and $K=5$ (bottom right). The black lines correspond to the ordering of the instances according to non-trimmed losses of Logit regression, the green lines to trimmed losses of Logit regression, the red lines to non-trimmed losses of robust binary classification, and the blue lines to trimmed losses of robust binary classification. All losses are test losses.}} \label{valcleanp20Xcontvallossglm}
\end{center}
\end{figure}

\subsubsection{$p=250$, classification, loss-based}

\begin{figure}[H]
\begin{center}
\includegraphics[width=5cm,height=4cm]{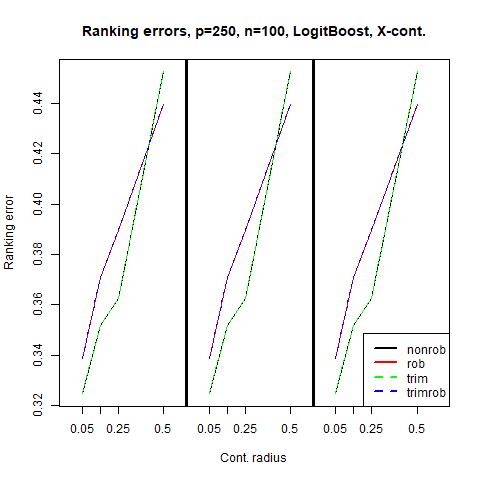} 
\includegraphics[width=5cm,height=4cm]{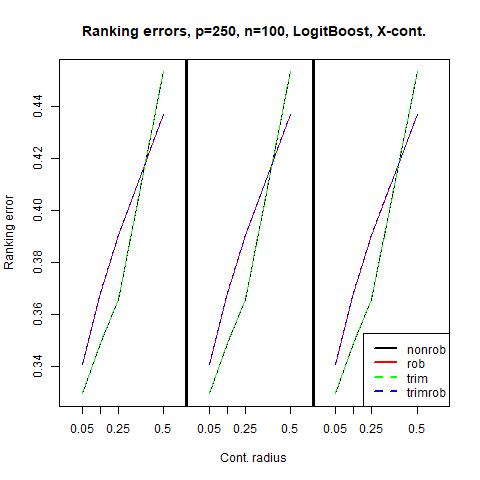} \\
\includegraphics[width=5cm,height=4cm]{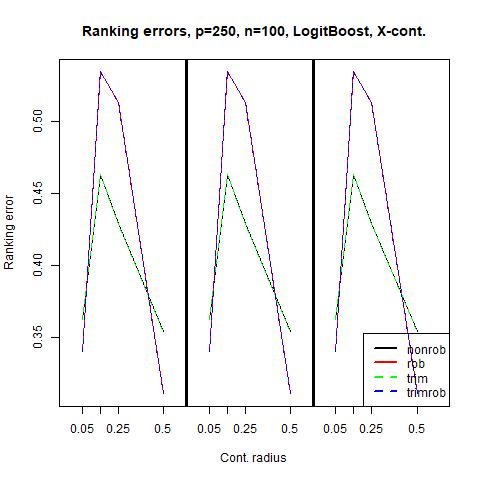} 
\includegraphics[width=5cm,height=4cm]{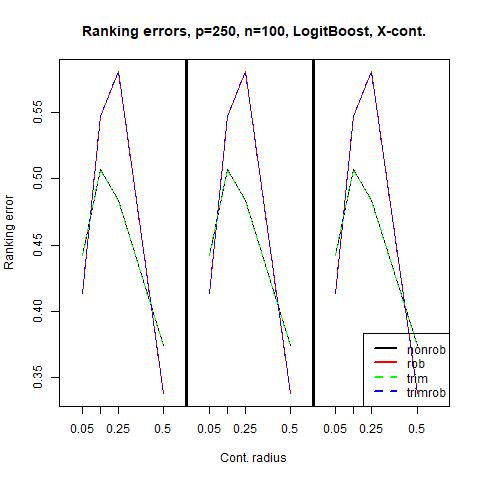} 
\caption{\tiny{Hard ranking errors in the context of $X$-contamination for randomized cross validation with 10 (upper left) and 100 batches (upper right) and for $k$-fold cross validation with $K=10$ (bottom left) and $K=5$ (bottom right). The black lines correspond to the ordering of the instances according to non-trimmed losses of LogitBoost, the green lines to trimmed losses of LogitBoost, the red lines to non-trimmed losses of AUC-Boosting, and the blue lines to trimmed losses of AUC-Boosting. All losses are training losses.}} \label{valcleanp250Xcontglm}
\end{center}
\end{figure}

\begin{figure}[H]
\begin{center}
\includegraphics[width=5cm,height=4cm]{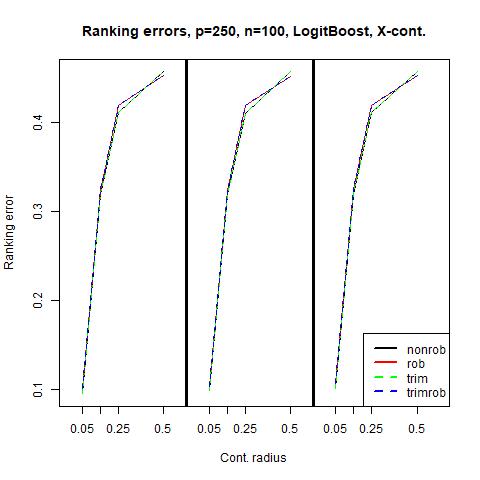} 
\includegraphics[width=5cm,height=4cm]{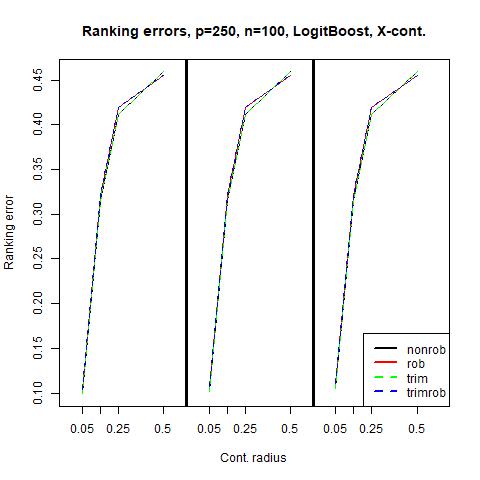} \\
\includegraphics[width=5cm,height=4cm]{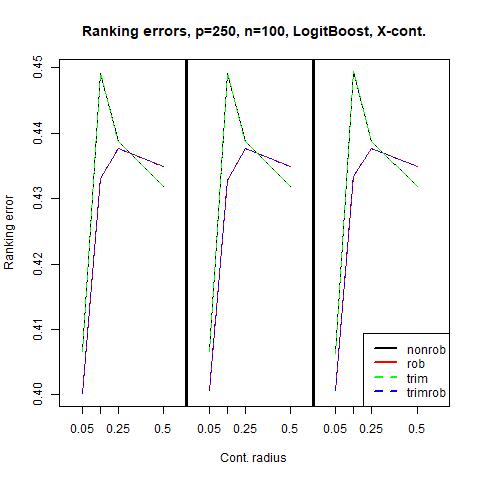} 
\includegraphics[width=5cm,height=4cm]{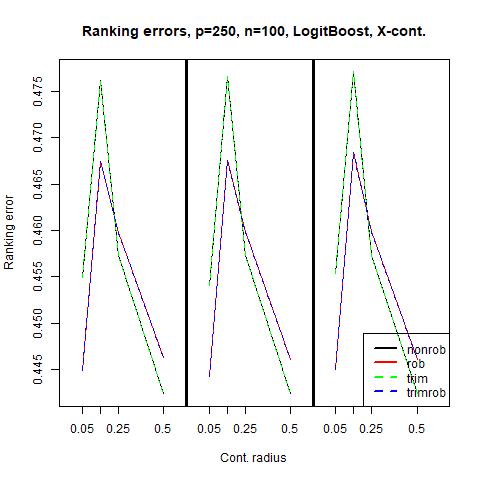} 
\caption{\tiny{Hard ranking errors in the context of $X$-contamination for randomized cross validation with 10 (upper left) and 100 batches (upper right) and for $k$-fold cross validation with $K=10$ (bottom left) and $K=5$ (bottom right). The black lines correspond to the ordering of the instances according to non-trimmed losses of LogitBoost, the green lines to trimmed losses of LogitBoost, the red lines to non-trimmed losses of AUC-Boosting, and the blue lines to trimmed losses of AUC-Boosting. All losses are test losses.} }\label{valcleanp250Xcontvallossglm}
\end{center}
\end{figure}

\subsubsection{$p=500$, classification, loss-based}

\begin{figure}[H]
\begin{center}
\includegraphics[width=5cm,height=4cm]{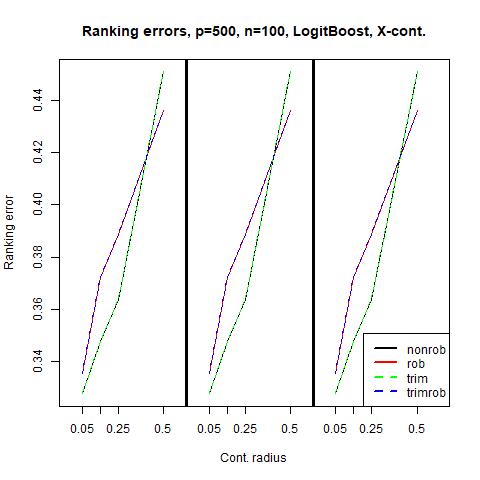} 
\includegraphics[width=5cm,height=4cm]{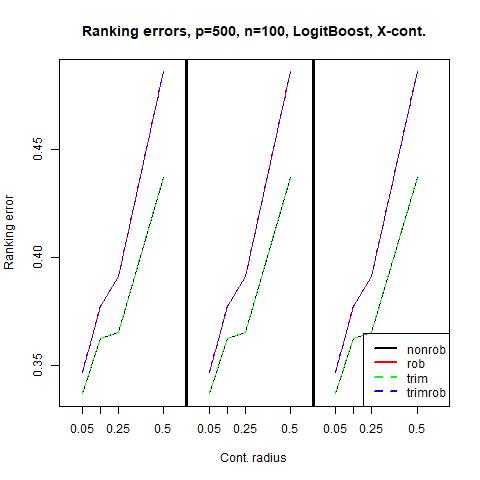} \\
\includegraphics[width=5cm,height=4cm]{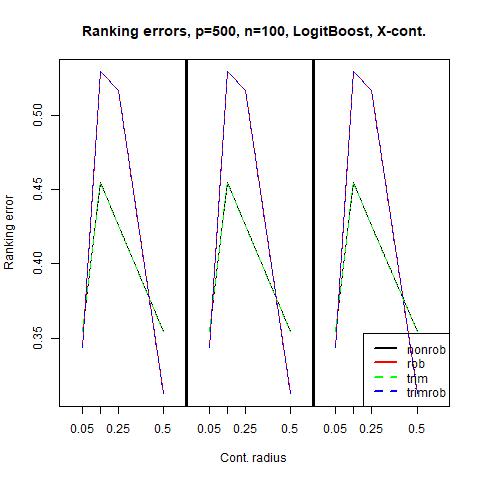} 
\includegraphics[width=5cm,height=4cm]{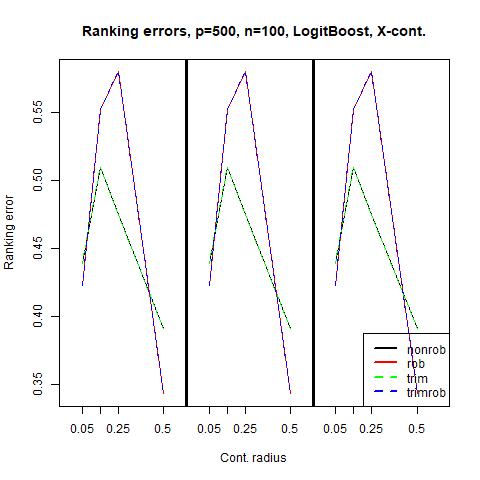} 
\caption{\tiny{Hard ranking errors in the context of $X$-contamination for randomized cross validation with 10 (upper left) and 100 batches (upper right) and for $k$-fold cross validation with $K=10$ (bottom left) and $K=5$ (bottom right). The black lines correspond to the ordering of the instances according to non-trimmed losses of LogitBoost, the green lines to trimmed losses of LogitBoost, the red lines to non-trimmed losses of AUC-Boosting, and the blue lines to trimmed losses of AUC-Boosting. All losses are training losses.}} \label{valcleanp500Xcontglm}
\end{center}
\end{figure}

\begin{figure}[H]
\begin{center}
\includegraphics[width=5cm,height=4cm]{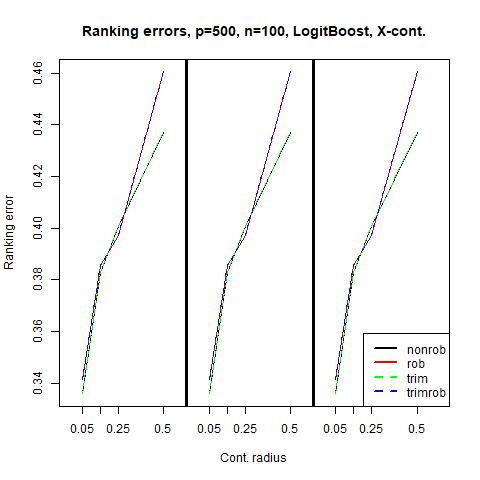} 
\includegraphics[width=5cm,height=4cm]{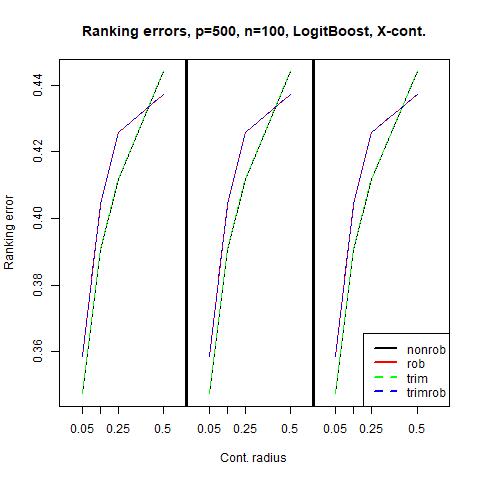} \\
\includegraphics[width=5cm,height=4cm]{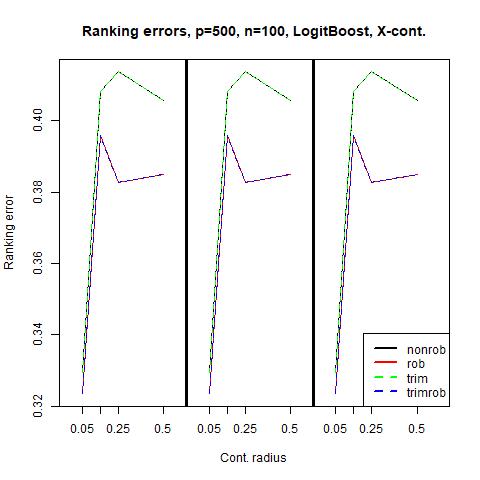} 
\includegraphics[width=5cm,height=4cm]{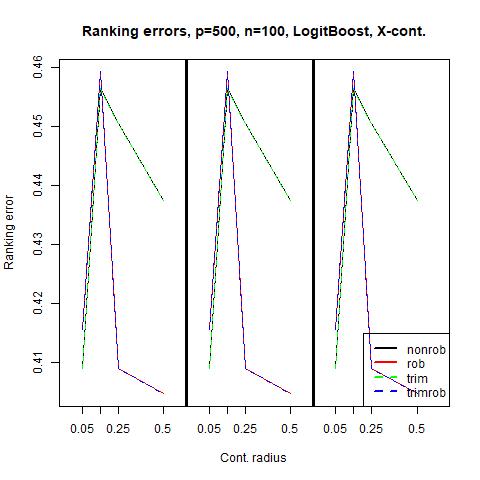} 
\caption{\tiny{Hard ranking errors in the context of $X$-contamination for randomized cross validation with 10 (upper left) and 100 batches (upper right) and for $k$-fold cross validation with $K=10$ (bottom left) and $K=5$ (bottom right). The black lines correspond to the ordering of the instances according to non-trimmed losses of LogitBoost, the green lines to trimmed losses of LogitBoost, the red lines to non-trimmed losses of AUC-Boosting, and the blue lines to trimmed losses of AUC-Boosting. All losses are test losses.}} \label{valcleanp500Xcontvallossglm}
\end{center}
\end{figure}

\subsubsection{Classification, coefficient-based}

\begin{figure}[H]
\begin{center}
\includegraphics[width=5cm,height=4cm]{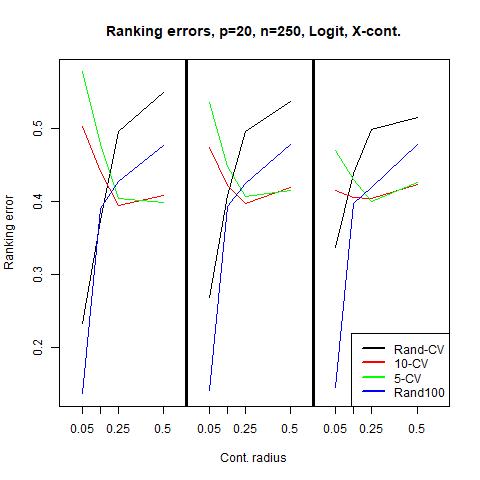} 
\includegraphics[width=5cm,height=4cm]{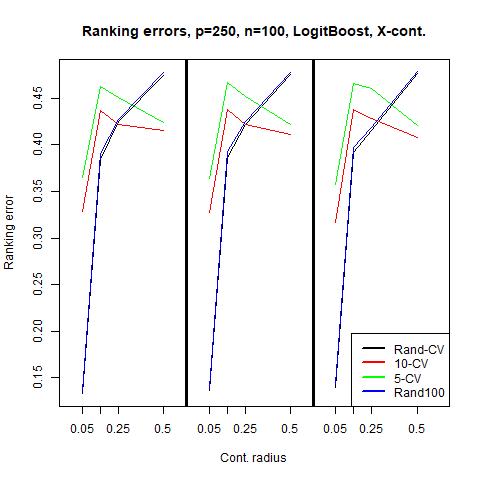} 
\includegraphics[width=5cm,height=4cm]{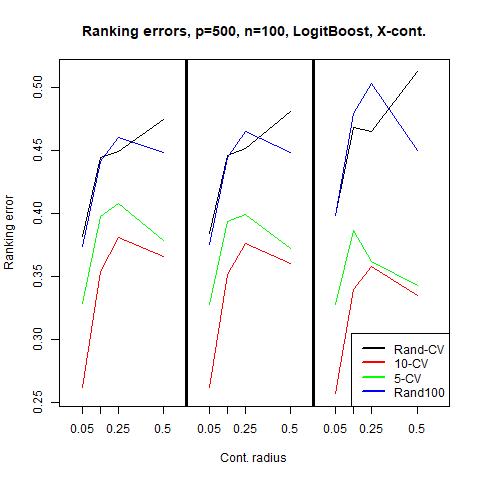} 
\caption{\tiny{Hard ranking errors in the context of $X$-contamination for randomized cross validation for $p=20$ (left), $p=250$ (middle), and $p=500$ (right) according to the differences in the coefficient vectors. The black lines correspond to the ordering of the instances according to the Euclidean distances of the coefficient vectors for randomized cross-validation with 10 batches, the blue lines for randomized cross validation with 100 batches, and the red and green lines for $k$-fold cross validation with $K=10$ and $K=5$, respectively.}} \label{valcleancoefdiffglm}
\end{center}
\end{figure}

\subsection{Contaminated training and test data, post-trimming}  \label{app:traintrim}

Here, due to the minor differences between the scenarios with an SNR of 5 and an SNR of 2, we restrict ourselves to an SNR of 5 (left part of the graphics), and an SNR of 0.5.

\subsubsection{$p=20$, regression, loss-based}

\begin{figure}[H]
\begin{center}
\includegraphics[width=5cm,height=4cm]{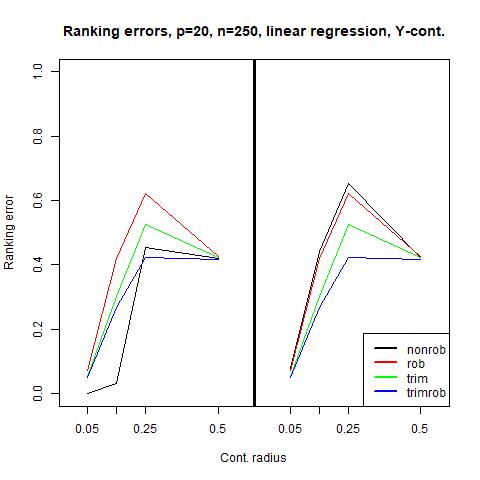}
\includegraphics[width=5cm,height=4cm]{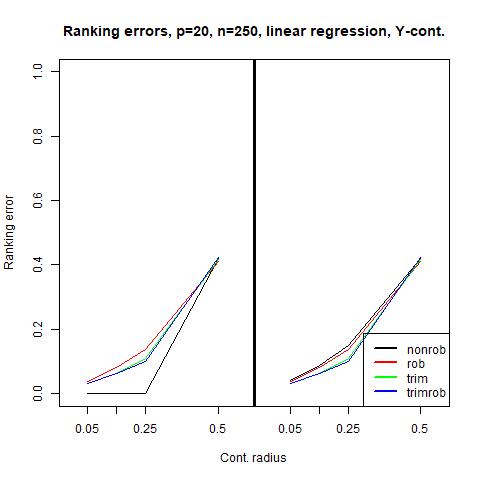} \\
\includegraphics[width=5cm,height=4cm]{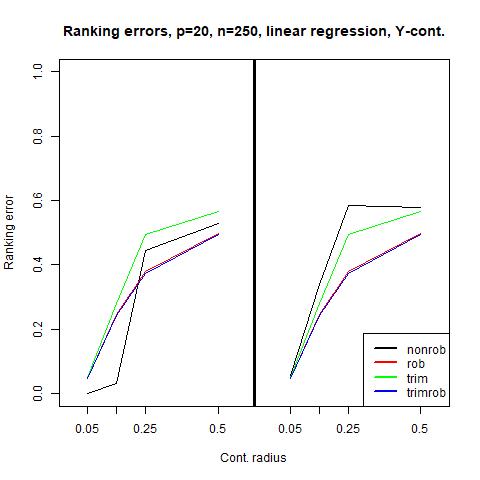} 
\includegraphics[width=5cm,height=4cm]{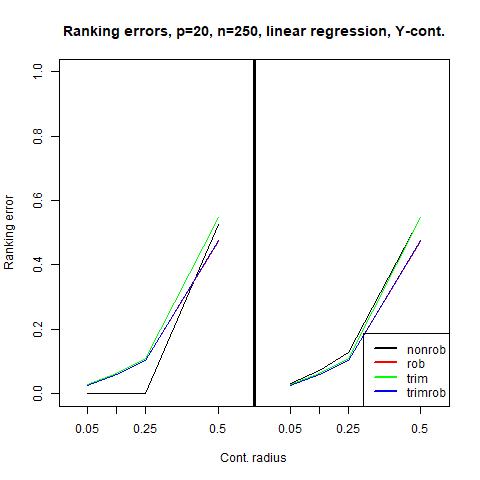} \\
\caption{Hard ranking errors in the context of $Y$-contamination for random cross-validation according to the training losses when trimming the test data according to a non-robust model (upper left) and a robust model (upper right), and according to the test losses when trimming the test data according to a non-robust model (bottom left) and a robust model (bottom right).} \label{traintrimp20Ycont}
\end{center}
\end{figure}

\begin{figure}[H]
\begin{center}
\includegraphics[width=5cm,height=4cm]{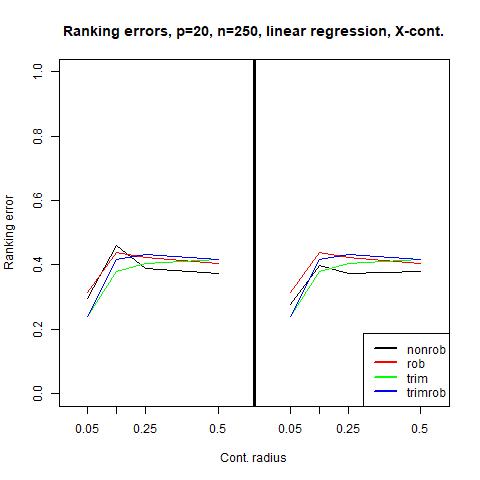}
\includegraphics[width=5cm,height=4cm]{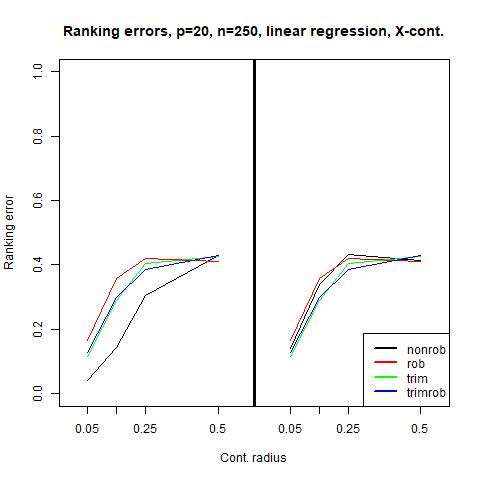} \\
\includegraphics[width=5cm,height=4cm]{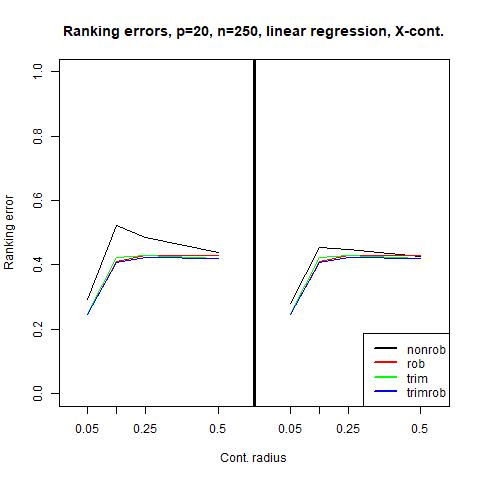} 
\includegraphics[width=5cm,height=4cm]{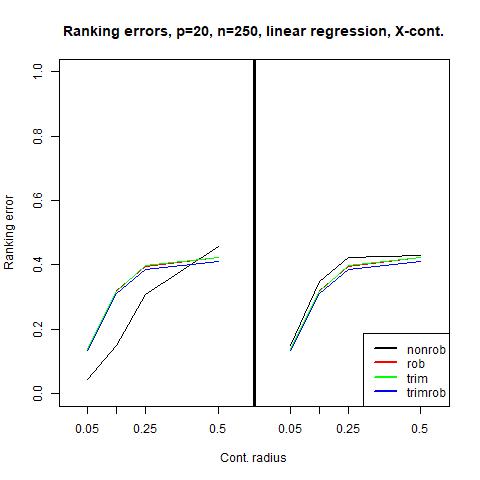} \\
\caption{Hard ranking errors in the context of $X$-contamination for random cross-validation according to the training losses when trimming the test data according to a non-robust model (upper left) and a robust model (upper right), and according to the test losses when trimming the test data according to a non-robust model (bottom left) and a robust model (bottom right).} \label{traintrimp20Xcont}
\end{center}
\end{figure}

\subsubsection{$p=250$, regression, loss-based}

\begin{figure}[H]
\begin{center}
\includegraphics[width=5cm,height=4cm]{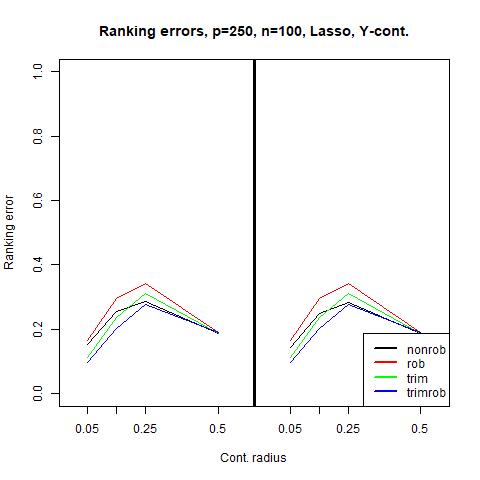}
\includegraphics[width=5cm,height=4cm]{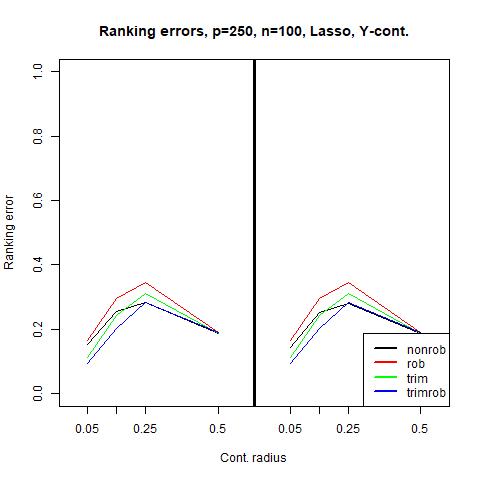} \\
\includegraphics[width=5cm,height=4cm]{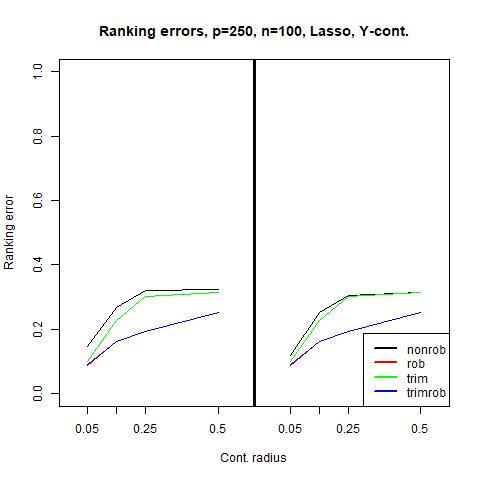} 
\includegraphics[width=5cm,height=4cm]{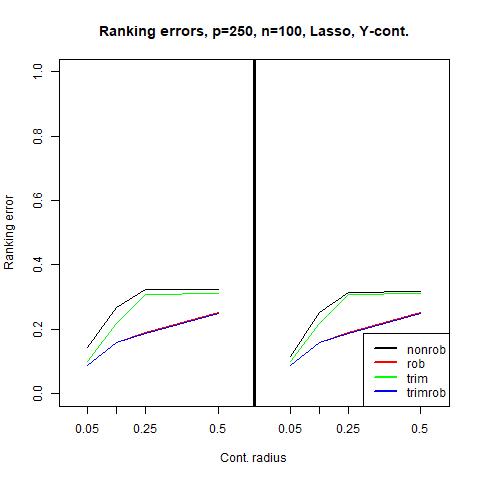} \\
\caption{Hard ranking errors in the context of $Y$-contamination for random cross-validation according to the training losses when trimming the test data according to a non-robust model (upper left) and a robust model (upper right), and according to the test losses when trimming the test data according to a non-robust model (bottom left) and a robust model (bottom right).} \label{traintrimp250Ycont}
\end{center}
\end{figure}

\begin{figure}[H]
\begin{center}
\includegraphics[width=5cm,height=4cm]{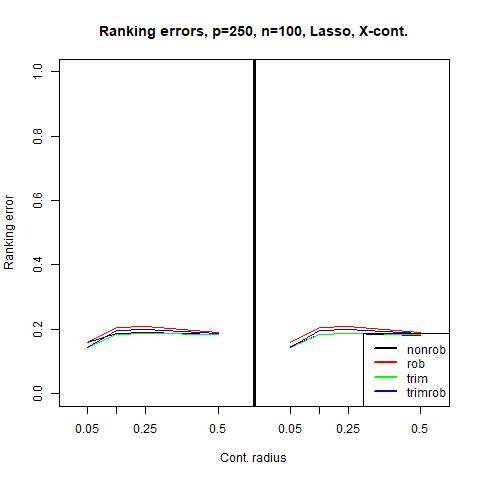}
\includegraphics[width=5cm,height=4cm]{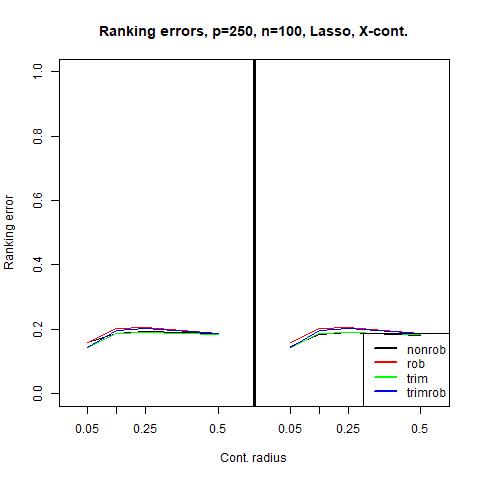} \\
\includegraphics[width=5cm,height=4cm]{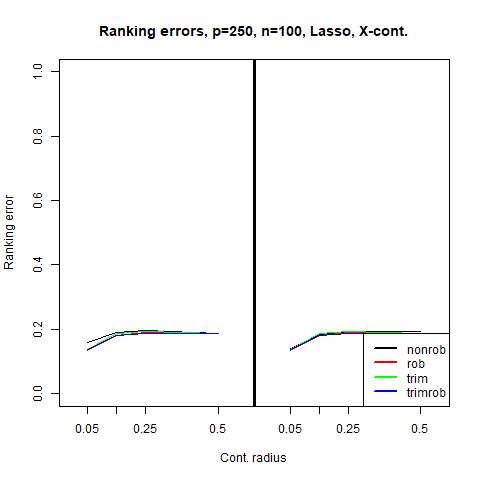} 
\includegraphics[width=5cm,height=4cm]{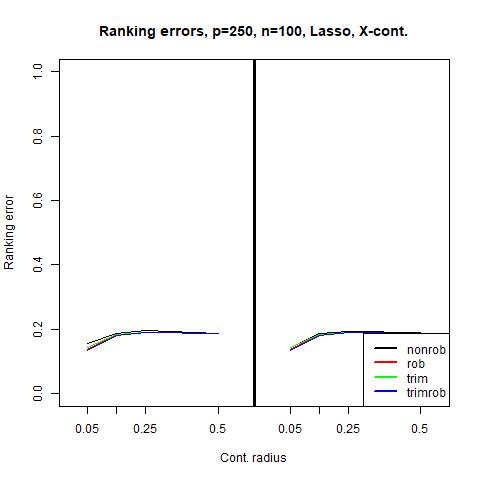} \\
\caption{Hard ranking errors in the context of $X$-contamination for random cross-validation according to the training losses when trimming the test data according to a non-robust model (upper left) and a robust model (upper right), and according to the test losses when trimming the test data according to a non-robust model (bottom left) and a robust model (bottom right).} \label{traintrimp250Xcont}
\end{center}
\end{figure}

\subsubsection{$p=500$, regression, loss-based}

\begin{figure}[H]
\begin{center}
\includegraphics[width=5cm,height=4cm]{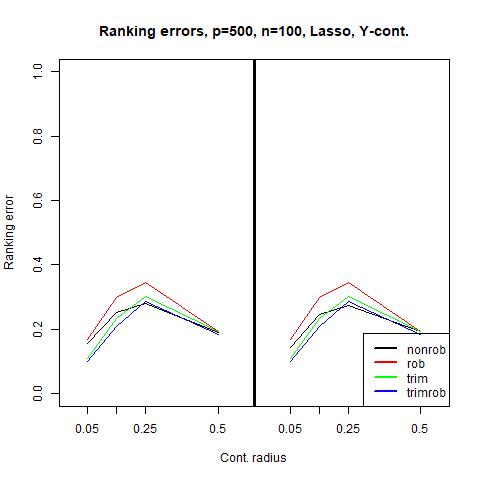}
\includegraphics[width=5cm,height=4cm]{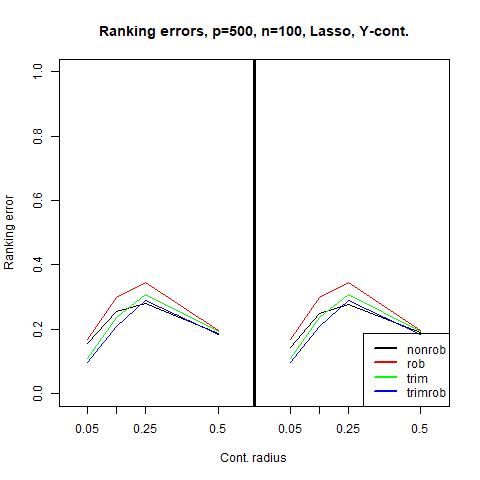} \\
\includegraphics[width=5cm,height=4cm]{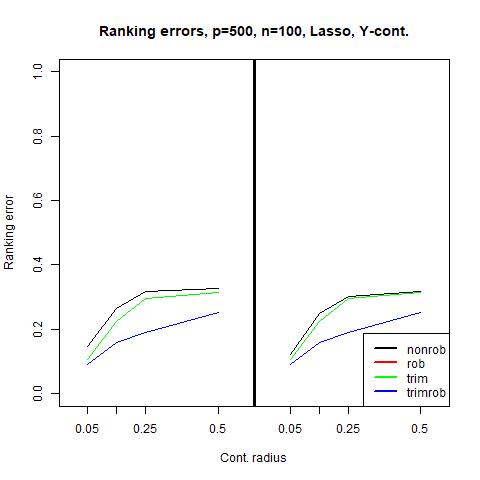} 
\includegraphics[width=5cm,height=4cm]{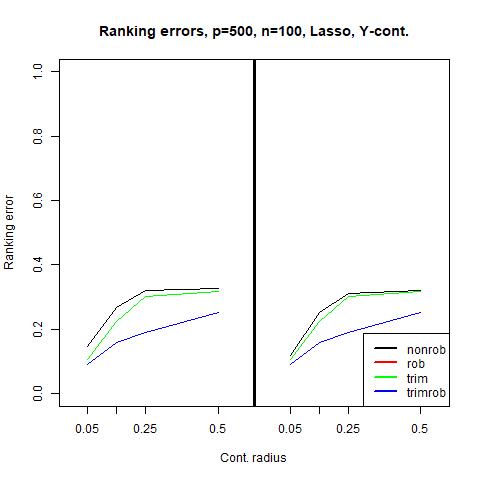} \\
\caption{Hard ranking errors in the context of $Y$-contamination for random cross-validation according to the training losses when trimming the test data according to a non-robust model (upper left) and a robust model (upper right), and according to the test losses when trimming the test data according to a non-robust model (bottom left) and a robust model (bottom right).} \label{traintrimp500Ycont}
\end{center}
\end{figure}

\begin{figure}[H]
\begin{center}
\includegraphics[width=5cm,height=4cm]{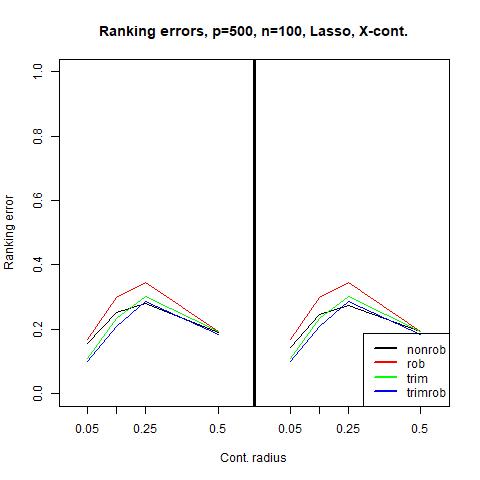}
\includegraphics[width=5cm,height=4cm]{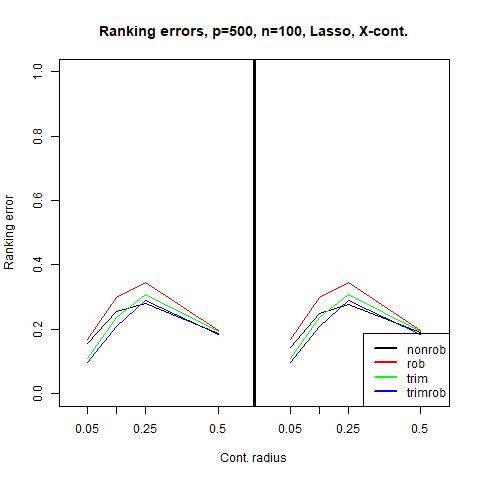} \\
\includegraphics[width=5cm,height=4cm]{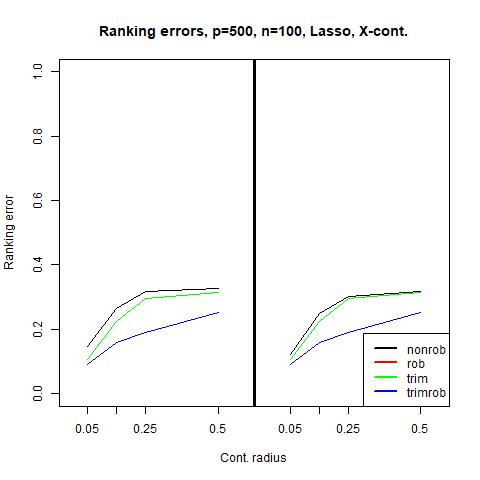} 
\includegraphics[width=5cm,height=4cm]{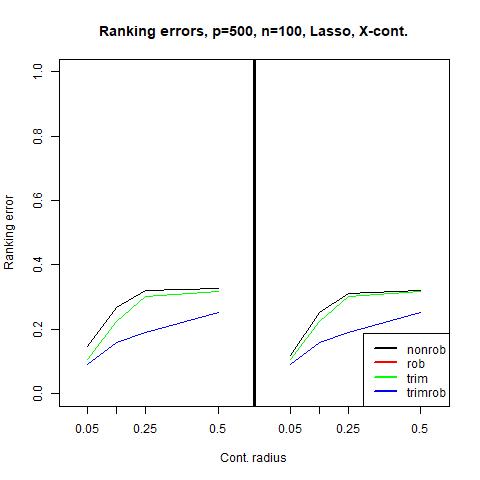} \\
\caption{Hard ranking errors in the context of $X$-contamination for random cross-validation according to the training losses when trimming the test data according to a non-robust model (upper left) and a robust model (upper right), and according to the test losses when trimming the test data according to a non-robust model (bottom left) and a robust model (bottom right).} \label{traintrimp500Xcont}
\end{center}
\end{figure}

\subsubsection{Regression, coefficient-based}

\begin{figure}[H]
\begin{center}
\includegraphics[width=5cm,height=4cm]{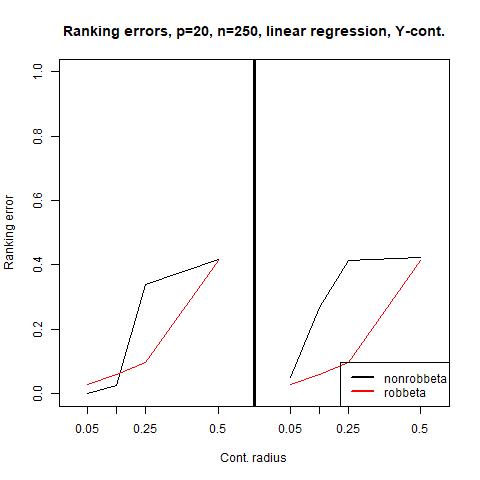}
\includegraphics[width=5cm,height=4cm]{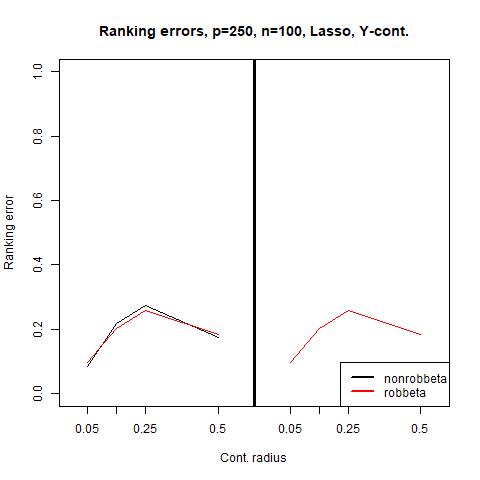}
\includegraphics[width=5cm,height=4cm]{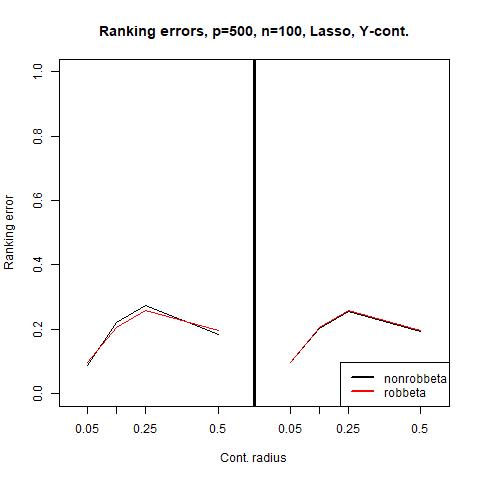}\\
\includegraphics[width=5cm,height=4cm]{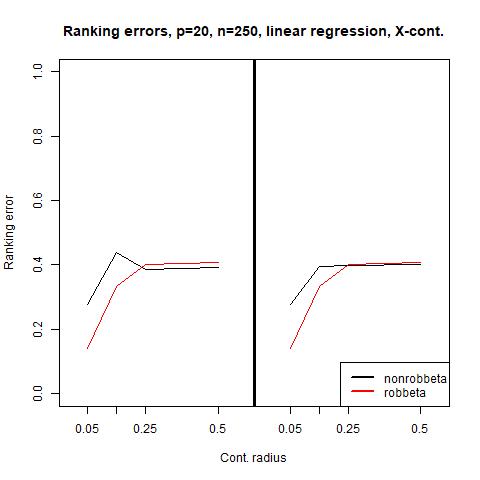}
\includegraphics[width=5cm,height=4cm]{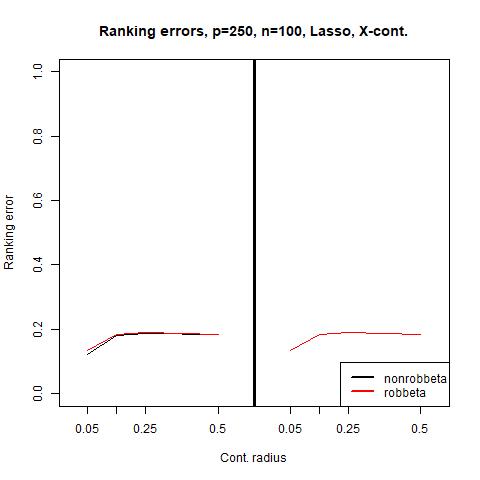}
\includegraphics[width=5cm,height=4cm]{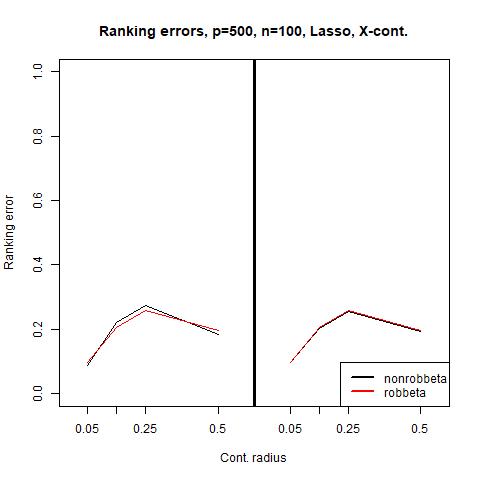}\\
\caption{Hard ranking errors in the context of $Y$-contamination (top row) and $X$-contamination (bottom row), respectively, for randomized cross validation for $p=20$ (left), $p=250$ (middle), and $p=500$ (right) according to the differences in the coefficient vectors of the robust and non-robust model. } \label{traintrimcoeffdiff}
\end{center}
\end{figure}

\subsubsection{$p=20$, classification, loss-based}

\vspace{-0.5cm}

\begin{figure}[H]
\begin{center}
\includegraphics[width=5cm,height=4cm]{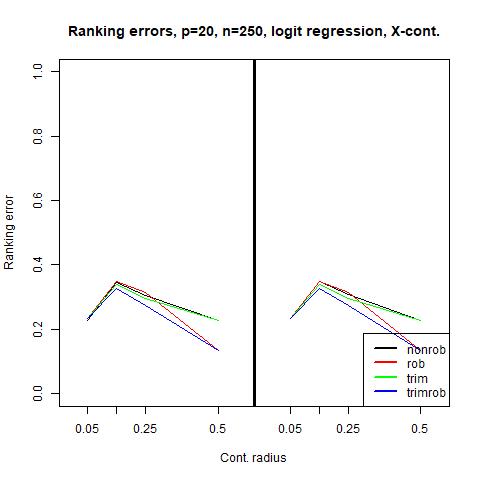}
\includegraphics[width=5cm,height=4cm]{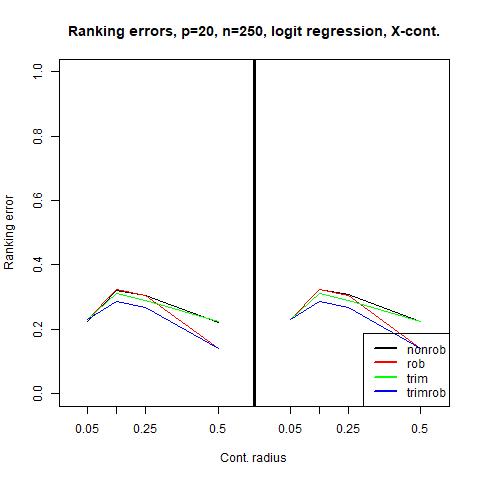} \\
\includegraphics[width=5cm,height=4cm]{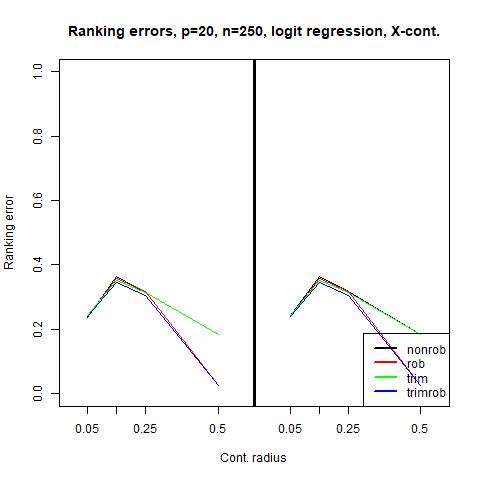} 
\includegraphics[width=5cm,height=4cm]{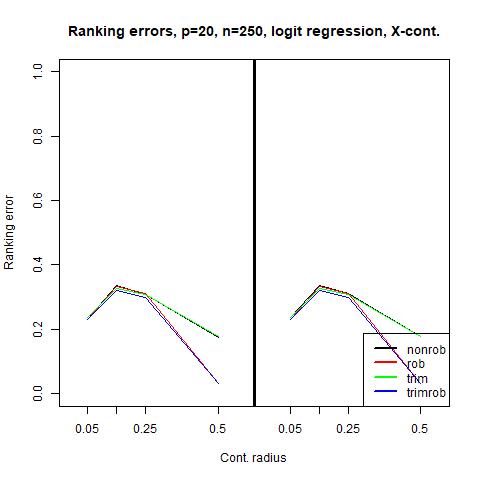} \\
\caption{Hard ranking errors in the context of $X$-contamination for random cross-validation according to the training losses when trimming the test data according to a non-robust model (upper left) and a robust model (upper right), and according to the test losses when trimming the test data according to a non-robust model (bottom left) and a robust model (bottom right).} \label{traintrimp20Xcontglm}
\end{center}
\end{figure}

\subsubsection{$p=250$, classification, loss-based}

\vspace{-0.5cm}

\begin{figure}[H]
\begin{center}
\includegraphics[width=5cm,height=4cm]{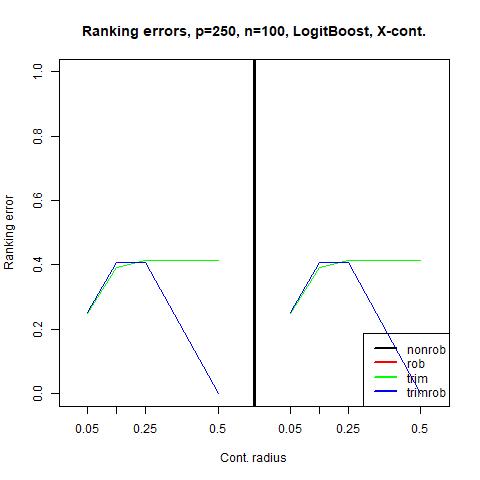}
\includegraphics[width=5cm,height=4cm]{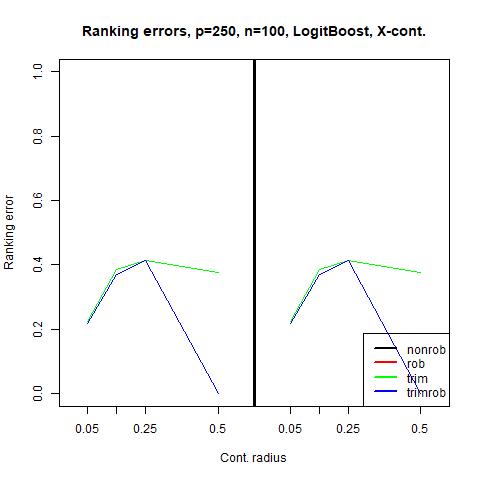} \\
\includegraphics[width=5cm,height=4cm]{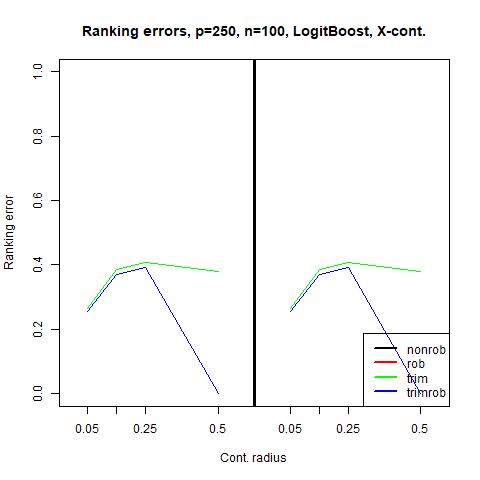} 
\includegraphics[width=5cm,height=4cm]{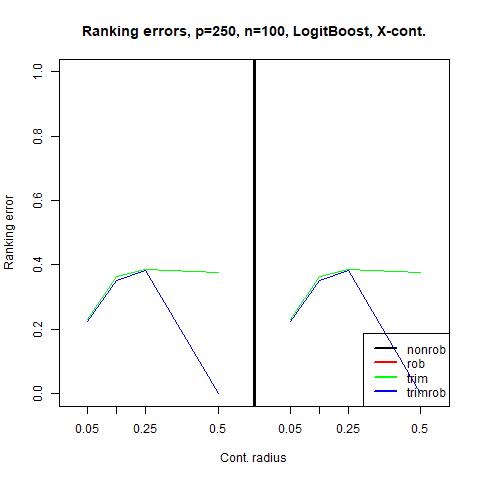} \\
\caption{Hard ranking errors in the context of $X$-contamination for random cross-validation according to the training losses when trimming the test data according to a model trained by LogitBoost (upper left) and AUC-Boosting (upper right), and according to the test losses when trimming the test data according to a model trained by LogitBoost (bottom left) and AUC-Boosting (bottom right), respectively.} \label{traintrimp250Xcontglm}
\end{center}
\end{figure}

\subsubsection{$p=500$, classification, loss-based}

\begin{figure}[H]
\begin{center}
\includegraphics[width=5cm,height=4cm]{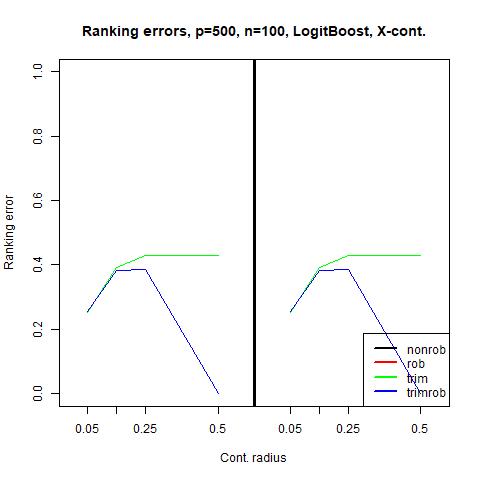}
\includegraphics[width=5cm,height=4cm]{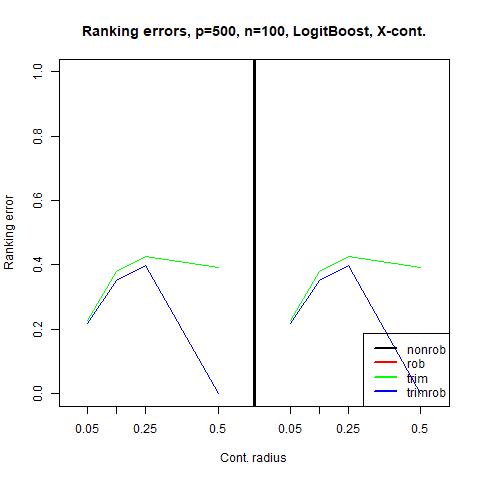} \\
\includegraphics[width=5cm,height=4cm]{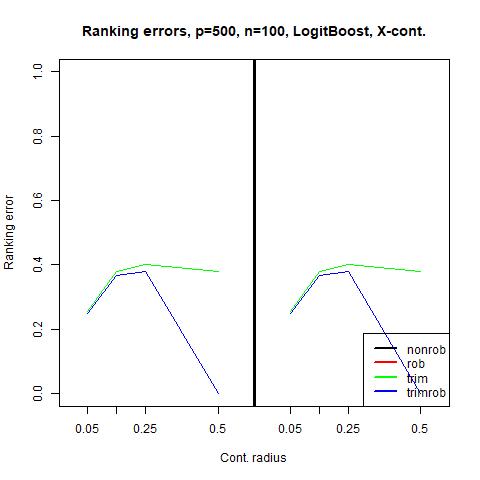} 
\includegraphics[width=5cm,height=4cm]{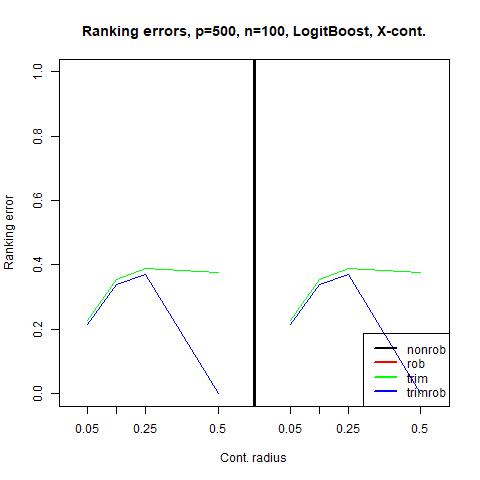} \\
\caption{Hard ranking errors in the context of $X$-contamination for random cross-validation according to the training losses when trimming the test data according to a model trained by LogitBoost (upper left) and AUC-Boosting (upper right), and according to the test losses when trimming the test data according to a model trained by LogitBoost (bottom left) and AUC-Boosting (bottom right), respectively.} \label{traintrimp500Xcontglm}
\end{center}
\end{figure}

\subsubsection{Classification, coefficient-based}

\begin{figure}[H]
\begin{center}
\includegraphics[width=5cm,height=4cm]{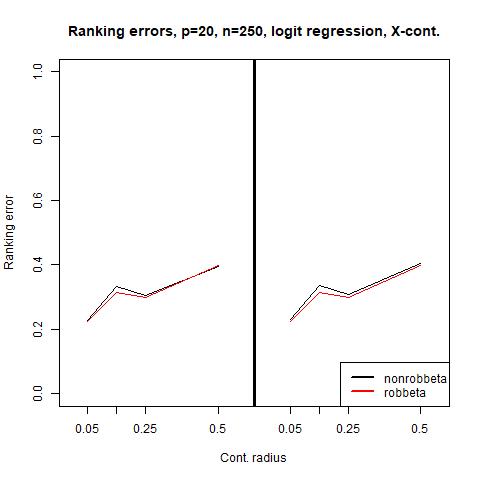}
\includegraphics[width=5cm,height=4cm]{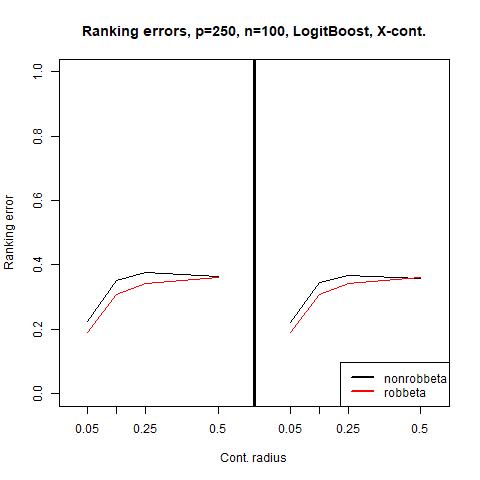}
\includegraphics[width=5cm,height=4cm]{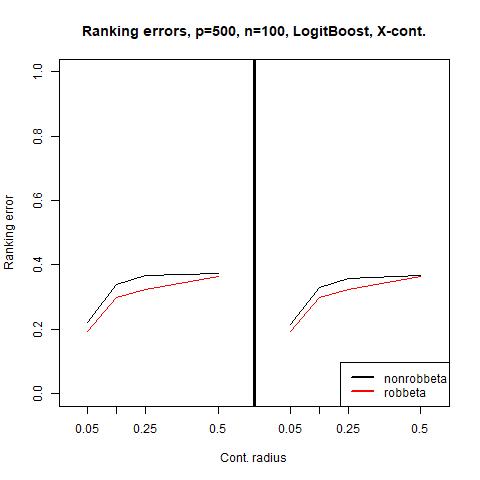}\\
\caption{Hard ranking errors in the context of $X$-contamination for randomized cross validation for $p=20$ (left), $p=250$ (middle), and $p=500$ (right) according to the differences in the coefficient vectors. } \label{traintrimglmcoeffdiff}
\end{center}
\end{figure}

\subsection{Contaminated test data, post-trimming}  \label{app:valtrim}

For the ranking errors according to the training losses, which we also evaluate here in order to enable a direct comparison, the results from the first set of experiments in Sec. \ref{app:valclean} are reproduced, as expected, since the test data are not used here. Note that the graphical presentation is different here than in Sec. \ref{app:valclean} due to the additional parameter $r_{val}$. Here, the left three columns in the figures correspond to $SNR=5$ or $\mu=3$, respectively, the right three columns to $SNR=0.5$ and $\mu=0.5$, respectively. The leftmost of the three respective columns corresponds to $r=0.1$, the middle one to $r=2.5$, the right one to $r=0.5$, respectively. The values on the $x$-axis in the individual columns represent $r_{val}$. For ranking errors ased on training losses, the results would therefore be constant w.r.t. $r_{val}$. The slight changes for varying $r_{val}$ only occur because we made separate simulations for each configuration $(r,r_{val})$.

\subsubsection{$p=20$, regression, loss-based}

\begin{figure}[H]
\begin{center}
\includegraphics[width=5cm,height=4cm]{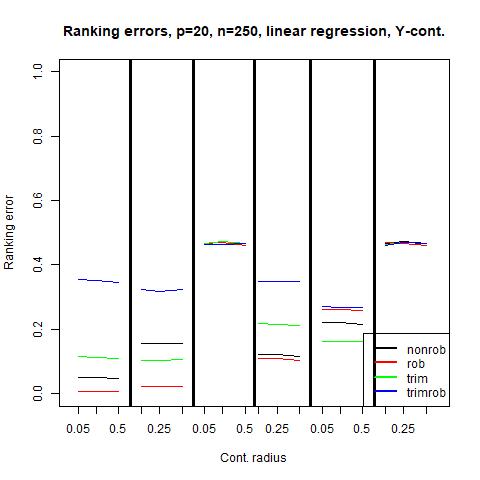}
\includegraphics[width=5cm,height=4cm]{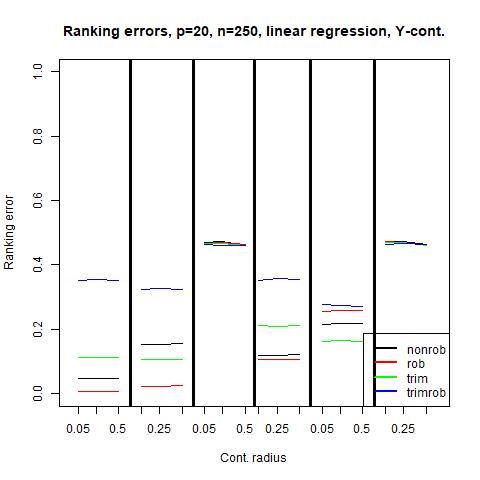} \\
\includegraphics[width=5cm,height=4cm]{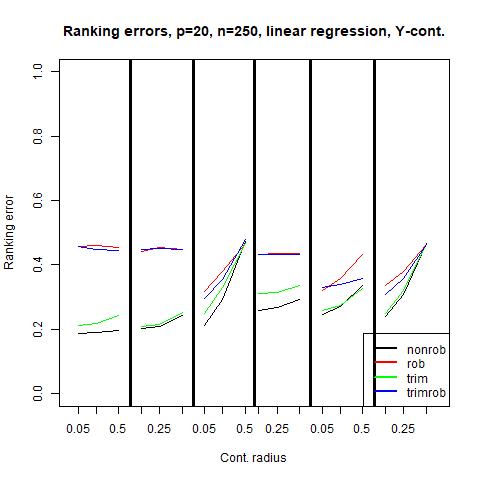} 
\includegraphics[width=5cm,height=4cm]{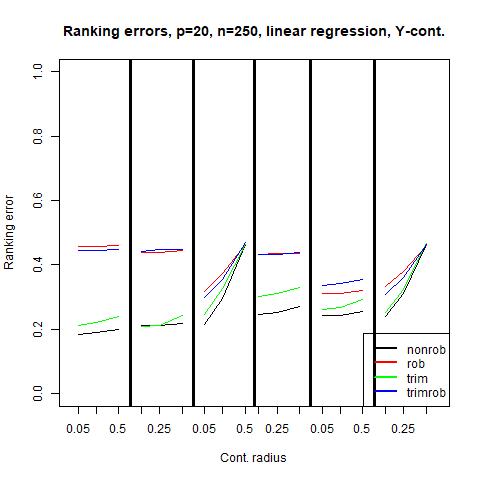} \\
\caption{Hard ranking errors in the context of $Y$-contamination for random cross-validation according to the training losses when trimming the test data according to a non-robust model (upper left) and a robust model (upper right), and according to the test losses when trimming the test data according to a non-robust model (bottom left) and a robust model (bottom right).} \label{valtrimp20Ycont}
\end{center}
\end{figure}

\begin{figure}[H]
\begin{center}
\includegraphics[width=5cm,height=4cm]{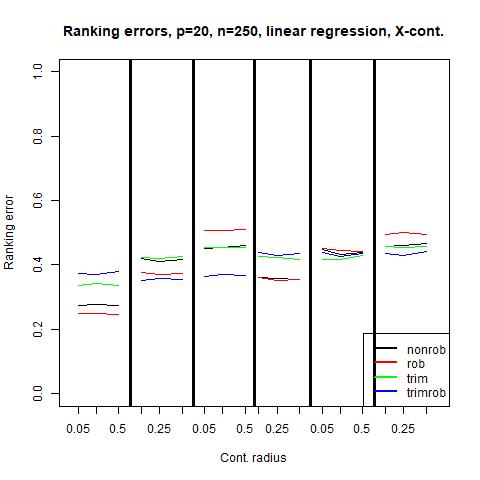}
\includegraphics[width=5cm,height=4cm]{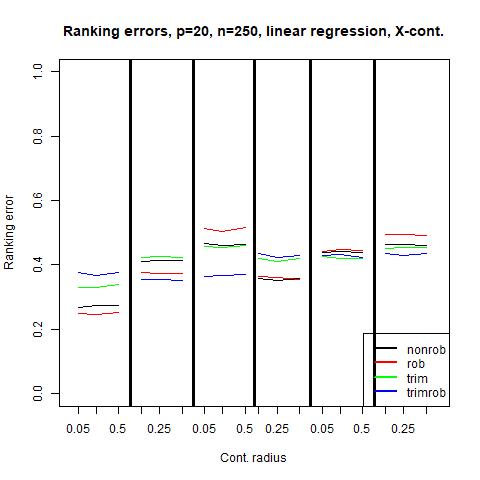} \\
\includegraphics[width=5cm,height=4cm]{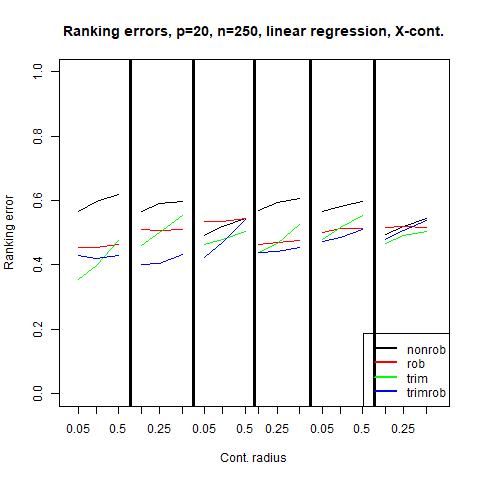} 
\includegraphics[width=5cm,height=4cm]{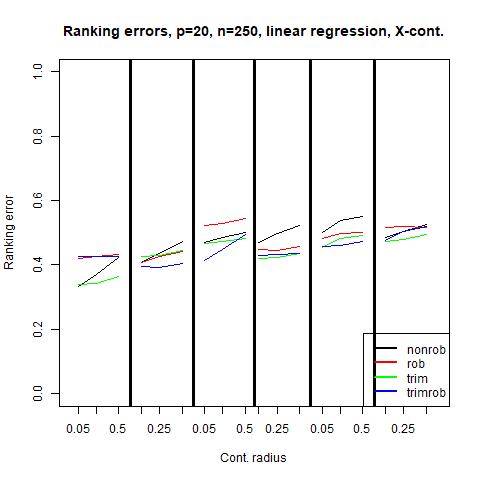} \\
\caption{Hard ranking errors in the context of $X$-contamination for random cross-validation according to the training losses when trimming the test data according to a non-robust model (upper left) and a robust model (upper right), and according to the test losses when trimming the test data according to a non-robust model (bottom left) and a robust model (bottom right).} \label{valtrimp20Xcont}
\end{center}
\end{figure}

\subsubsection{$p=250$, regression, loss-based}

\begin{figure}[H]
\begin{center}
\includegraphics[width=5cm,height=4cm]{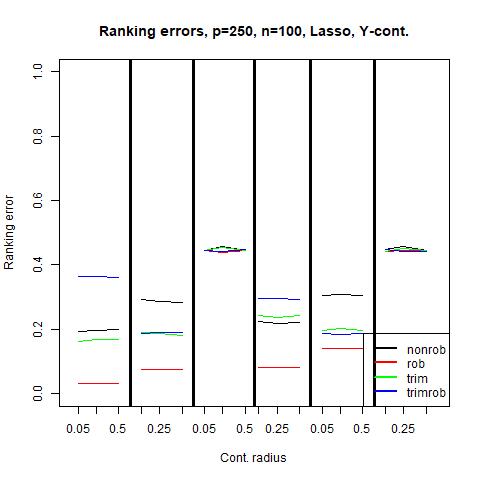}
\includegraphics[width=5cm,height=4cm]{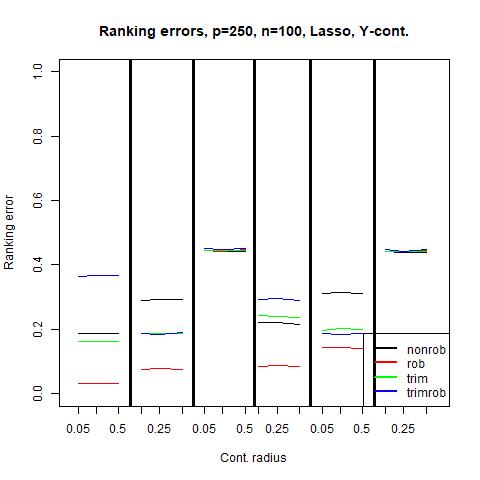} \\
\includegraphics[width=5cm,height=4cm]{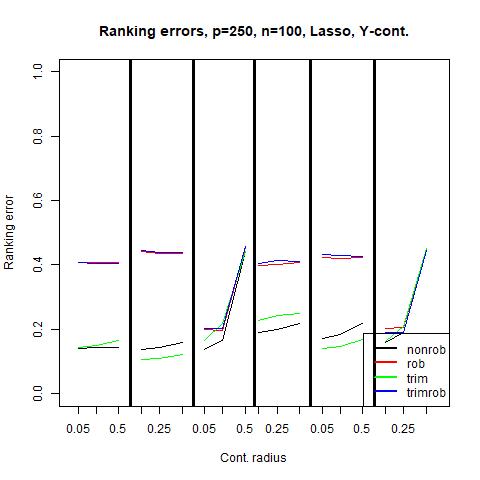} 
\includegraphics[width=5cm,height=4cm]{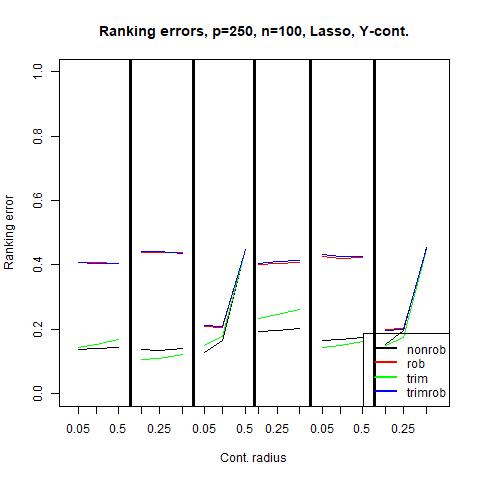} \\
\caption{Hard ranking errors in the context of $Y$-contamination for random cross-validation according to the training losses when trimming the test data according to a non-robust model (upper left) and a robust model (upper right), and according to the test losses when trimming the test data according to a non-robust model (bottom left) and a robust model (bottom right).} \label{valtrimp250Ycont}
\end{center}
\end{figure}

\begin{figure}[H]
\begin{center}
\includegraphics[width=5cm,height=4cm]{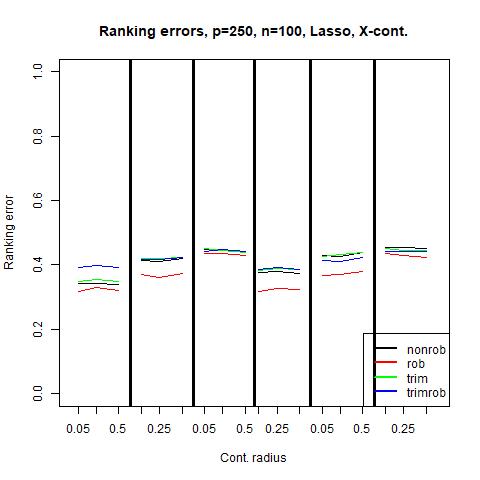}
\includegraphics[width=5cm,height=4cm]{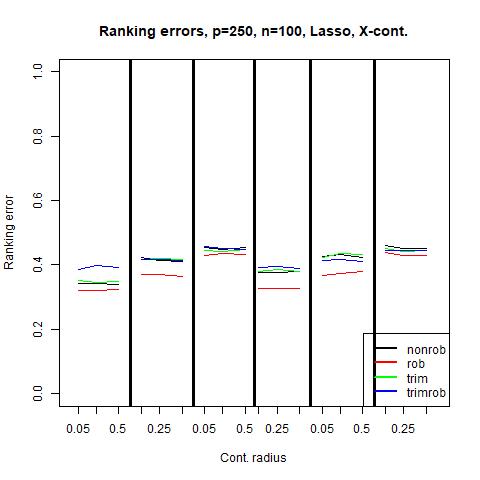} \\
\includegraphics[width=5cm,height=4cm]{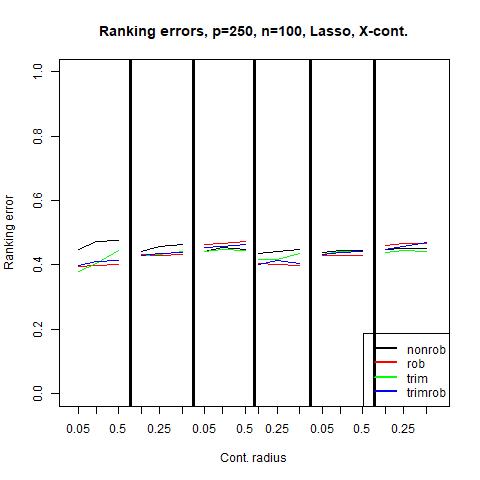} 
\includegraphics[width=5cm,height=4cm]{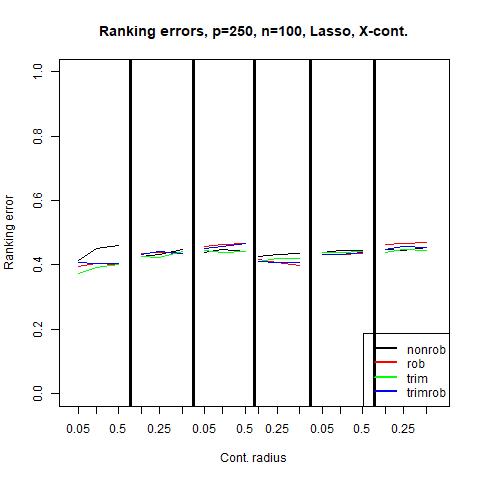} \\
\caption{Hard ranking errors in the context of $X$-contamination for random cross-validation according to the training losses when trimming the test data according to a non-robust model (upper left) and a robust model (upper right), and according to the test losses when trimming the test data according to a non-robust model (bottom left) and a robust model (bottom right).} \label{valtrimp250Xcont}
\end{center}
\end{figure}

\subsubsection{$p=500$, regression, loss-based}

\begin{figure}[H]
\begin{center}
\includegraphics[width=5cm,height=4cm]{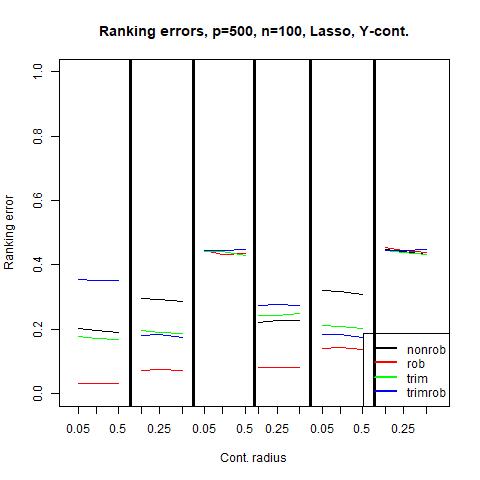}
\includegraphics[width=5cm,height=4cm]{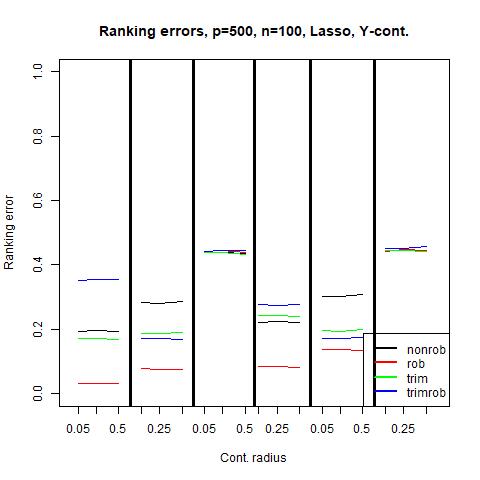} \\
\includegraphics[width=5cm,height=4cm]{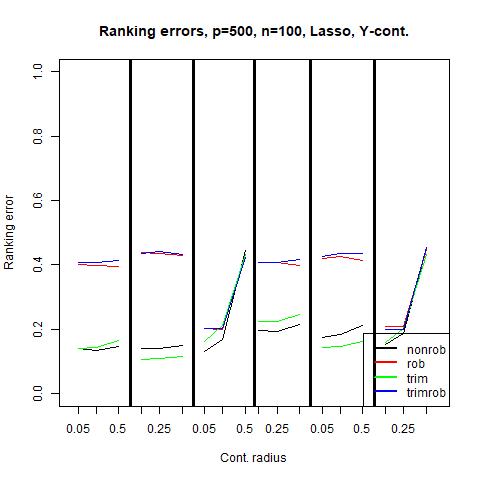} 
\includegraphics[width=5cm,height=4cm]{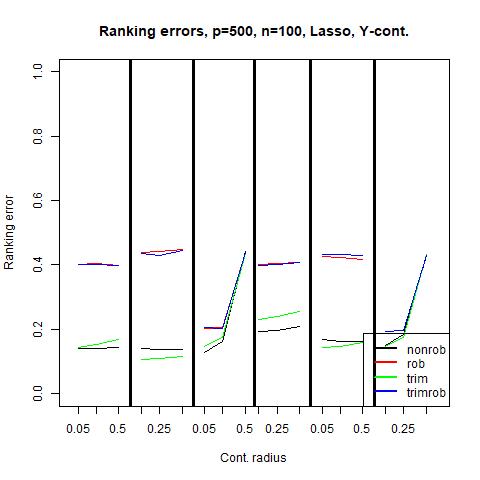} \\
\caption{Hard ranking errors in the context of $Y$-contamination for random cross-validation according to the training losses when trimming the test data according to a non-robust model (upper left) and a robust model (upper right), and according to the test losses when trimming the test data according to a non-robust model (bottom left) and a robust model (bottom right).} \label{valtrimp500Ycont}
\end{center}
\end{figure}

\begin{figure}[H]
\begin{center}
\includegraphics[width=5cm,height=4cm]{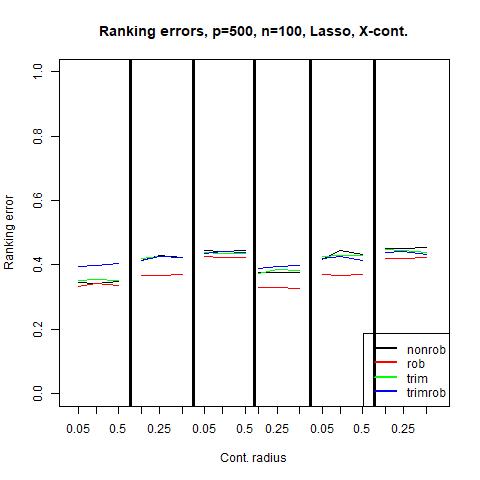}
\includegraphics[width=5cm,height=4cm]{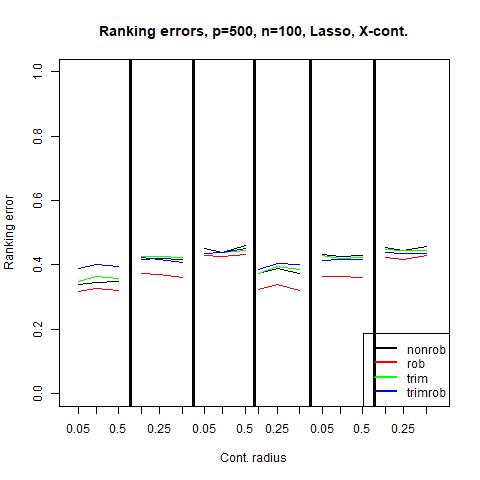} \\
\includegraphics[width=5cm,height=4cm]{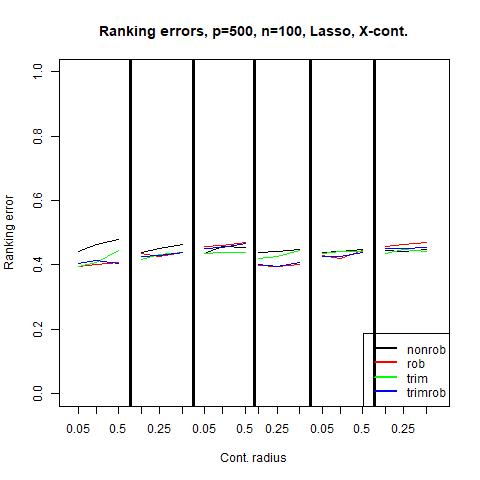} 
\includegraphics[width=5cm,height=4cm]{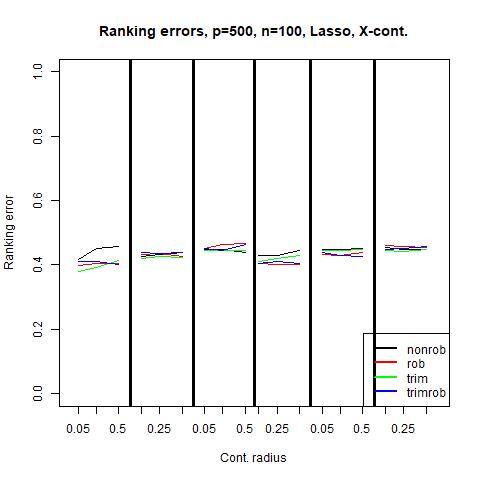} \\
\caption{Hard ranking errors in the context of $X$-contamination for random cross-validation according to the training losses when trimming the test data according to a non-robust model (upper left) and a robust model (upper right), and according to the test losses when trimming the test data according to a non-robust model (bottom left) and a robust model (bottom right).} \label{valtrimp500Xcont}
\end{center}
\end{figure}

\subsubsection{Regression, coefficient-based}

\begin{figure}[H]
\begin{center}
\includegraphics[width=5cm,height=4cm]{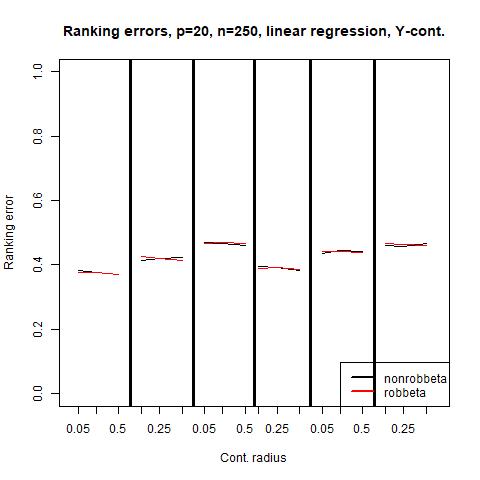}
\includegraphics[width=5cm,height=4cm]{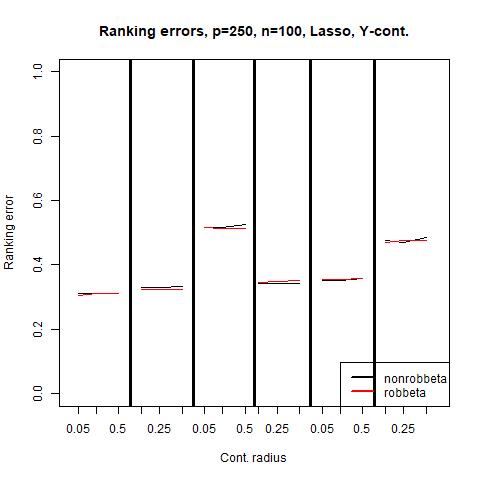}
\includegraphics[width=5cm,height=4cm]{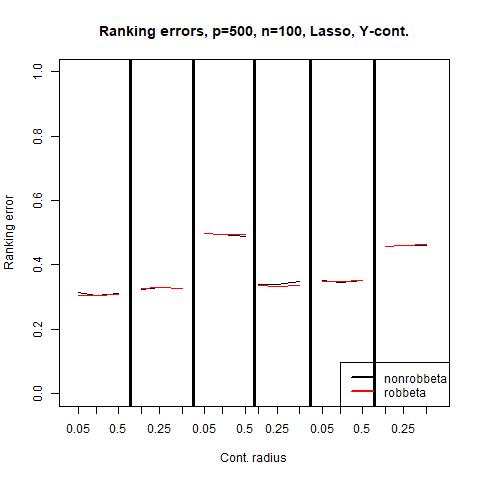}\\
\includegraphics[width=5cm,height=4cm]{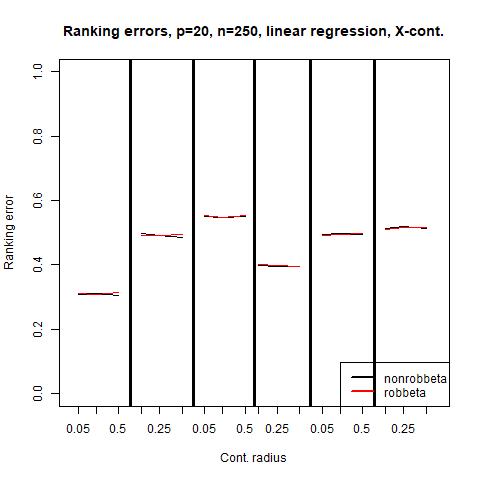}
\includegraphics[width=5cm,height=4cm]{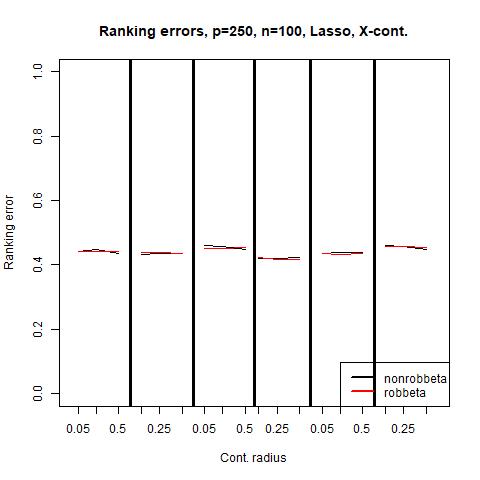}
\includegraphics[width=5cm,height=4cm]{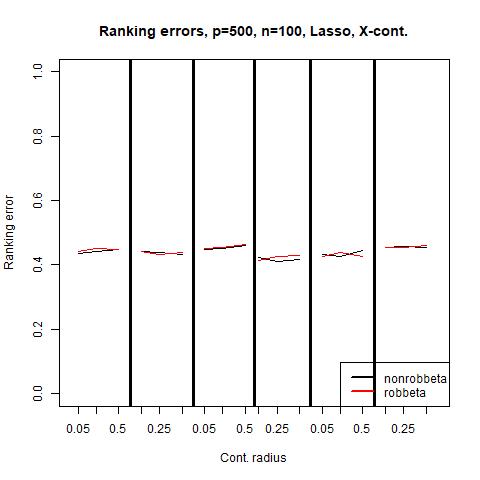}\\
\caption{Hard ranking errors in the context of $Y$-contamination (top row) and $X$-contamination (bottom row), respectively, for randomized cross validation for $p=20$ (left), $p=250$ (middle), and $p=500$ (right) according to the differences in the coefficient vectors. } \label{valtrimcoeffdiff}
\end{center}
\end{figure}

\subsubsection{$p=20$, classification, loss-based}

\vspace{-0.5cm}
\begin{figure}[H]
\begin{center}
\includegraphics[width=5cm,height=4cm]{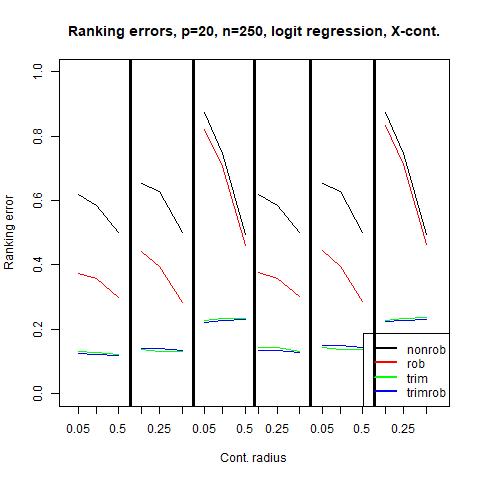}
\includegraphics[width=5cm,height=4cm]{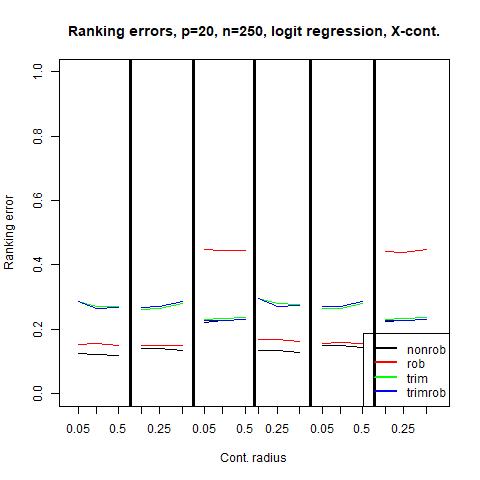} \\
\includegraphics[width=5cm,height=4cm]{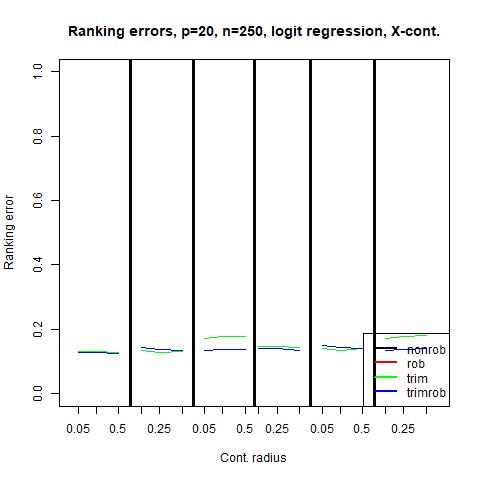} 
\includegraphics[width=5cm,height=4cm]{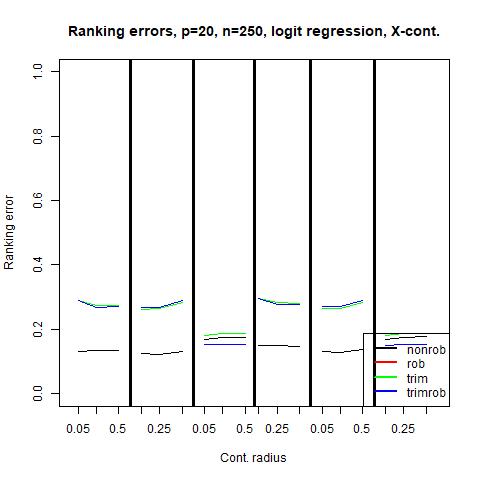} \\
\caption{Hard ranking errors in the context of $X$-contamination for random cross-validation according to the training losses when trimming the test data according to a non-robust (upper left) and a robust binary classification model (upper right), and according to the test losses when trimming the test data according to a non-robust (bottom left) and a robust classification model (bottom right).} \label{valtrimp20Xcontglm}
\end{center}
\end{figure}

\subsubsection{$p=250$, classification, loss-based}

\vspace{-0.5cm}
\begin{figure}[H]
\begin{center}
\includegraphics[width=5cm,height=4cm]{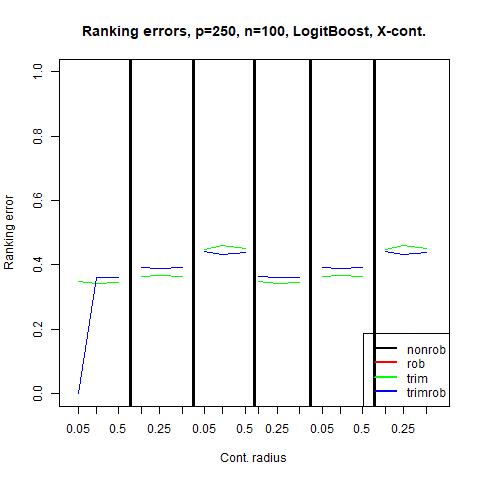}
\includegraphics[width=5cm,height=4cm]{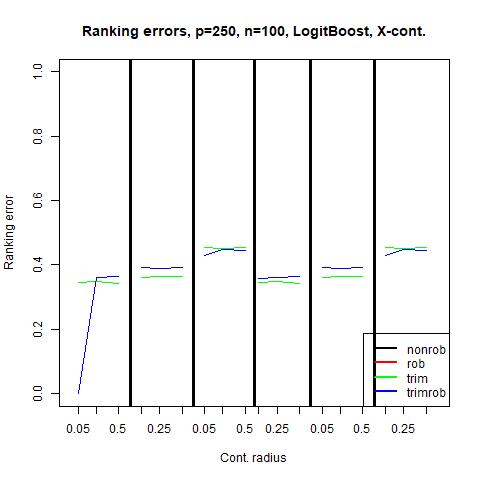} \\
\includegraphics[width=5cm,height=4cm]{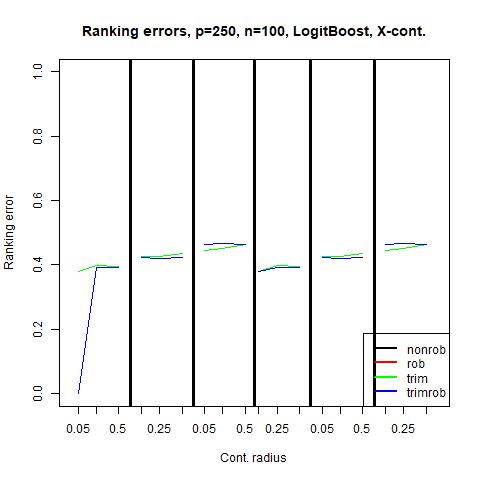} 
\includegraphics[width=5cm,height=4cm]{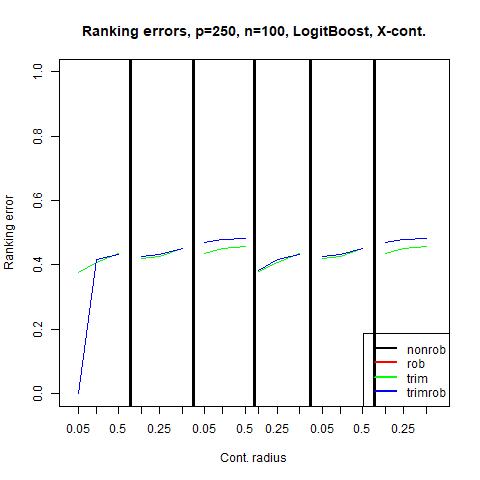} \\
\caption{Hard ranking errors in the context of $X$-contamination for random cross-validation according to the training losses when trimming the test data according to a LogitBoost model (upper left) and an AUC-Boosting model (upper right), and according to the test losses when trimming the test data according to a LogitBoost model (bottom left) and an AUC-Boosting model (bottom right).} \label{valtrimp250Xcontglm}
\end{center}
\end{figure}

\subsubsection{$p=500$, classification, loss-based}

\begin{figure}[H]
\begin{center}
\includegraphics[width=5cm,height=4cm]{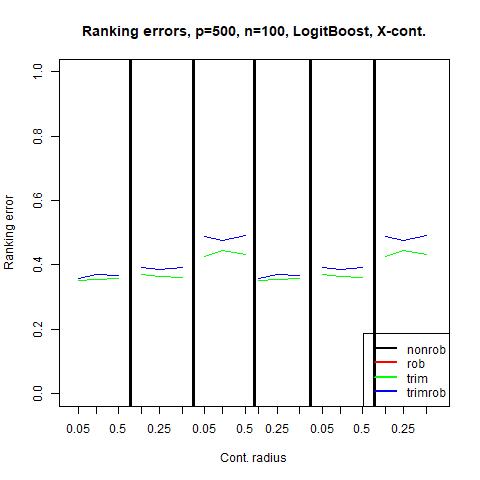}
\includegraphics[width=5cm,height=4cm]{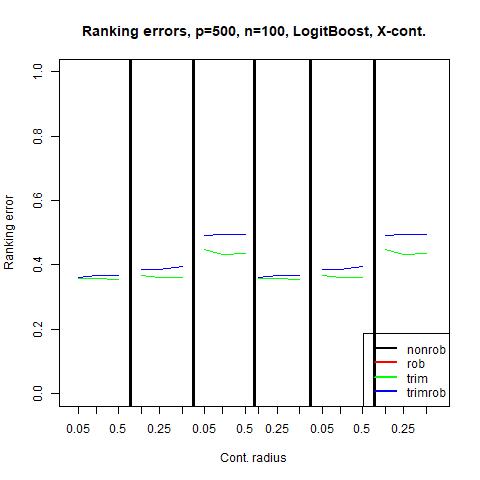} \\
\includegraphics[width=5cm,height=4cm]{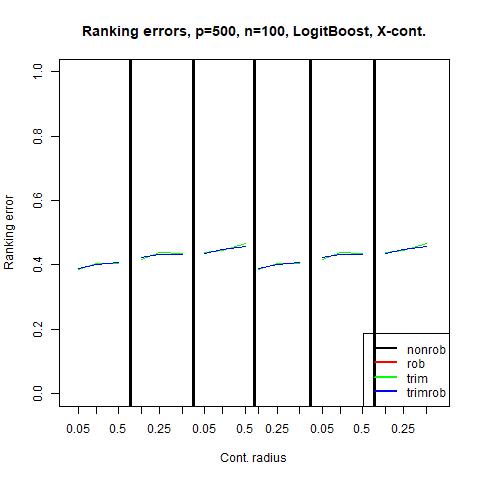} 
\includegraphics[width=5cm,height=4cm]{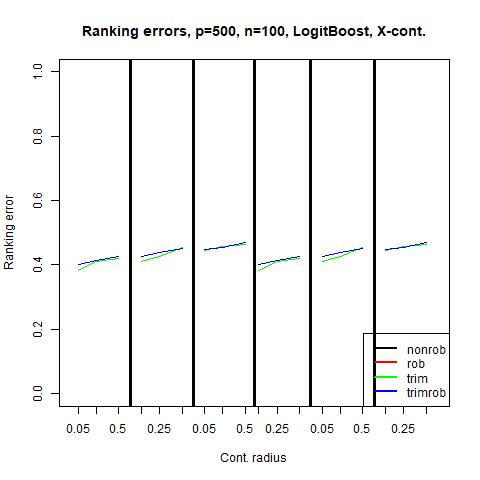} \\
\caption{Hard ranking errors in the context of $X$-contamination for random cross-validation according to the training losses when trimming the test data according to a LogitBoost model (upper left) and an AUC-Boosting model (upper right), and according to the test losses when trimming the test data according to a LogitBoost model (bottom left) and an AUC-Boosting model (bottom right).} \label{valtrimp500Xcontglm}
\end{center}
\end{figure}

\subsubsection{Classification, coefficient-based}

\begin{figure}[H]
\begin{center}
\includegraphics[width=5cm,height=4cm]{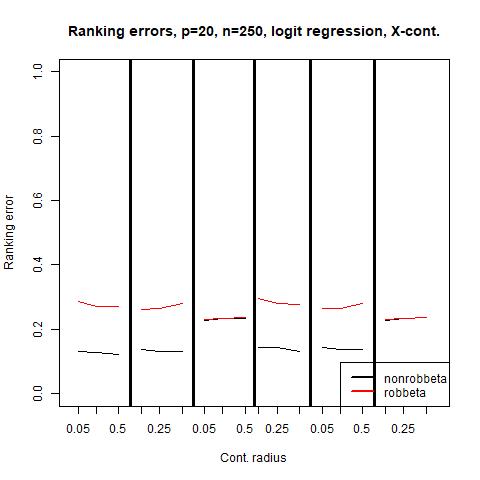}
\includegraphics[width=5cm,height=4cm]{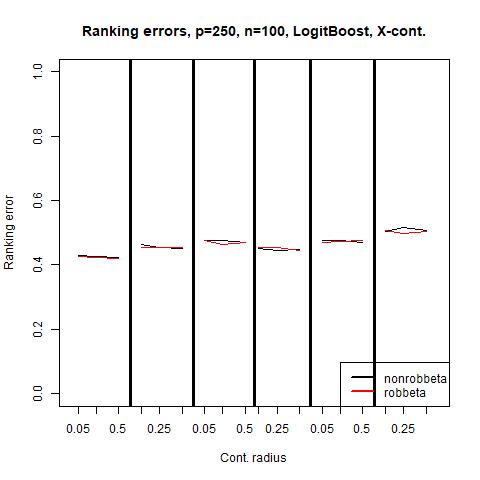}
\includegraphics[width=5cm,height=4cm]{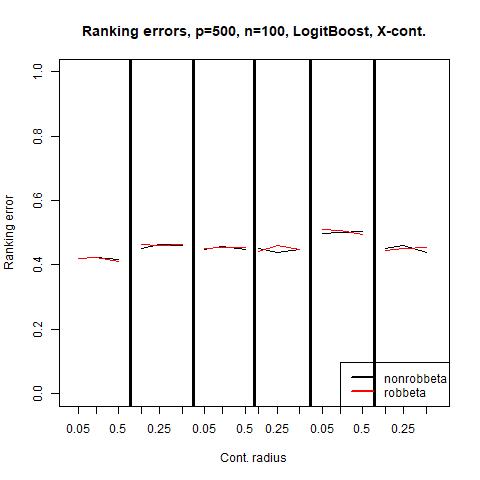}\\
\caption{Hard ranking errors in the context of $X$-contamination for randomized cross validation for $p=20$ (left), $p=250$ (middle), and $p=500$ (right) according to the differences in the coefficient vectors. } \label{valtrimglmcoeffdiff}
\end{center}
\end{figure}

\section{Trimming single instances} \label{app:triminstance}

\subsection{Contaminated training and test data, pre-trimming} \label{app:valcont}

In the following, we depict the boxplots of the weak ranking errors for the identified validation and training instances, respectively. Moreover, we evaluate the mean empirical BDP of elicitability, which is the mean number of repetitions where the average weak ranking error for the identification of test instances is non-zero.

Here, due to already having 9 configurations of $(r,r_{val})$, we again only use an SNR of 5 and an SNR of 0.5 for regression, and $\mu=3$ and $\mu=0.5$ for classification. In the figures, the three left parts correspond to an SNR of 5 or $\mu=3$, respectively, the three right parts to an SNR of 0.5 and $\mu=0.5$, respectively. The contamination radius of the training data is $0.1$, $0.25$, and $0.5$, from left to right, respectively, and the $x$-axis denotes $r_{val}$. For example, the right boxplot in the fifth part of a figure corresponds to $r=0.25$, $r_{val}=0.5$, and $SNR=0.5$ or $\mu=0.5$, respectively.

\subsubsection{$p=20$, regression}

\begin{figure}[H]
\begin{center}
\includegraphics[width=5cm,height=4cm]{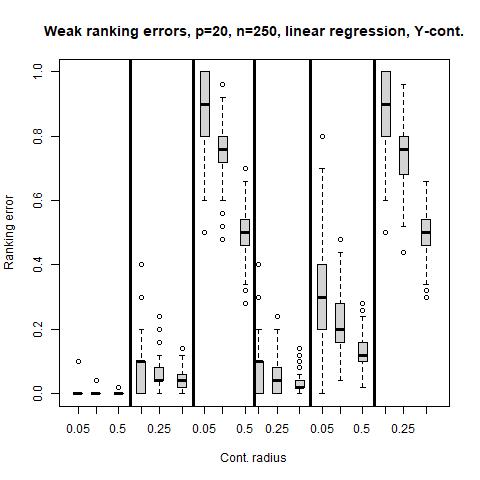}
\includegraphics[width=5cm,height=4cm]{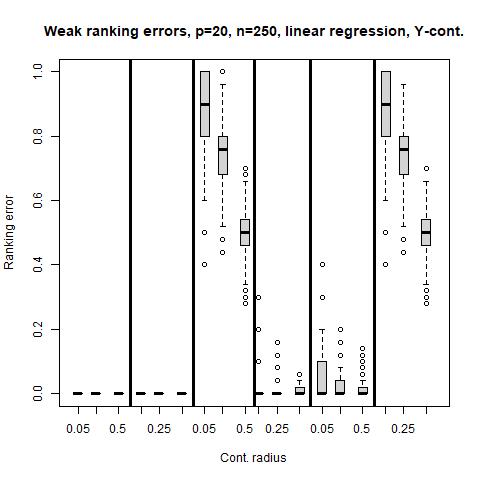} \\
\includegraphics[width=5cm,height=4cm]{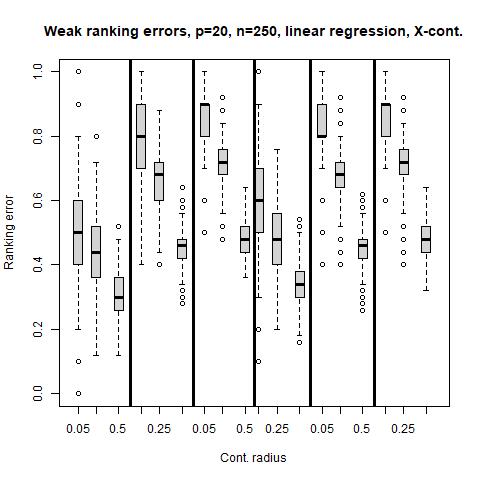} 
\includegraphics[width=5cm,height=4cm]{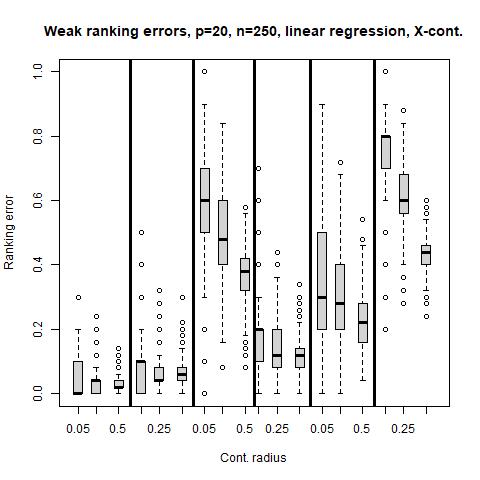} \\
\caption{Boxplots of the weak ranking errors corresponding to the model-based identification of the most outlying test instances. Upper row: $Y$-contamination, bottom row: $X$-contamination, left column: non-robust regression, right column: robust regression.} \label{valcontp20val}
\end{center}
\end{figure}

\begin{figure}[H]
\begin{center}
\includegraphics[width=5cm,height=4cm]{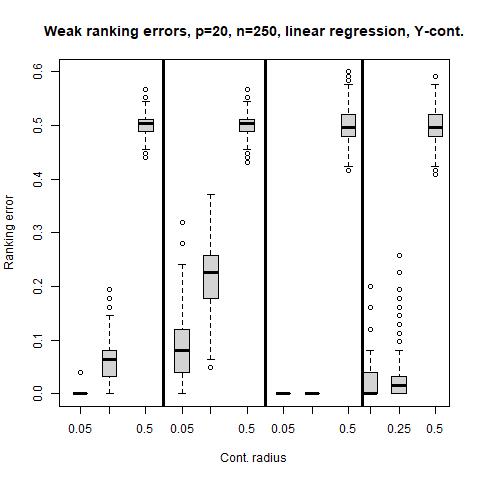}
\includegraphics[width=5cm,height=4cm]{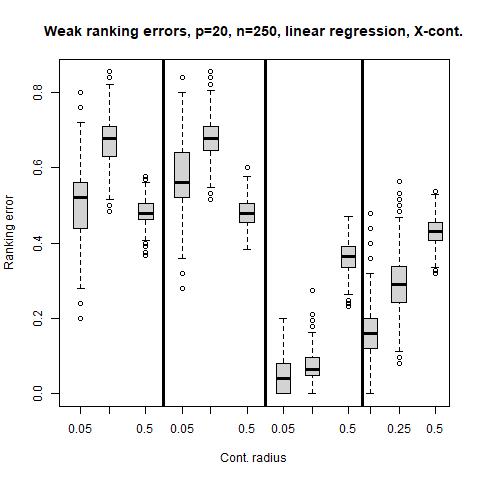} 
\caption{Boxplots of the weak ranking errors corresponding to the leave-one-out model-based identification of the most outlying training instances.  Left figure: $Y$-contamination, right figure: $X$-contamination. Left two columns: Non-robust regression, right two columns: robust regression.} \label{valcontp20tr}
\end{center}
\end{figure}

\begin{figure}[H]
\begin{center}
\includegraphics[width=5cm,height=4cm]{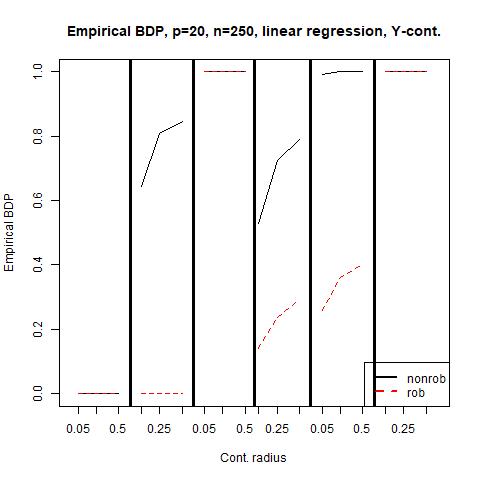}
\includegraphics[width=5cm,height=4cm]{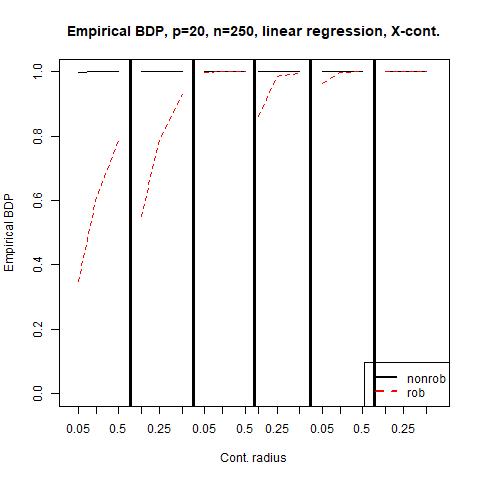} 
\caption{Mean empirical BDP of elicitability, i.e., mean relative number of repetitions where the weak ranking error corresponding to the identification of the contaminated test instances is non-zero.} \label{valcontp20bdp}
\end{center}
\end{figure}

\subsubsection{$p=250$, regression}

\begin{figure}[H]
\begin{center}
\includegraphics[width=5cm,height=4cm]{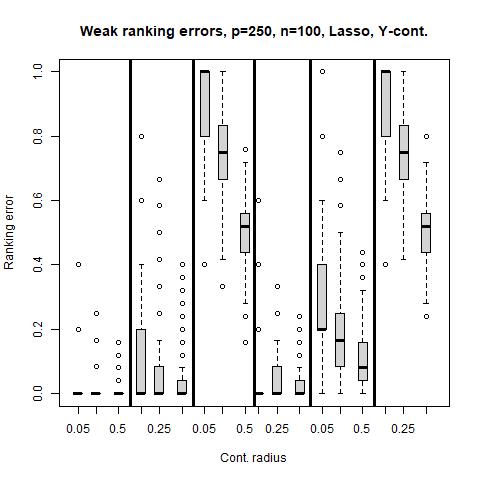}
\includegraphics[width=5cm,height=4cm]{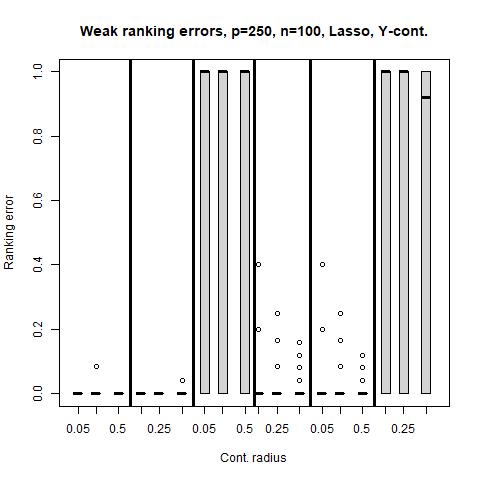} \\
\includegraphics[width=5cm,height=4cm]{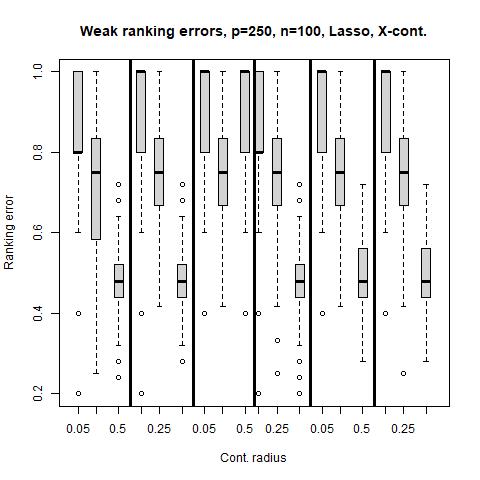} 
\includegraphics[width=5cm,height=4cm]{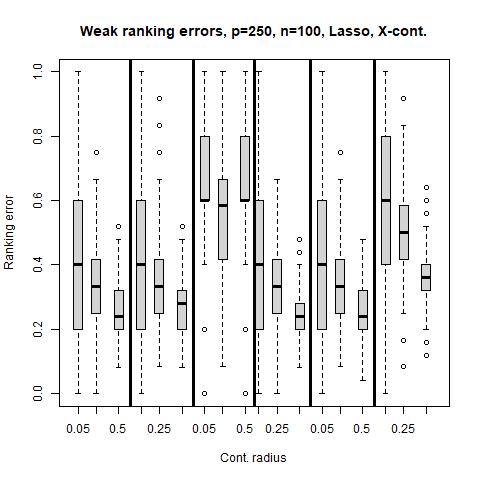} \\
\caption{Boxplots of the weak ranking errors corresponding to the model-based identification of the most outlying test instances. Upper row: $Y$-contamination, bottom row: $X$-contamination, left column: non-robust regression, right column: robust regression.} \label{valcontp250val}
\end{center}
\end{figure}

\begin{figure}[H]
\begin{center}
\includegraphics[width=5cm,height=4cm]{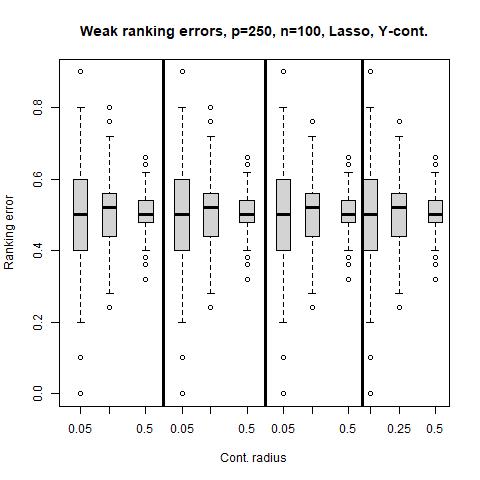}
\includegraphics[width=5cm,height=4cm]{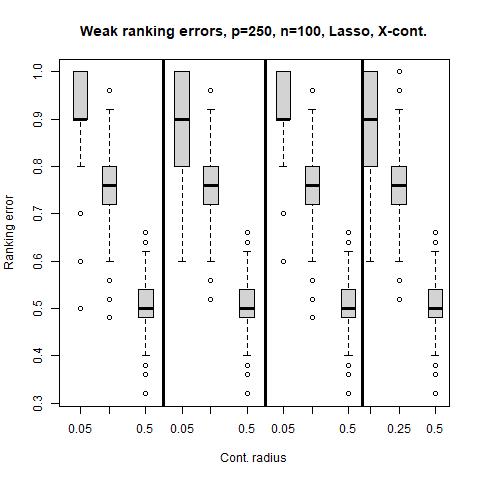} 
\caption{Boxplots of the weak ranking errors corresponding to the leave-one-out model-based identification of the most outlying training instances.  Left figure: $Y$-contamination, right figure: $X$-contamination. Left two columns: Non-robust regression, right two columns: robust regression.} \label{valcontp250tr}
\end{center}
\end{figure}

\begin{figure}[H]
\begin{center}
\includegraphics[width=5cm,height=4cm]{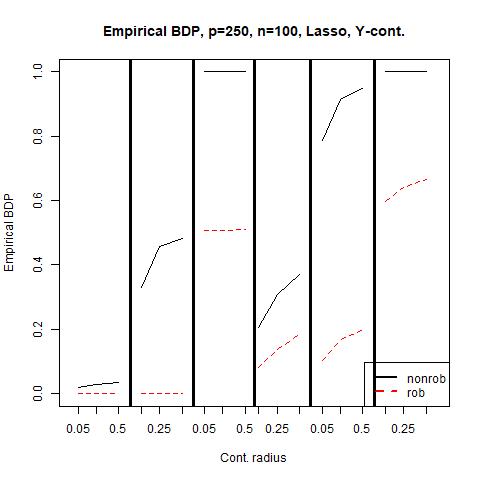}
\includegraphics[width=5cm,height=4cm]{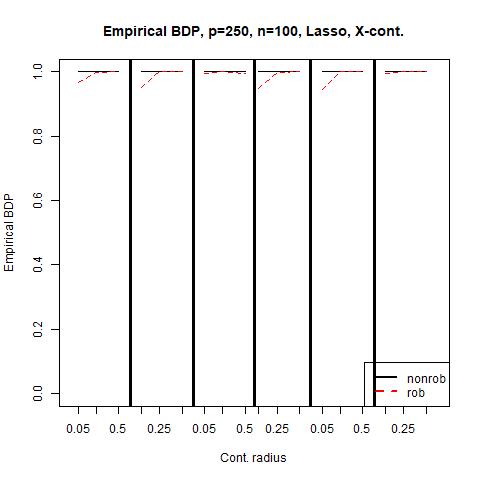} 
\caption{Mean empirical BDP of elicitability, i.e., mean relative number of repetitions where the weak ranking error corresponding to the identification of the contaminated test instances is non-zero.} \label{valcontp250bdp}
\end{center}
\end{figure}

\subsubsection{$p=500$, regression}

\begin{figure}[H]
\begin{center}
\includegraphics[width=5cm,height=4cm]{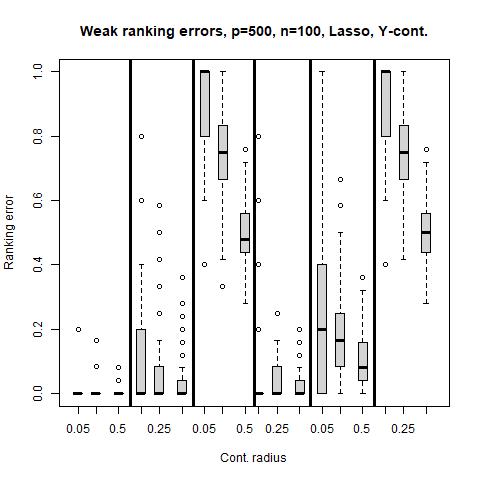}
\includegraphics[width=5cm,height=4cm]{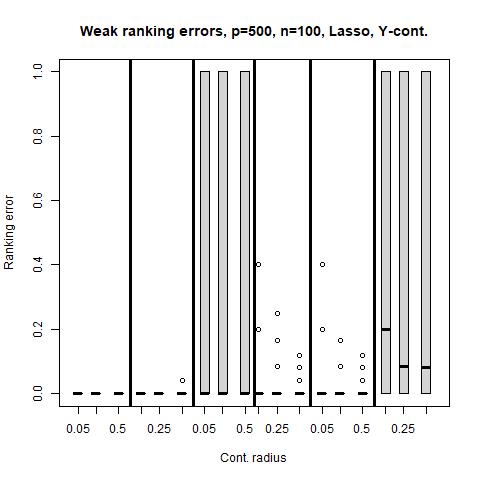} \\
\includegraphics[width=5cm,height=4cm]{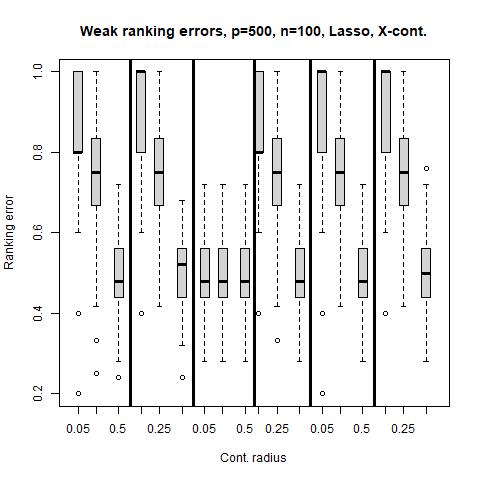} 
\includegraphics[width=5cm,height=4cm]{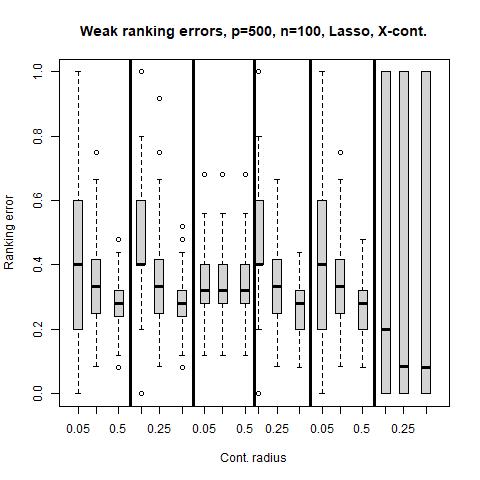} \\
\caption{Boxplots of the weak ranking errors corresponding to the model-based identification of the most outlying test instances. Upper row: $Y$-contamination, bottom row: $X$-contamination, left column: non-robust regression, right column: robust regression.} \label{valcontp500val}
\end{center}
\end{figure}

\begin{figure}[H]
\begin{center}
\includegraphics[width=5cm,height=4cm]{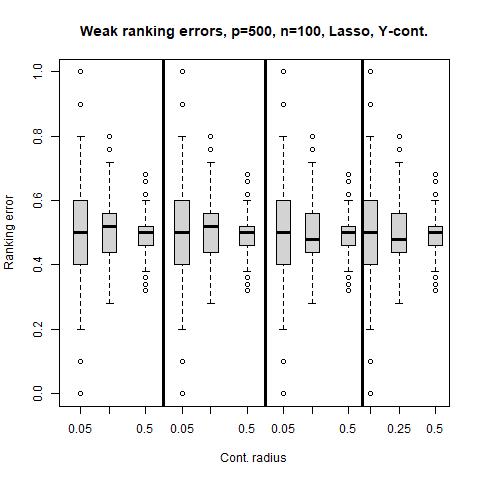}
\includegraphics[width=5cm,height=4cm]{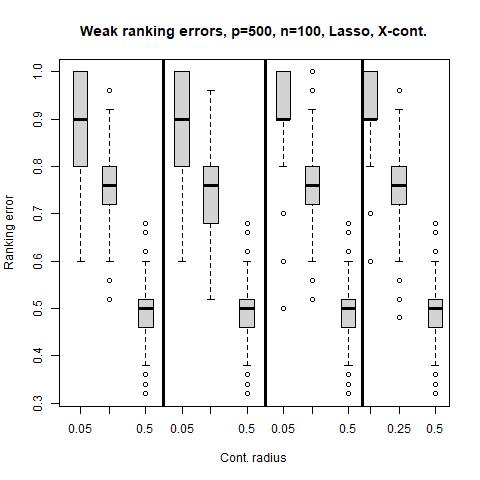} 
\caption{Boxplots of the weak ranking errors corresponding to the leave-one-out model-based identification of the most outlying training instances.  Left figure: $Y$-contamination, right figure: $X$-contamination. Left two columns: Non-robust regression, right two columns: robust regression.} \label{valcontp500tr}
\end{center}
\end{figure}

\begin{figure}[H]
\begin{center}
\includegraphics[width=5cm,height=4cm]{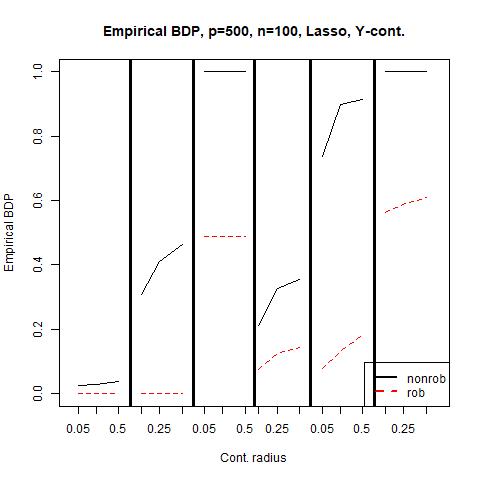}
\includegraphics[width=5cm,height=4cm]{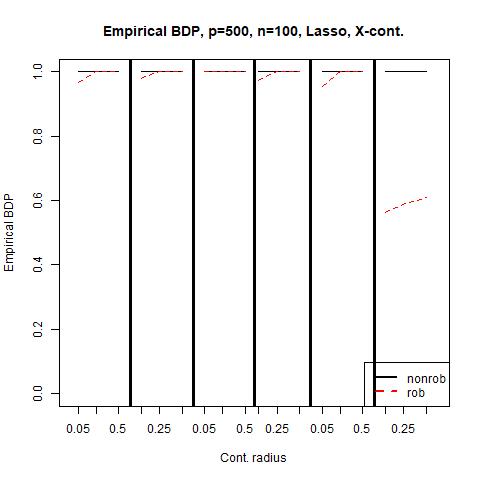} 
\caption{Mean empirical BDP of elicitability, i.e., mean relative number of repetitions where the weak ranking error corresponding to the identification of the contaminated test instances is non-zero.} \label{valcontp500bdp}
\end{center}
\end{figure}

\subsubsection{Classification}

\begin{figure}[H]
\begin{center}
\includegraphics[width=5cm,height=4cm]{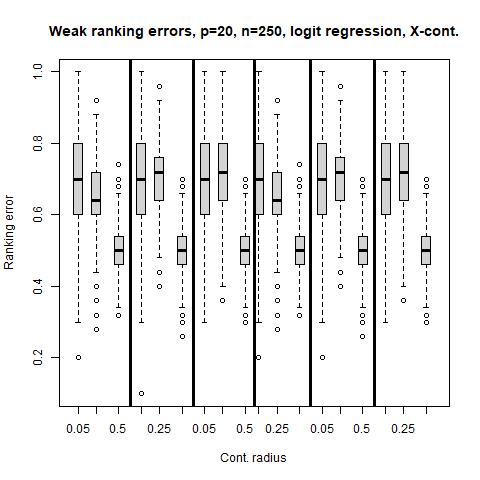} 
\includegraphics[width=5cm,height=4cm]{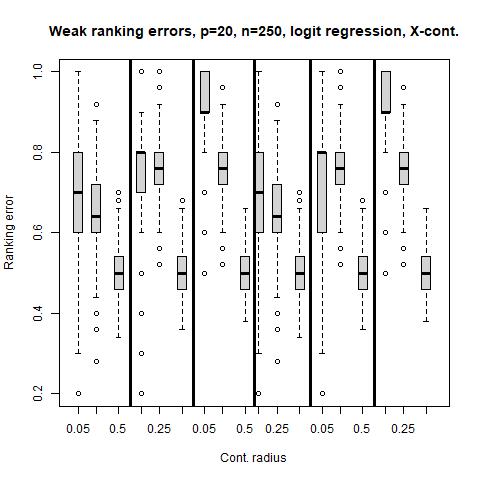} \\
\includegraphics[width=5cm,height=4cm]{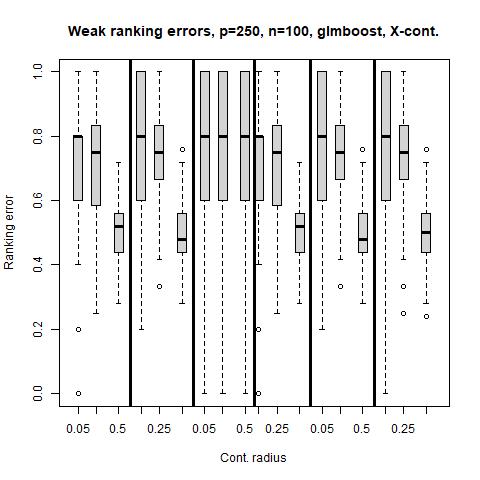} 
\includegraphics[width=5cm,height=4cm]{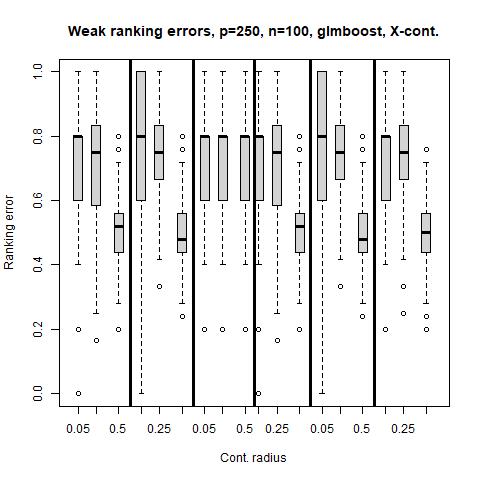} \\
\includegraphics[width=5cm,height=4cm]{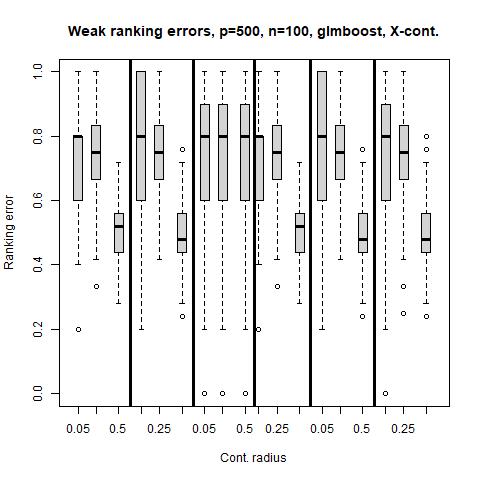} 
\includegraphics[width=5cm,height=4cm]{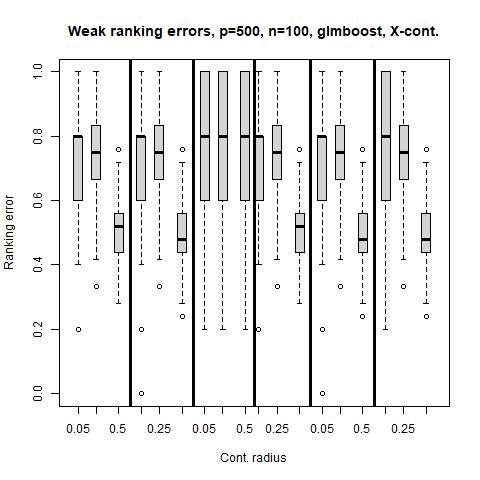} \\
\caption{Boxplots of the weak ranking errors corresponding to the model-based identification of the most outlying test instances. Left column: LogitBoost, right column: AUC-Boosting.} \label{valcontglmval}
\end{center}
\end{figure}

\begin{figure}[H]
\begin{center}
\includegraphics[width=5cm,height=4cm]{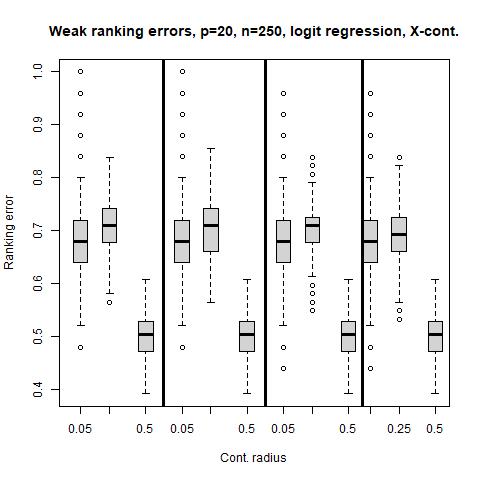} 
\includegraphics[width=5cm,height=4cm]{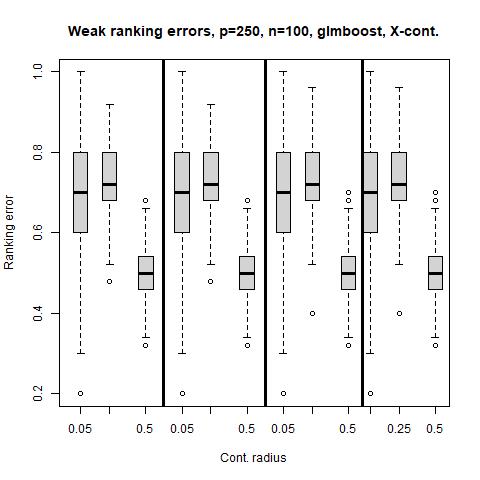} \\
\includegraphics[width=5cm,height=4cm]{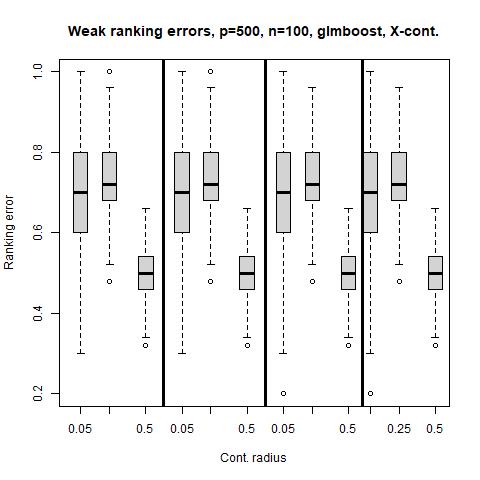} 
\caption{Boxplots of the weak ranking errors corresponding to the leave-one-out model-based identification of the most outlying training instances. Left two columns per figure: LogitBoost; right two columns: AUC-Boosting.} \label{valcontglmtr}
\end{center}
\end{figure}

\begin{figure}[H]
\begin{center}
\includegraphics[width=5cm,height=4cm]{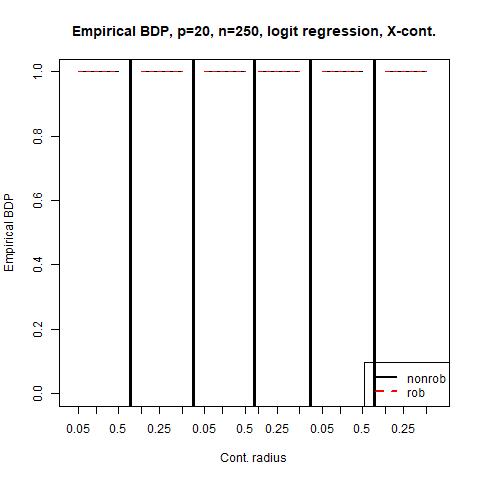} 
\includegraphics[width=5cm,height=4cm]{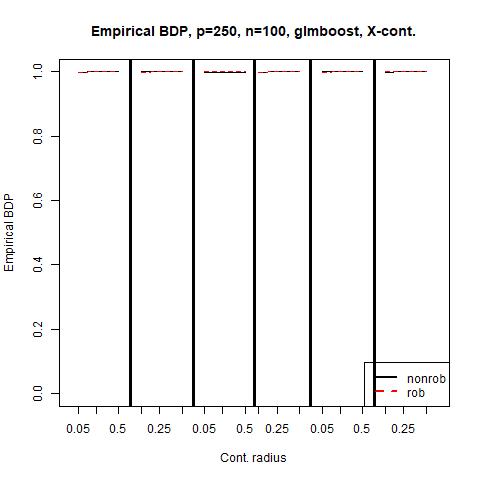} \\
\includegraphics[width=5cm,height=4cm]{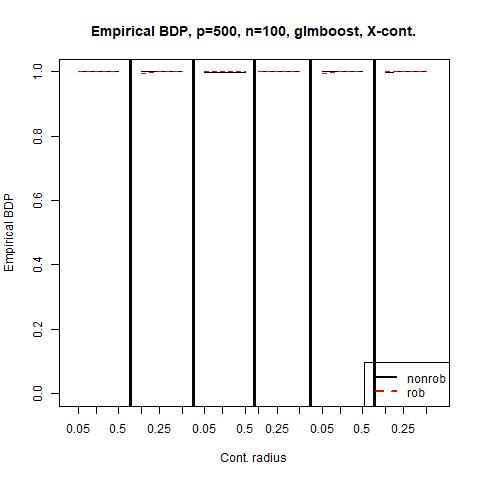} 
\caption{Mean empirical BDP of elicitability, i.e., mean relative number of repetitions where the weak ranking error corresponding to the identification of the contaminated test instances is non-zero.} \label{valcontglmbdp}
\end{center}
\end{figure}

\subsection{Contaminated training and test data, post-trimming} \label{app:valconttraintrim}

\subsubsection{$p=20$, regression}

\begin{figure}[H]
\begin{center}
\includegraphics[width=5cm,height=4cm]{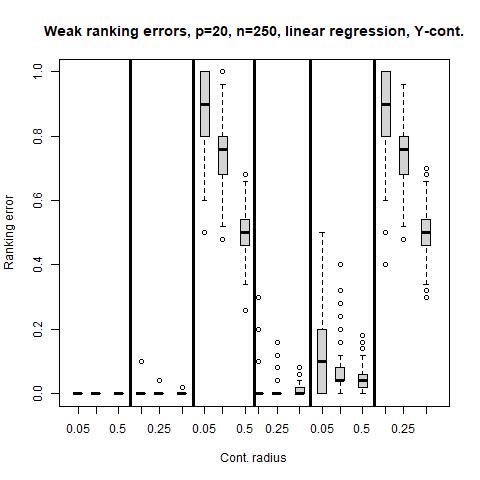}
\includegraphics[width=5cm,height=4cm]{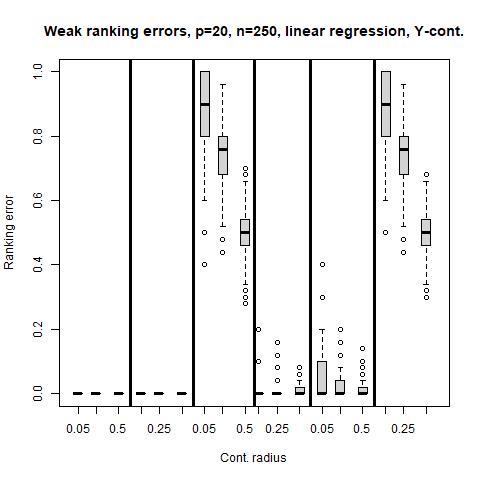} \\
\includegraphics[width=5cm,height=4cm]{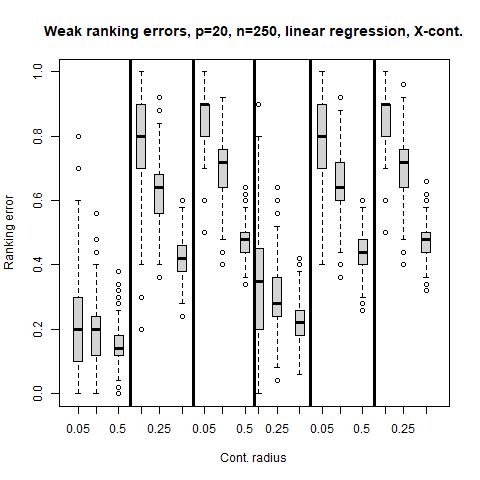} 
\includegraphics[width=5cm,height=4cm]{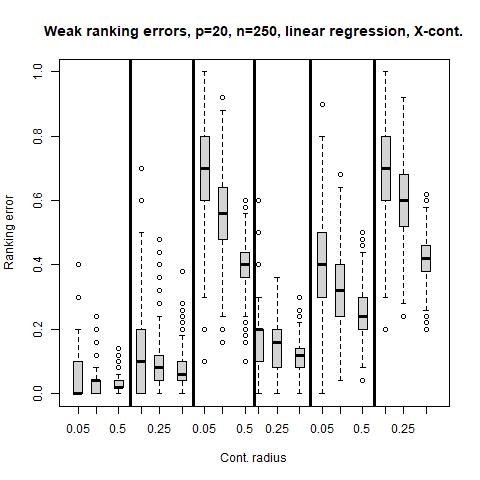} 
\caption{Boxplots of the weak ranking errors corresponding to the model-based identification of the most outlying test instances. Upper row: $Y$-contamination, bottom row: $X$-contamination, left column: non-robust regression, right column: robust regression.} \label{valconttraintrimp20val}
\end{center}
\end{figure}

\begin{figure}[H]
\begin{center}
\includegraphics[width=5cm,height=4cm]{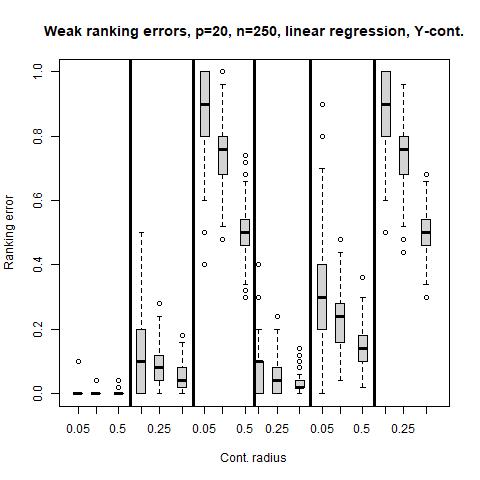}
\includegraphics[width=5cm,height=4cm]{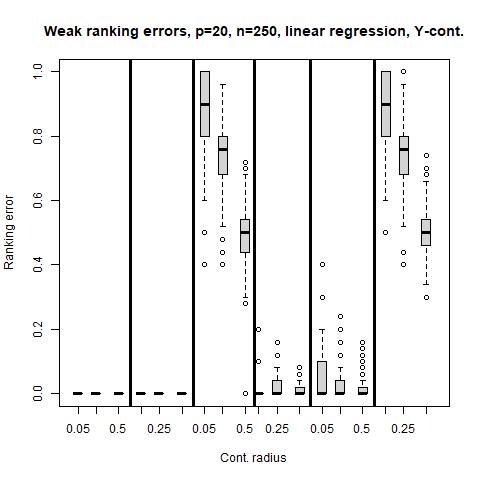} \\
\includegraphics[width=5cm,height=4cm]{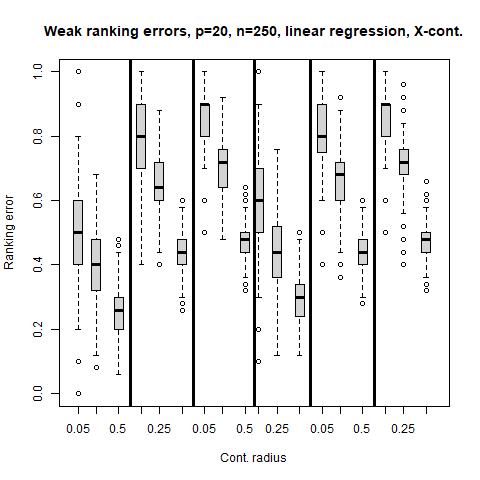} 
\includegraphics[width=5cm,height=4cm]{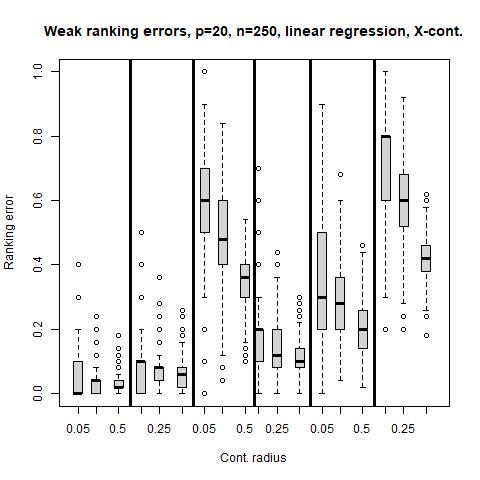} 
\caption{Boxplots of the weak ranking errors corresponding to the model-based identification of the most outlying test instances. Upper row: $Y$-contamination, bottom row: $X$-contamination, left column: non-robust regression, right column: robust regression.} \label{valconttraintrimp20valbetarob}
\end{center}
\end{figure}

\begin{figure}[H]
\begin{center}
\includegraphics[width=5cm,height=4cm]{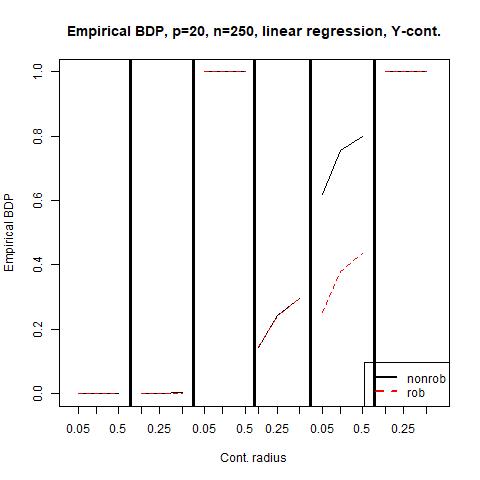}
\includegraphics[width=5cm,height=4cm]{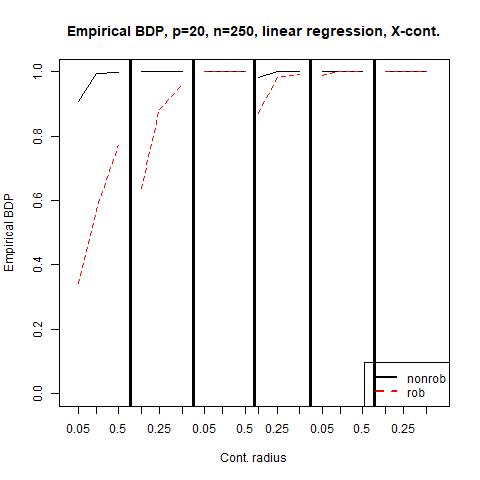} \\
\includegraphics[width=5cm,height=4cm]{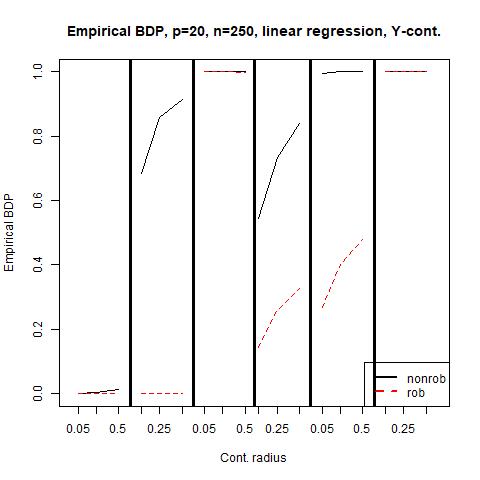}
\includegraphics[width=5cm,height=4cm]{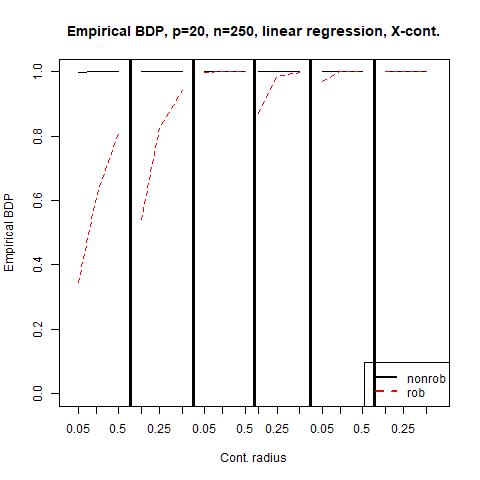} 
\caption{Mean empirical BDP of elicitability, i.e., mean relative number of repetitions where the weak ranking error corresponding to the identification of the contaminated test instances is non-zero.} \label{valconttraintrimp20bdp}
\end{center}
\end{figure}

\subsubsection{$p=250$, regression}

\begin{figure}[H]
\begin{center}
\includegraphics[width=5cm,height=4cm]{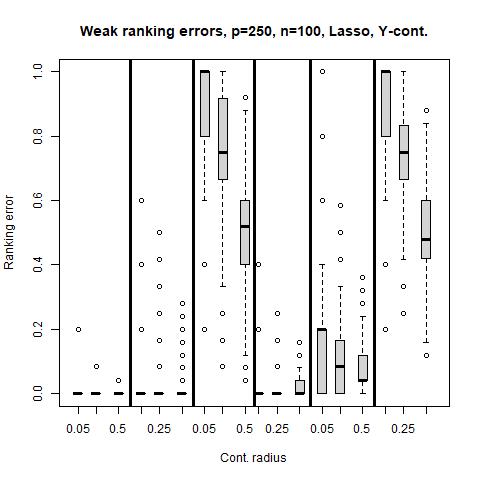}
\includegraphics[width=5cm,height=4cm]{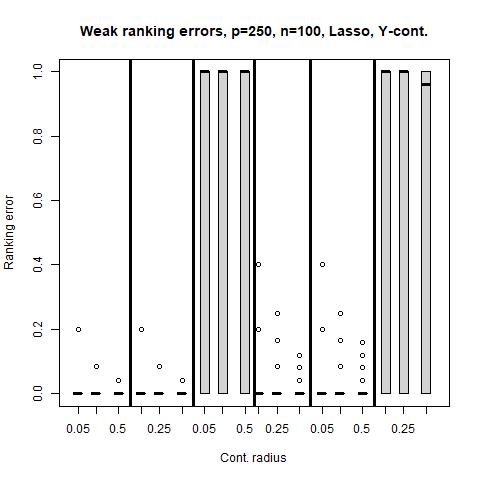} \\
\includegraphics[width=5cm,height=4cm]{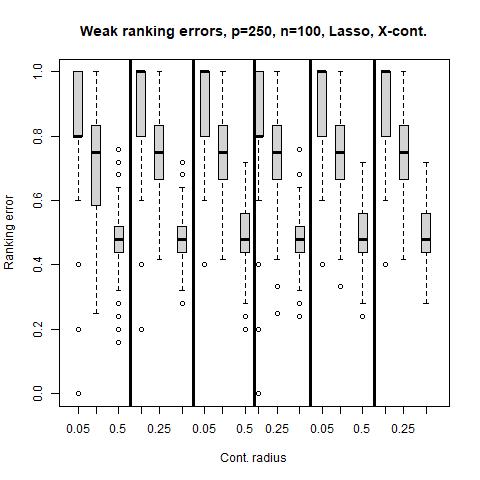} 
\includegraphics[width=5cm,height=4cm]{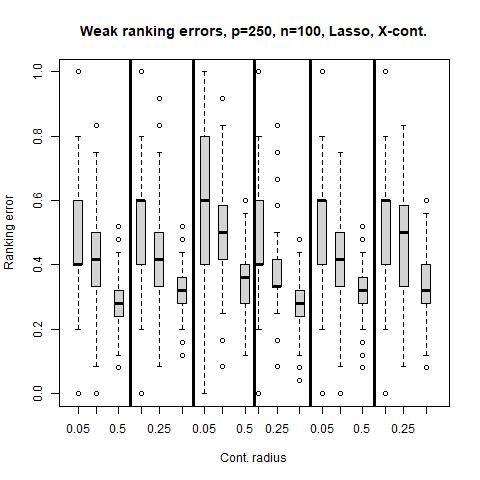} 
\caption{Boxplots of the weak ranking errors corresponding to the model-based identification of the most outlying test instances. Upper row: $Y$-contamination, bottom row: $X$-contamination, left column: non-robust regression, right column: robust regression.} \label{valconttraintrimp250val}
\end{center}
\end{figure}

\begin{figure}[H]
\begin{center}
\includegraphics[width=5cm,height=4cm]{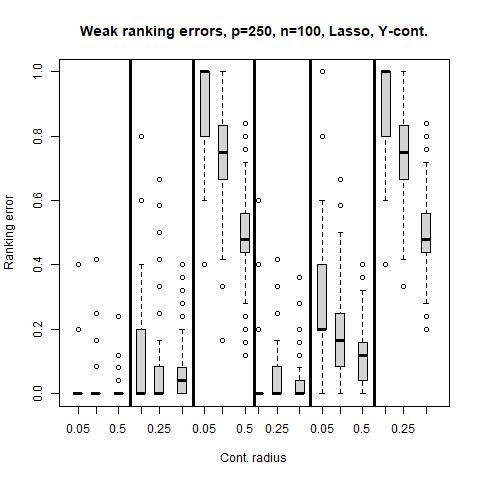}
\includegraphics[width=5cm,height=4cm]{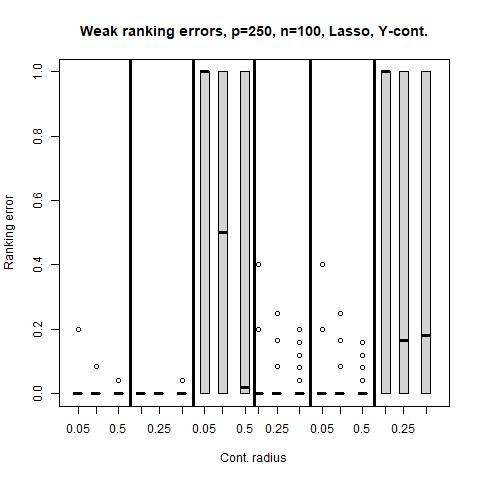} \\
\includegraphics[width=5cm,height=4cm]{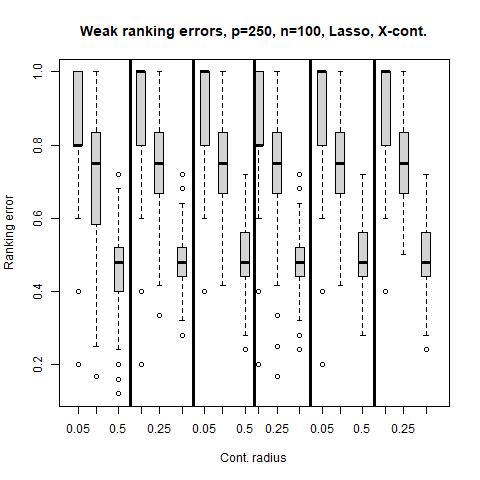} 
\includegraphics[width=5cm,height=4cm]{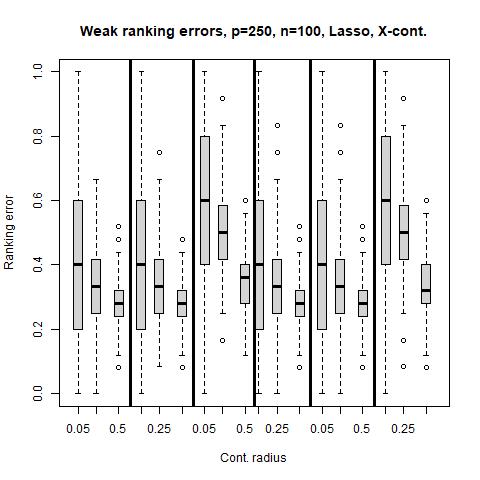} 
\caption{Boxplots of the weak ranking errors corresponding to the model-based identification of the most outlying test instances. Upper row: $Y$-contamination, bottom row: $X$-contamination, left column: non-robust regression, right column: robust regression.} \label{valconttraintrimp250valbetarob}
\end{center}
\end{figure}

\begin{figure}[H]
\begin{center}
\includegraphics[width=5cm,height=4cm]{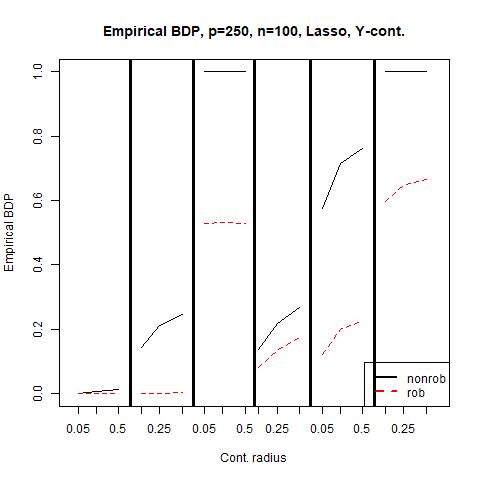}
\includegraphics[width=5cm,height=4cm]{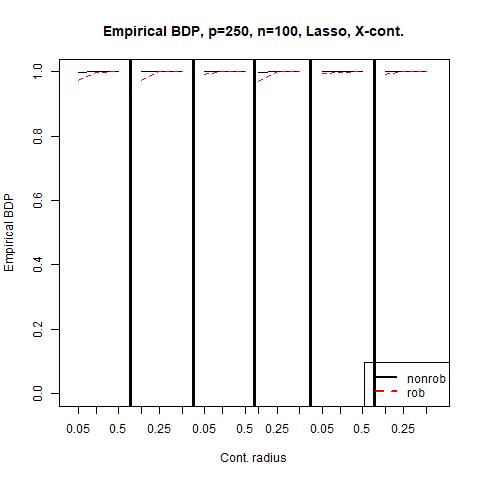} \\
\includegraphics[width=5cm,height=4cm]{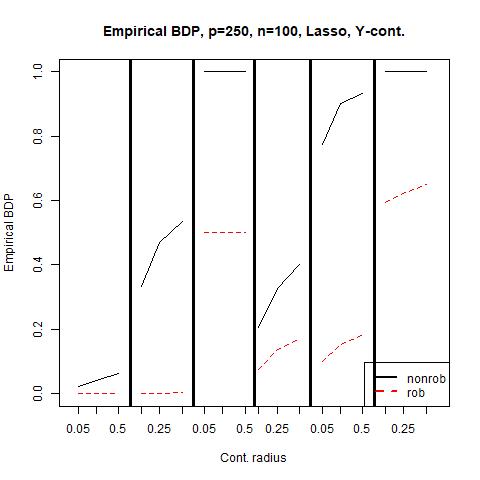}
\includegraphics[width=5cm,height=4cm]{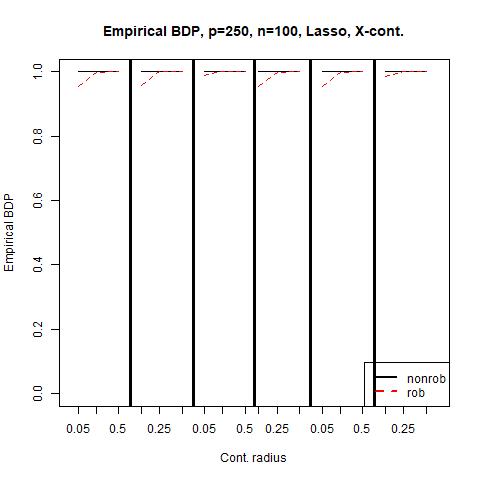} 
\caption{Mean empirical BDP of elicitability, i.e., mean relative number of repetitions where the weak ranking error corresponding to the identification of the contaminated test instances is non-zero.} \label{valconttraintrimp250bdp}
\end{center}
\end{figure}

\subsubsection{$p=500$, regression}

\begin{figure}[H]
\begin{center}
\includegraphics[width=5cm,height=4cm]{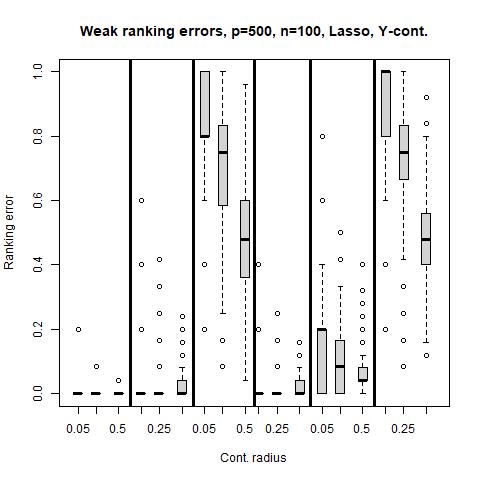}
\includegraphics[width=5cm,height=4cm]{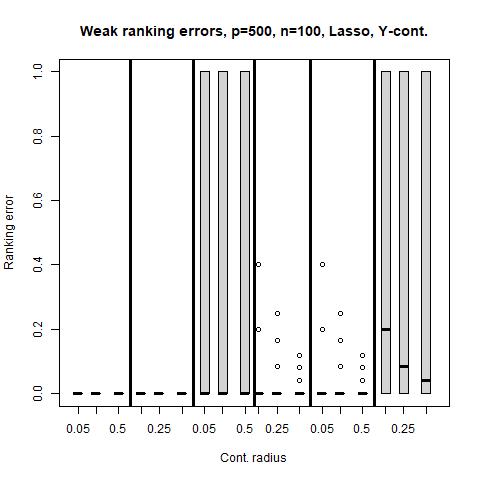} \\
\includegraphics[width=5cm,height=4cm]{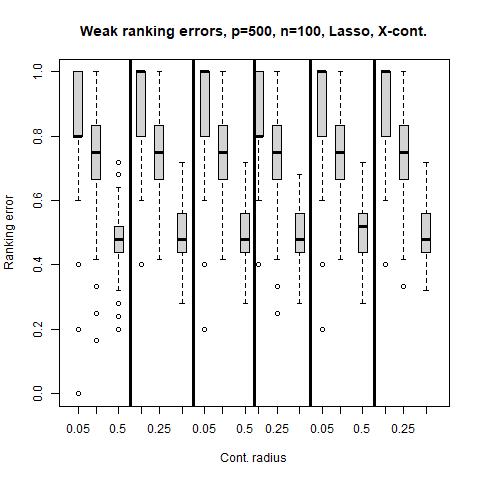} 
\includegraphics[width=5cm,height=4cm]{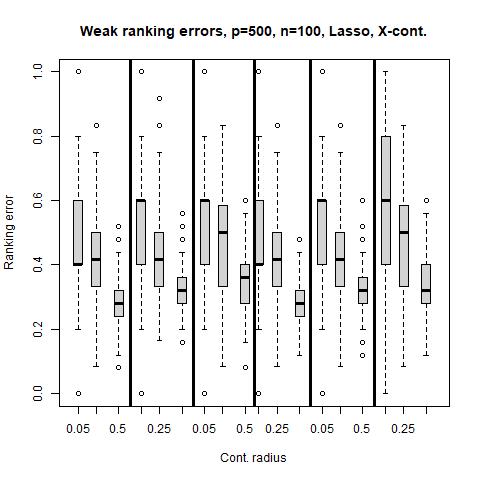} 
\caption{Boxplots of the weak ranking errors corresponding to the model-based identification of the most outlying test instances. Upper row: $Y$-contamination, bottom row: $X$-contamination, left column: non-robust regression, right column: robust regression.} \label{valconttraintrimp500val}
\end{center}
\end{figure}

\begin{figure}[H]
\begin{center}
\includegraphics[width=5cm,height=4cm]{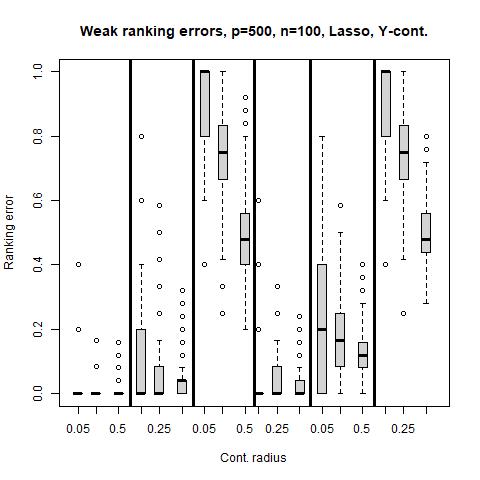}
\includegraphics[width=5cm,height=4cm]{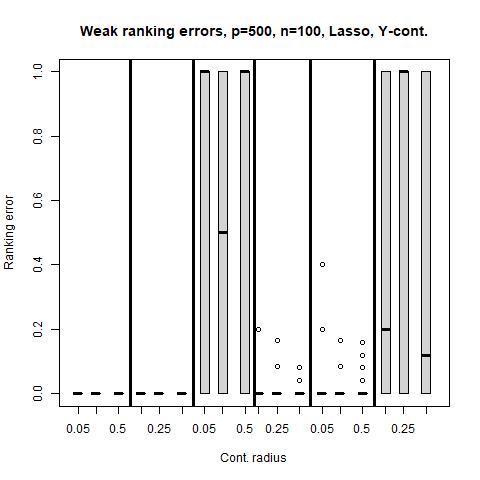} \\
\includegraphics[width=5cm,height=4cm]{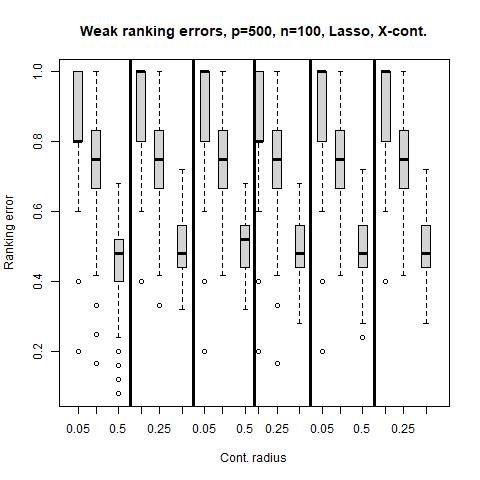} 
\includegraphics[width=5cm,height=4cm]{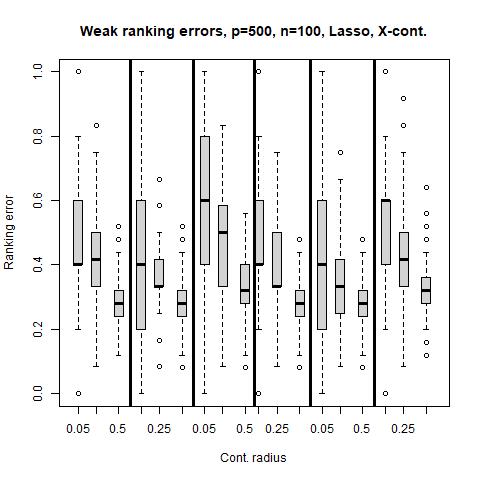} 
\caption{Boxplots of the weak ranking errors corresponding to the model-based identification of the most outlying test instances. Upper row: $Y$-contamination, bottom row: $X$-contamination, left column: non-robust regression, right column: robust regression.} \label{valconttraintrimp500valbetarob}
\end{center}
\end{figure}

\begin{figure}[H]
\begin{center}
\includegraphics[width=5cm,height=4cm]{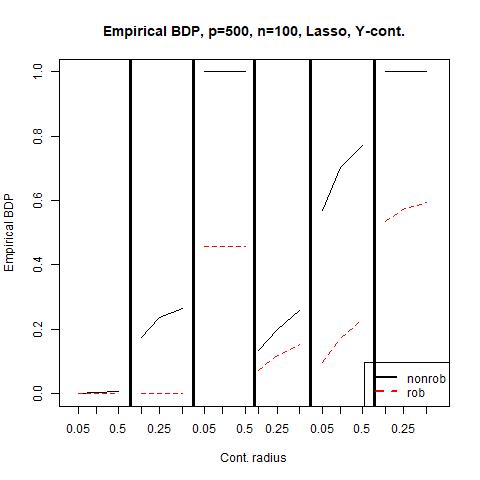}
\includegraphics[width=5cm,height=4cm]{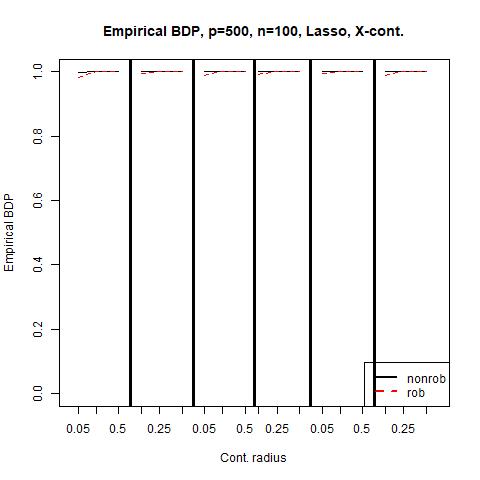} \\
\includegraphics[width=5cm,height=4cm]{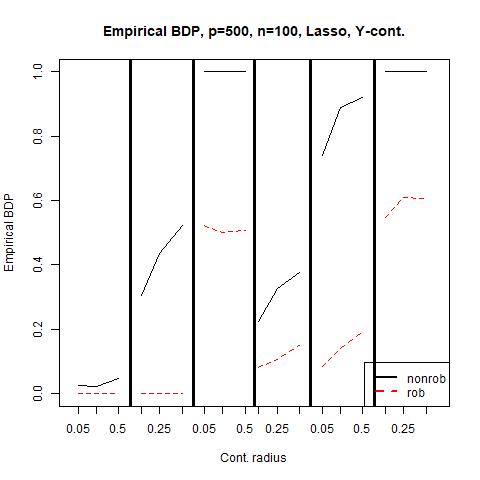}
\includegraphics[width=5cm,height=4cm]{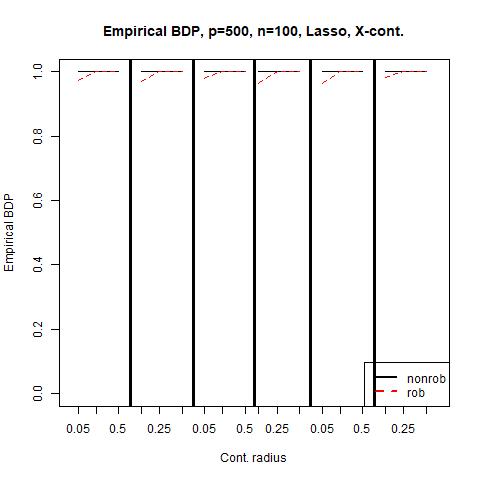} 
\caption{Mean empirical BDP of elicitability, i.e., mean relative number of repetitions where the weak ranking error corresponding to the identification of the contaminated test instances is non-zero.} \label{valconttraintrimp500bdp}
\end{center}
\end{figure}

\subsubsection{Classification}

\begin{figure}[H]
\begin{center}
\includegraphics[width=5cm,height=4cm]{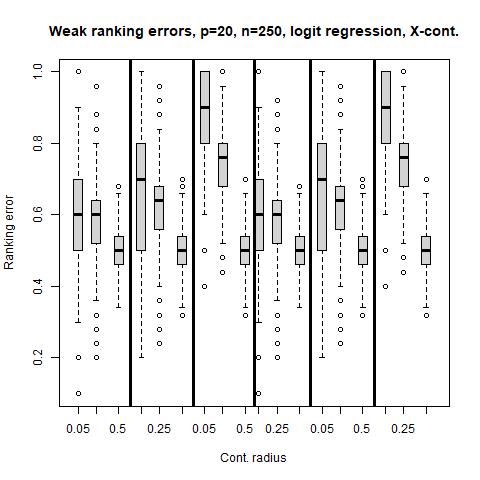} 
\includegraphics[width=5cm,height=4cm]{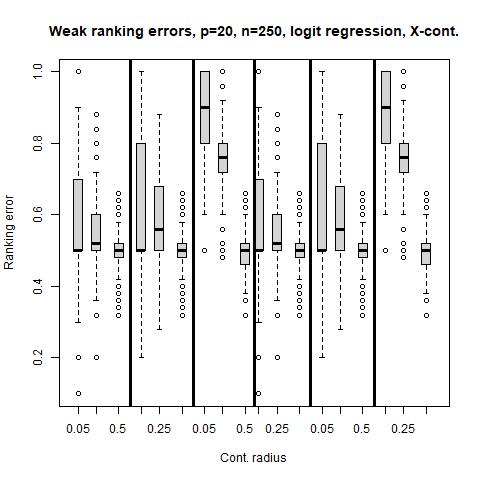} \\
\includegraphics[width=5cm,height=4cm]{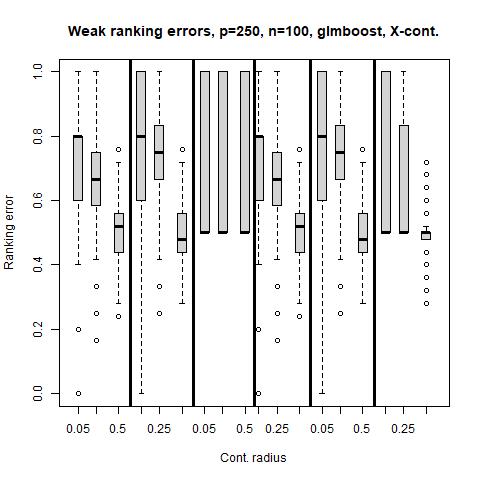} 
\includegraphics[width=5cm,height=4cm]{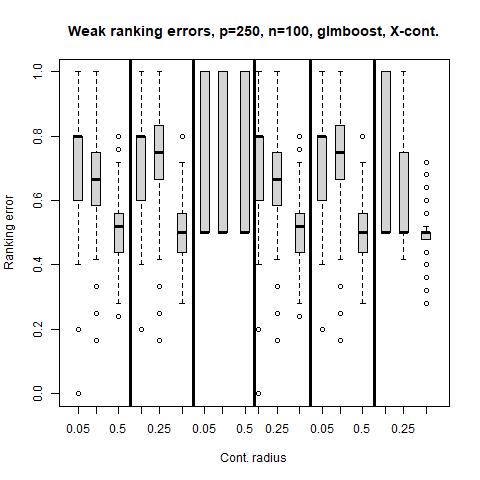} \\
\includegraphics[width=5cm,height=4cm]{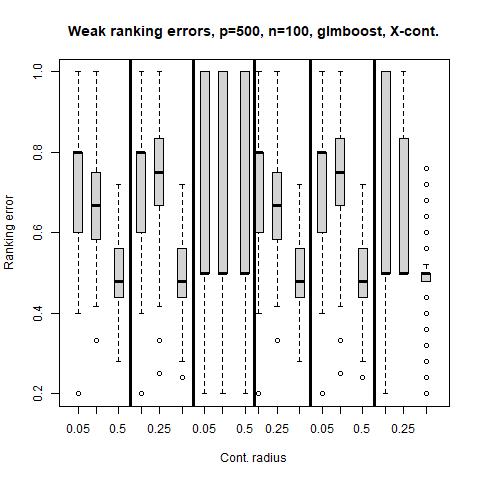} 
\includegraphics[width=5cm,height=4cm]{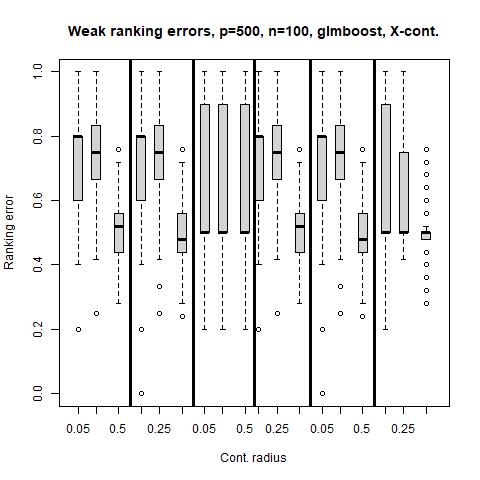} \\
\caption{Boxplots of the weak ranking errors corresponding to the model-based identification of the most outlying test instances. Left column: LogitBoost, right column: AUC-Boosting.} \label{valconttraintrimglmval}
\end{center}
\end{figure}

\begin{figure}[H]
\begin{center}
\includegraphics[width=5cm,height=4cm]{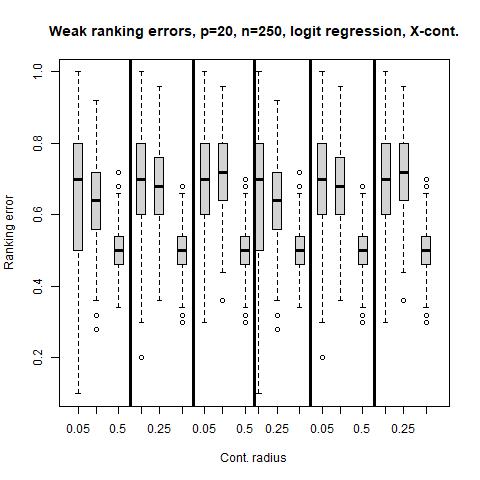} 
\includegraphics[width=5cm,height=4cm]{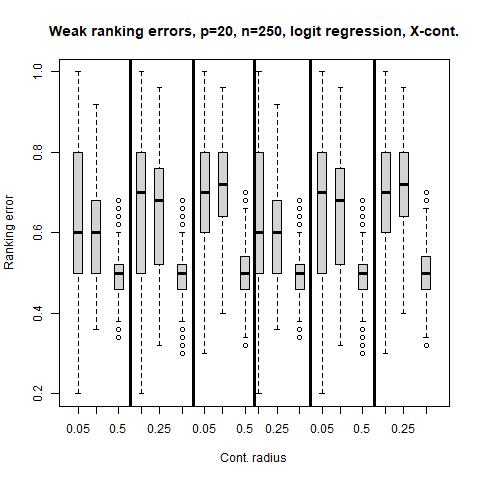} \\
\includegraphics[width=5cm,height=4cm]{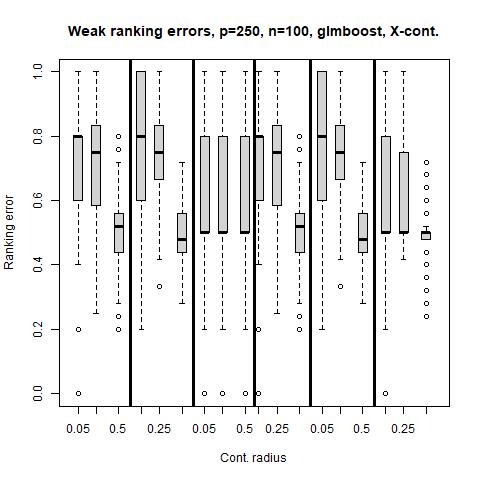} 
\includegraphics[width=5cm,height=4cm]{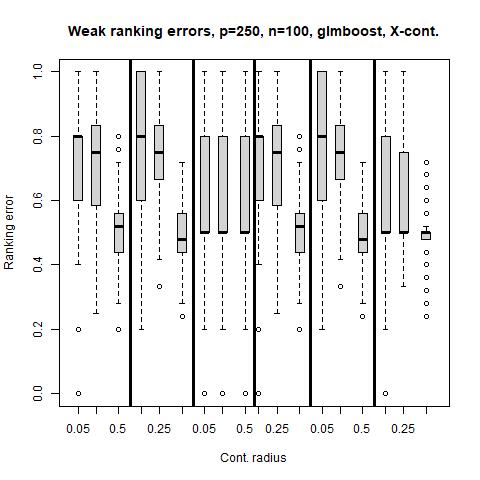} \\
\includegraphics[width=5cm,height=4cm]{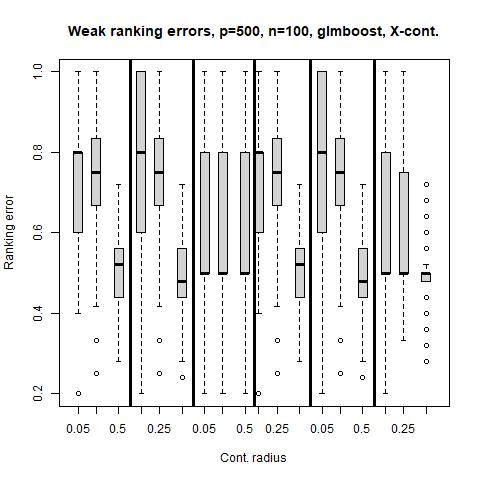} 
\includegraphics[width=5cm,height=4cm]{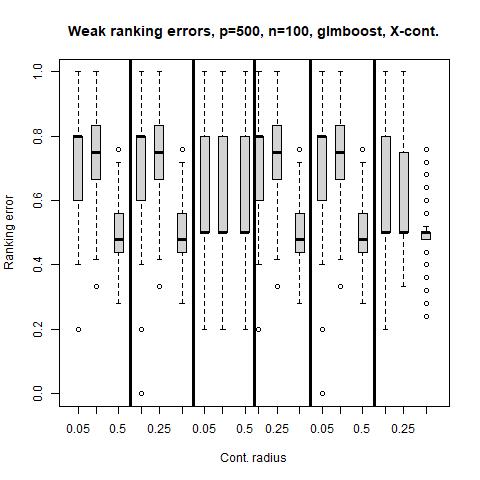} \\
\caption{Boxplots of the weak ranking errors corresponding to the model-based identification of the most outlying test instances. Left column: LogitBoost, right column: AUC-Boosting.} \label{valconttraintrimglmvalbetarob}
\end{center}
\end{figure}

\begin{figure}[H]
\begin{center}
\includegraphics[width=5cm,height=4cm]{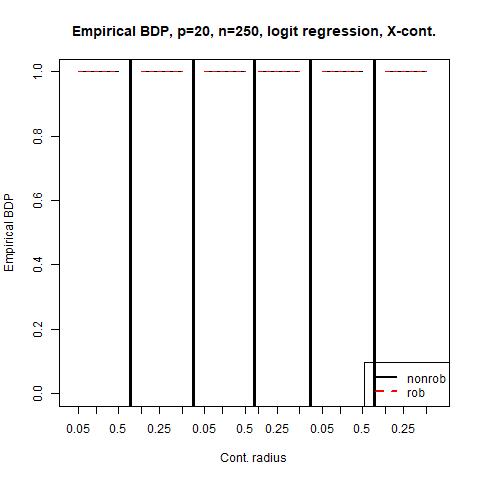} 
\includegraphics[width=5cm,height=4cm]{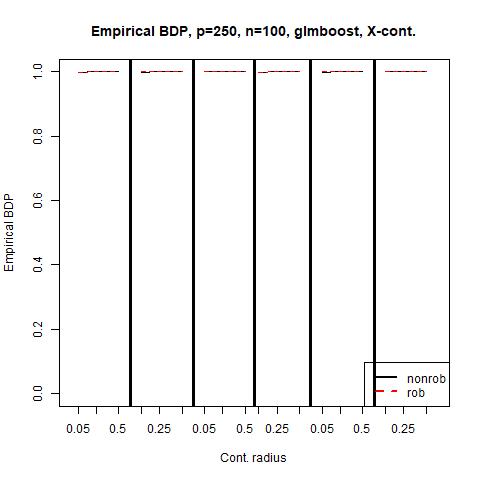} \\
\includegraphics[width=5cm,height=4cm]{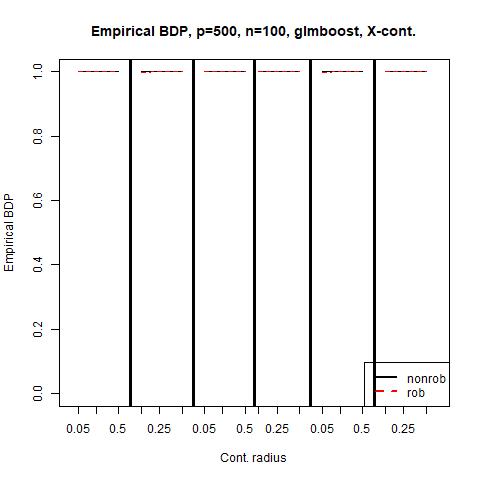} 
\includegraphics[width=5cm,height=4cm]{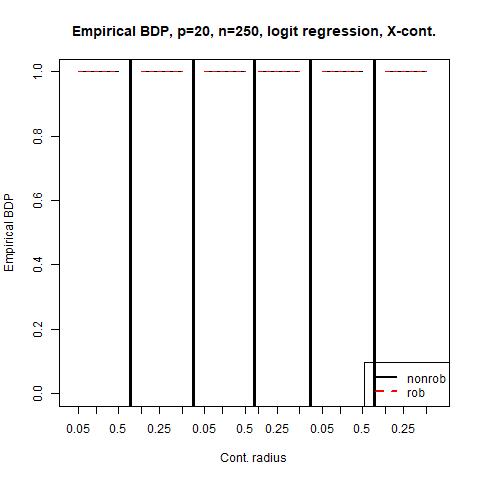} \\
\includegraphics[width=5cm,height=4cm]{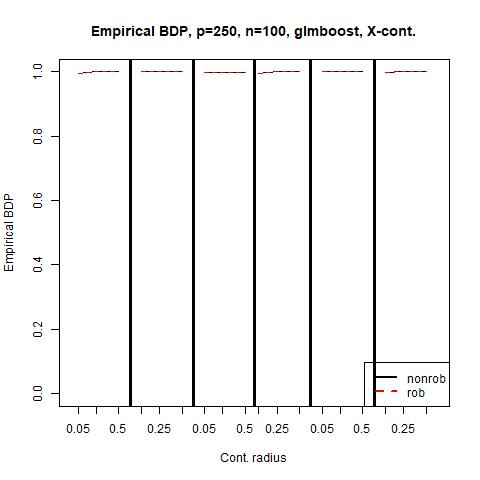} 
\includegraphics[width=5cm,height=4cm]{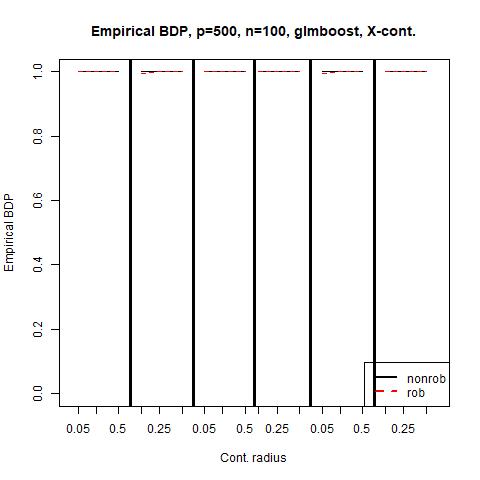} 
\caption{Mean empirical BDP of elicitability, i.e., mean relative number of repetitions where the weak ranking error corresponding to the identification of the contaminated test instances is non-zero.} \label{valconttraintrimglmbdp}
\end{center}
\end{figure}

\end{document}